\documentclass[twoside,11pt, preprint]{article}

\usepackage[abbrvbib]{jmlr2e}

\usepackage{lastpage}
\usepackage{amsmath}
\usepackage{array}
\usepackage{comment}
\usepackage{color}
\usepackage{srcltx}
\usepackage[mathscr]{eucal}
\usepackage[math]{easyeqn}
\usepackage{etoolbox}
\usepackage{hyperref}
\usepackage{nameref}
\usepackage{comment}
\usepackage{mathrsfs}
\usepackage{algorithm}
\usepackage[noend]{algpseudocode}
\usepackage{wasysym}
\usepackage{bbm}
\usepackage{nicefrac}
\usepackage{enumitem}
\usepackage{tikz}
\usetikzlibrary{positioning,shadows,arrows,trees,shapes,fit}
\usetikzlibrary{calc,intersections}

\def\cond{\, | \,}

\renewcommand\leq\leqslant
\renewcommand\geq\geqslant

\newcommand\eps\varepsilon

\newcommand{\E}{\mathbb E}
\newcommand{\1}{\mathbbm 1}
\newcommand{\p}{\mathbb P}
\newcommand{\R}{\mathbb R}
\newcommand{\X}{\mathbb X}
\newcommand{\Y}{\mathbb Y}

\newcommand{\A}{\mathcal A}
\newcommand{\B}{\mathcal B}
\newcommand{\C}{\mathcal C}
\newcommand{\Exp}{\mathcal E}
\newcommand{\F}{\mathcal F}

\newcommand{\J}{\mathcal J}
\newcommand{\K}{\mathcal K}

\newcommand{\M}{\mathcal M}
\newcommand{\N}{\mathcal N}

\newcommand{\Rcal}{\mathcal R}
\newcommand{\T}{\mathcal T}

\newcommand{\Z}{\mathcal Z}

\newcommand{\bid}{\boldsymbol I}

\newcommand{\bpi}{\boldsymbol\Pi}
\newcommand{\bsigma}{\boldsymbol\Sigma}
\newcommand{\bxi}{\boldsymbol\Xi}

\newcommand\ind[1]{^{(#1)}}

\newcommand{\mclass}{\mathscr M}
\newcommand{\pset}{\mathscr P}

\newcommand\reach[1]{\text{reach}\left(#1\right)}
\newcommand\proj[2]{\pi_{#1}\left(#2\right)}

\newcommand\kl{\mathcal K\mathcal L}

\newtheorem{Prop}{Proposition}
\newtheorem{Lem}{Lemma}
\newtheorem{Th}{Theorem}

\newtheorem{Rem}{Remark}

\ShortHeadings{Structure-adaptive Manifold Estimation}{Puchkin and Spokoiny}
\firstpageno{1}

\begin{document}

\title{Structure-Adaptive Manifold Estimation}

\author{\name Nikita Puchkin \email npuchkin@hse.ru \\
       \addr National Research University Higher School of Economics,\\
       Pokrovsky boulevard 11, 109028 Moscow, Russian Federation\\
       and\\
       \addr Institute for Information Transmission Problems RAS,\\
       Bolshoy Karetny per. 19, build.1, 127051 Moscow, Russian Federation
       \AND
       \name Vladimir Spokoiny \email spokoiny@wias-berlin.de \\
       \addr Weierstrass Institute and Humboldt University,\\
       Mohrenstrasse 39, 10117 Berlin, Germany\\
       and\\
       \addr National Research University Higher School of Economics,\\
       Pokrovsky boulevard 11, 109028 Moscow, Russian Federation\\
       and\\
       \addr Institute for Information Transmission Problems RAS,\\
       Bolshoy Karetny per. 19, build.1, 127051 Moscow, Russian Federation
}

\maketitle

\begin{abstract}%
	We consider a problem of manifold estimation from noisy observations. Many manifold learning procedures locally approximate a manifold by a weighted average over a small neighborhood. However, in the presence of large noise, the assigned weights become so corrupted that the averaged estimate shows very poor performance. We suggest a structure-adaptive procedure, which simultaneously reconstructs a smooth manifold and estimates projections of the point cloud onto this manifold. The proposed approach iteratively refines the weights on each step, using the structural information obtained at previous steps. After several iterations, we obtain nearly ``oracle'' weights, so that the final estimates are nearly efficient even in the presence of relatively large noise. In our theoretical study, we establish tight lower and upper bounds proving  asymptotic optimality of the method for manifold estimation under the Hausdorff loss, provided that the noise degrades to zero fast enough.
\end{abstract}


\section{Introduction}

We consider a problem of manifold learning, that is, to recover a low dimensional 
manifold from a cloud of points in a high dimensional space.
This problem is of great theoretical and practical interest.
For instance, if one deals with a problem of supervised or semi-supervised regression, 
the feature vectors, though lying in a very high-dimensional space, may occupy only a 
low-dimensional subset.
In this case, one can hope to obtain a rate of prediction which depends on the intrinsic 
dimension of the data rather than on the ambient one and escape the 
curse of dimensionality.
At the beginning of the century, the popularity of manifold learning gave rise to several 
novel nonlinear dimension reduction procedures, such as Isomap \citep{tsl00}, locally 
linear embedding \citep[LLE]{rs00} and its modification \citep{zw06}, Laplacian 
eigenmaps \citep{bn03}, and 
t-SNE \citep{vdmh08}.
More recent works include interpolation on manifolds via geometric multi-resolution 
analysis \citep{mms16}, local polynomial estimators \citep{al19} and numerical solution 
of PDE \citep{ss17}.
It is worth mentioning that all these works assume that the data points either lie exactly 
on the manifold or in its very small vicinity (which shrinks as the sample size \( n \) tends 
to infinity), so the noise \( \eps \) is so negligible that it may be ignored and put into a 
remainder term in Taylor's expansion.
However, in practice, this assumption can be too resrictive.
and the observed data do not exactly lie on a manifold.
One may think of this situation as there are unobserved \textquotedblleft 
true\textquotedblright \, features that lie exactly on the manifold and the learner 
observes its corrupted versions.
Such noise corruption leads to a dramatic decrease in the quality of manifold 
reconstruction for those algorithms which misspecify the model and assume that the 
data lies exactly on the manifold.
Therefore, one has to do a preliminary step, which is called manifold denoising (see e.g. 
\citep{hm06, wcp10, gsm10}), to first project the data onto the manifold.
Such methods usually act locally, i.e. consider a set of small neighborhoods, determined 
by a smoothing parameter (e.g. a number of neighbors or a radius \( h \)), and construct 
local approximations based on these neighborhoods.
The problem of this approach is that the size of the neighborhood must be large 
compared to the noise magnitude \( M \), which may lead to a non-optimal choice of the 
smoothing parameter.
The exclusion is the class of procedures, based on an optimization problem, such as 
mean-shift \citep{fh75, cheng95} and its variants \citep{wcp10, oe11, gppvw14}.
The mean-shift algorithm may be viewed as a generalized EM algorithm applied to 
the kernel density estimate (see \citep{cp07}).
This algorithm and its modifications were extensively studied in the literature 
\citep{cm02, hm06, li07, gppvw14, acmp16}.
For a comprehensive review on mean-shift algorithms, a reader is referred to \citep{cp15}. 
Though mean-shift algorithm was initially proposed for mode seeking and clustering, it 
found its applications in manifold denoising (see e.g. \citep{hm06, wcp10, oe11, gppvw14, 
cp15}).
If the observations lie around a smooth manifold, then few iterations of the mean-shift 
algorithm move the data towards the manifold.
However, since the mean shift algorithm and its variants (for example, 
subspace-constrained mean-shift \citep{oe11, gppvw14} which is based on density ridges \citep{eberly}) approximate the true density of $Y_1, \dots, Y_n$ by the kernel density estimate, they may suffer from the curse of 
dimensionality and the rates of convergence we found in the literature depend on the 
ambient dimension rather than on the intrinsic one in the noisy case.
To our best knowledge, only papers \citep{gppvw12a, gppvw12b} consider 
the case, when the noise magnitude does not tend to zero as \( n \) grows.
However, the approach in \citep{gppvw12a, gppvw12b} assumes that the noise 
distribution is known and has a very special structure.
For instance, considered in \citep{gppvw12a}, the noise has a uniform distribution in 
the direction orthogonal to the manifold tangent space.
Without belittling a significant impact of this paper, the assumption about the uniform 
distribution is unlikely to hold in practice.
Moreover, the authors point out that their goal was to establish minimax rates rather 
than propose a practical estimator.
Thus, there are two well studied extremal situations in manifold learning.
The first one corresponds to the case of totally unknown noise distribution but extremely 
small noise magnitude, and the other one corresponds to the case of large noise, which 
distribution is completely known.
This paper aims at studying the problem of manifold recovery under weak and realistic 
assumptions on the noise.

Below we focus on a model with additive noise.
Suppose we are given an i.i.d. sample \( \Y_{n} = (Y_{1}, \ldots, Y_{n}) \), where \( Y_{i} \) 
are independent copies of a random vector \( Y \) in \( \R^D \), generated from the model
\begin{equation}
	\label{model}
	Y = X + \eps.
\end{equation}
Here \( X \) is a random element whose distribution is supported on 
a low-dimensional manifold \( \M^{*} \subset \R^D \), \( \dim(\M^{*}) = d < D \), and \( \eps \) 
is a full dimensional noise. 
The goal of a statistician is to recover the corresponding unobserved variables \( \X_{n} 
= \{X_{1}, \dots, X_{n}\} \), which lie on the manifold \( \M^{*} \), and estimate \( \M^{*} \) 
itself.
Assumptions on the noise are crucial for the quality of estimation. 
One usually assumes that the noise is not too large, that is, \( \|\eps\| \leq M \) almost 
surely for some relatively small noise magnitude \( M \).
If the value \( M \) is smaller than the reach\footnote{A reader is referred to Section 
\ref{sec_model} for the definition.} of the manifold then the noise can be naturally 
decomposed in a component aligned with the manifold tangent space and another 
component describing the departure from the manifold.
It is clear that the impact of these two components is different, and it is natural to 
consider an anisotropic noise. 
For this purpose, we introduce a free parameter $b$ which controls the norm of the 
tangent component of the noise; see \eqref{a2} for the precise definition.
The pair of parameters \( (M, b) \) characterizes the noise structure more precisely than 
just the noise magnitude \( M \) and allows us to understand the influence of the noise 
anisotropy on the rates of convergence.
In our work we are particularly interested in situations when \( b \) is of order \( 1 \) 
(non-orthogonal noise) and when \( b \) is small (\emph{nearly} orthogonal noise) but our 
theoretical study is also valid for intermediate values of \( b \) (see Equation 
\ref{a3} below). 
We still have to assume that the noise magnitude \( M = M(n) \) tends to zero as \( n \) 
tends to infinity but aim at describing the best possible rate of convergence still 
ensuring a consistent estimation.
More precisely, if \( b \) is sufficiently small we allow \( M \) to be of order \( n^{-2/(3d+8)} 
\), which is much slower than, for instance, \( (\log n / n)^{2/d} \) and \( (\log n / n)^{1/d} \), 
considered in \citep{al19} and \citep{al18}, respectively (see the assumption \eqref{a3'} for 
the precise statement). 
To the best of our knowledge, this is the first paper which provides a rigorous theoretical 
study in this setup as well as the setup for intermediate \( b \).

As already mentioned, most of the existing manifold denoising procedures involve some 
nonparametric local smoothing methods
with a corresponding bandwidth.
The use of isotropic smoothing leads to the constraint that the noise magnitude is significantly smaller than the 
width of local neighborhoods; see e.g. \citep{hm06, mms16, ldmm,  al19}.
Similar problem arises even the case of effective dimension reduction in regression 
corresponding 
to the case of linear manifolds. 
The use of anisotropic smoothing helps to overcome this difficulty and to build efficient 
and asymptotically optimal estimation procedures; see e.g. \citep{Xia2002} or \citep{hjs01b}.
This paper extends the idea of \emph{structural adaptation} proposed in \citep{hjs01a, hjs01b}.
In these papers, the authors suggested to use anisotropic elliptic neighborhoods with axes shinking 
in the direction of the estimated effective dimension reduction (e.d.r.) subspace and 
stretching in the orthogonal directions to estimate the e.d.r. subspace.
As the shape of the local neighborhoods depends on the unknown e.d.r. structure, 
the procedure learns this structure from the data using iterations.
This explains the name ``structural adaptation''.
The use of anisotropic smoothing allows to obtain semiparametrically efficient 
root-n consistent estimates of the e.d.r. space \citep{Xia2002,hjs01b}.
In our method, we construct cylindric neighborhoods, which are stretched in a normal 
direction to the manifold.
However, our paper is not a formal generalization of \citep{hjs01a} and \citep{hjs01b}.
Those papers considered a regression setup, while our study focuses on a 
special unsupervised learning problem.
This requires to develop essentially different technique and use different mathematical 
tools for theoretical study
and substantially modify of the procedure.
Also to mention that a general manifold learning is much more involved than just linear dimension reduction, 
and a straightforward extension from the linear case is not possible.

Now we briefly describe our procedure.
Many manifold denoising procedures (see, for instance, \citep{hm06, gppvw14, ldmm, al18}) act in an iterative manner and our procedure is not an exception.
We start with some guesses \(\smash{\widehat{\bpi}\mathstrut\ind{0}_{1}, \dots, \widehat{\bpi}\mathstrut\ind 0_{n}} \) 
of the projectors onto the tangent spaces of \( \M^{*} \) at the points \( X_{1}, \dots, X_{n} \), respectively.
These guesses may be very poor, in fact.
Nevertheless, they give a bit of information, which can be used to construct initial 
estimates \(\smash{\widehat{X}\mathstrut_{1}\ind 0, \dots, \widehat{X}\mathstrut_{n}\ind 0} \).
On the other hand, the estimates \( \smash{\widehat{X}\mathstrut_{1}\ind 0, \dots, \widehat{X}\mathstrut_{n}\ind 0} \) 
help to construct the estimates \( \smash{\widehat{\bpi}\mathstrut_{1}\ind 1, \dots, \widehat{\bpi}\mathstrut_{n}\ind 1} \) of the 
projectors onto the tangent spaces of \( \M^{*} \) at the points \( X_{1}, \dots, X_{n} \), respectively, which are better than \( \widehat{\bpi}\mathstrut_{1}\ind 0, \dots, \widehat{\bpi}\mathstrut_{n}\ind 0 \).
One can repeat these two steps to iteratively refine the estimates of \( X_{1}, \dots, 
X_{n}\) and of the manifold \( \M^{*} \) itself.
We call this approach a \emph{structure-adaptive manifold estimation} (SAME).  
We show that SAME constructs such estimates \( \smash{\widehat{X}_{1}, \dots, \widehat{X}_{n}} \) of 
\( X_{1}, \dots, X_{n} \) and a manifold estimate \( \smash{\widehat{\M}} \) of \( \M^{*} \), such that
\begin{equation}
	\tag{Theorem \ref{th1}}
	\max\limits_{1\leq i \leq n} \|\widehat{X}_{i} - X_{i} \| \lesssim \frac{Mb \vee M h \vee 
	h^2}\varkappa + \sqrt{\frac{D(h^2 \vee M^2) \log n}{nh^d}},
\end{equation}
\begin{equation}
	\tag{Theorem \ref{th2}}
	d_H(\widehat{\M}, \M^{*}) \lesssim \left( \frac{M^2 b^2}{\varkappa^3} \vee 
	\frac{h^2}\varkappa \right) + \sqrt{\frac{D (h^{4}/\varkappa^2 \vee M^{2}) \log 
	n}{nh^d}},
\end{equation}
provided that \( h \gtrsim \left((D\log n/n)^{1/d} \vee (DM^{2}\varkappa^2\log 
n/n)^{1/(d+4)}\right) \) and \( M \) and, possibly, $b$ degrade to zero fast enough, and 
both inequalities hold with an overwhelming probability.
Here \( h \) is the width of a cylindrical neighborhood, which we are able to 
control, $\varkappa$ is a lower bound for the reach of $\M^*$ (see Section 
\ref{sec_model} for the definition of reach). 
Moreover, our algorithm estimates projectors \( \bpi(X_1), \dots, \bpi(X_n) \) onto tangent 
spaces at \( X_1, \dots, X_n \).
It produces estimates \( \widehat{\bpi}_{1}, \dots, \widehat{\bpi}_{n} \), such that
\begin{equation}
	\tag{Theorem \ref{th1}}
	\max\limits_{1\leq i \leq n} \|\widehat{\bpi}_{i} - \bpi(X_i) \| \lesssim 
	\frac{h}\varkappa 
	+ h^{-1} \sqrt{\frac{D (h^{4}/\varkappa^2 \vee M^{2}) \log n}{nh^d}}
\end{equation}
with high probability.
Here, for any matrix $\boldsymbol{A}$, $\|\boldsymbol{A}\|$ denotes its spectral norm.
The notation \( f(n)\lesssim g(n) \) means that 
there exists a constant \( c > 0 \), which does not depend on \( n \), such that \( f(n) \leq 
c g(n) \).
\( d_H(\cdot, \cdot) \) denotes the Hausdorff distance and it is defined as follows:
\[
	d_H(\M_{1}, \M_{2}) = \inf\left\{ \eps > 0 : \M_{1} \subseteq \M_{2} \oplus \B(0, \eps), \, 
	\M_{2} \subseteq \M_{1} \oplus \B(0, \eps) \right\},
\]
where \( \oplus \) stands for the Minkowski sum and \( \B(0, r) \) is a Euclidean ball in \( 
\R^D \) of radius \( r \).

The optimal choice of \( h \) yields
\[
	\max\limits_{1\leq i \leq n} \|\widehat{X}_{i} - X_{i} \| \lesssim \frac{Mb}\varkappa \vee 
	\frac1\varkappa \left( \frac{D\varkappa^2 \log n}n\right)^{\frac2{d+2}} \vee 
	\frac{M}\varkappa \left(\frac{DM^2\varkappa^2 \log n}n \right)^{\frac1{d+4}},
\]
\[
	d_H(\widehat{\M}, \M^{*}) \lesssim \frac{M^2 b^2}{\varkappa^3} \vee \frac1\varkappa 
	\left(\frac{D\log n}n \right)^{\frac{2}{d}} \vee 
	\frac1\varkappa \left(\frac{DM^{2}\varkappa^2 \log n}n\right)^{\frac2{d+4}}
\]
and
\[
	\max\limits_{1\leq i \leq n} \|\widehat{\bpi}_{i} - \bpi(X_i) \| \lesssim \frac1\varkappa 
	\left(\frac{D\log n}n \right)^{\frac1d} \vee
	\frac1\varkappa \left(\frac{DM^{2}\varkappa^2 \log n}n\right)^{\frac1{d+4}}.
\]
Note that the optimal choice of \( h \) is much smaller than a possible value \( 
n^{-2/(3d+8)} \) of the noise magnitude \( M \).
As pointed out in \citep{gppvw12b}, the manifold estimation can be considered 
as a particular case of the error-in-variables regression problem.
Then the rate \( (M^2/n \log n)^{2/(d+4)} \) makes sense since it corresponds to an 
optimal accuracy of locally linear estimation with respect to \( \|\cdot\|_\infty \)-norm in a 
nonparametric regression problem (which is also \( (M^2/n \log n)^{2/(d+4)} \)).
Besides, we prove a lower bound
\begin{equation}
	\tag{Theorem \ref{lower_bound}}
	\inf\limits_{\widehat\M}\sup\limits_{\M^*} \E d_H(\widehat\M, \M^*) \gtrsim 
	\frac{M^2 b^2}{\varkappa^3} \vee \varkappa^{-1}\left(\frac{M^2 \varkappa^2 \log 
	n}n\right)^{\frac 2{d+4}}
\end{equation}
which has never appeared in the manifold learning literature.
Here \( \widehat\M \) is an arbitrary estimate of \( \M^* \) and \( \M^* \) fulfills some 
regularity conditions, which are precisely specified in Theorem \ref{lower_bound}. 
Theorem \ref{lower_bound}, together with Theorem 1 from \citep{kz15}, where the 
authors managed to obtain the lower bound \( \inf\limits_{\widehat\M}\sup\limits_{\M^*} 
\E d_H(\widehat\M, \M^*) \gtrsim (\log n/n)^{2/d} \), claims optimality of our method.

The rest of this paper is organized as follows.
In Section \ref{sec_model}, we formulate model assumptions and introduce notations.
In Section \ref{sec_algorithm}, we provide our algorithm for manifold denoising and then 
illustrate its performance in Section \ref{sec_numerical}.
Finally, in Section \ref{sec_theoretical}, we give a theoretical justification of the 
algorithm and discuss its optimality.
The proofs of the main results are collected in Section \ref{proofs}.
Many technical details are contained in Appendix.

\section{Model Assumptions}
\label{sec_model}

Let us remind that we consider the model \eqref{model},
where \( X \) belongs to  the manifold \( \M^* \) and the 
distribution of the error vector \( \eps \) will be described a bit later in this section.
First, we require regularity of the underlying manifold \( \M^{*} \).
We assume that it belongs to a class \( \mclass_\varkappa^{d} \) of twice differentiable, 
compact, connected manifolds without a boundary, contained in a ball 
\( \B(0, R) \), with a reach, bounded below by \( \varkappa \), and dimension \( d \):
\begin{align}
	\label{a1}
	\notag
	\M^{*}\in
	\mclass_\varkappa^{d}
	= \big\{
	\M \subset \R^{D} : \M\text{ is a compact, connected manifold}
	\\\tag{A1}
	\text{without a boundary}, \M\in\C^2, \M \subseteq \B(0, R),
	\\\notag
	\reach\M \geq \varkappa, \text{dim}(\M) = d < D \big\} \, .
\end{align}
The reach of a manifold \( \M \) is defined as a supremum of such \( r \) that any point in 
\( \M \oplus \B(0, r) \) has a unique (Euclidean) projection onto \( \M \).
Here \( \oplus \) stands for the Minkowski sum and \( \B(0, r) \) is a Euclidean ball in \( 
\R^{D} \) of radius \( r \).
One can also use the following equivalent definition of the reach (see \citep[Section 2.1]{gppvw12a}).
For a point \( x \in \M \), let \( \T_{x}\M \) stand for a tangent space of \( \M \) at \( 
x \), i. e. a linear space spanned by the derivative vectors of smooth curves on the 
manifold passing through \( x \), and define a 
fiber
\[
	F_{r}(x) = \left( \{x\} \oplus \big( \T_{x}\M \big)^\perp \right) \cap \B(x, r),
\]
where \( ( \T_{x}\M )^\perp \) is an orthogonal complement of \( \T_{x}\M \).
Then \( \reach\M \) is a supremum of such \( r > 0 \) that for any \( x, x' \in \M \), \( x \neq x' \), 
the sets \( F_{r}(x) \) and \( F_{r}(x') \) do not intersect: 
\[
	\reach\M = \sup \left\{ r > 0 : \forall \, x, x' \in \M, x \neq x', \, F_{r}(x) \cap F_{r}(x') = 
	\emptyset \right\} \, .
\]
The requirement that the reach is bounded away from zero prevents $\M^*$ from having 
a large curvature.
In fact, if the reach of $\M^*$ is at least $\varkappa$, then the curvature of any 
geodesic on $\M^*$ is bounded by $1/\varkappa$ (see \cite[Lemma 3]{gppvw12a}).

Second, the density \( p(x) \) of \( X \) (with respect to the \( d \)-dimensional Hausdorff 
measure on \( \M^* \)) meets the following condition:
\begin{align}
	\label{a1'}
	&
	\tag{A2}
	\exists \, p_1 \geq p_0 > 0 : \forall x \in \M^* \quad p_0 \leq p(x) \leq p_1,
	\\&\notag
	\exists L \geq 0 : \forall \, x, x' \in \M^* \quad |p(x) - p(x')| \leq \frac{L\|x - 
	x'\|}\varkappa.
\end{align}

Besides the aforementioned conditions on \( \M^{*} \) and \( X \), we require some 
properties of the noise \( \eps \).
We suppose that, given \( X \in \M^{*} \), the conditional distribution \( (\eps \cond X) \) 
fulfils the following assumption: there exist \( 0 \leq M < \varkappa \) and \( 0 \leq b \leq 
\varkappa \), such that
\begin{align}
	\label{a2}
	\tag{A3}
	&
	\E( \eps \cond X ) = 0, \, \|\eps\| \leq M < \varkappa, \,
	\\&\notag
	\|\bpi(X)\eps\| \leq \frac{M b}\varkappa \quad \text{\( \p(\cdot \cond 
	X) \)-almost surely},
\end{align}
where \( \bpi(X) \) is the projector onto the tangent space \( \T_X\M^* \).
The model with manifold \( \M^* \in \mclass_\varkappa^d \) and the bounded noise has 
been extensively studied in literature (see \citep{gppvw12a, mms16, al18, al19, tsay19}).
In \citep{fikln18}, the authors consider the Gaussian noise, which is unbounded, but they 
restrict themselves on the event \( \max\limits_{1\leq i \leq n} \|\eps_i\| \leq \varkappa \), 
which is essentially similar to the case of bounded noise.
In our work, we introduce an additional parameter \( b \in [0, \varkappa] \), which 
characterises maximal deviation in tangent direction.

The pair of parameters \( (M, b) \) determines the noise structure more precisely than just 
the noise magnitude \( M \).
If \( b=0 \), we deal with perpendicular noise, which was studied in \citep{gppvw12a, al19}.
The case \( b=\varkappa \) corresponds to the bounded noise, which is not constrained to 
be orthogonal.
Such model was considered, for instance, in \citep{al18}.
In our work, we provide upper bounds on accuracy of manifold estimation for all pairs 
\( (M, b) \) satisfying the following conditions:
\begin{equation}
	\label{a3}
	\tag{A4}
	\begin{cases}
		M \leq A n^{-\frac{2}{3d+8}},\\
		M^3 b^2 \leq \alpha \varkappa \left[ \left( \frac{D \log n}n \right)^{\frac4d} \vee 
		\left( \frac{D M^2 \varkappa^2 \log n}n \right)^{\frac4{d+4}} \right],
	\end{cases}
\end{equation}
where \( A \) and \( \alpha \) are some positive constants.
Among all the pairs \( (M, b) \), satisfying \eqref{a3}, we can highlight two cases.
The first one is the case of maximal admissible magnitude:
\begin{align}
	\label{a3'}
	\tag{A4.1}
	&
	M = M(n) \leq A n^{-\frac2{3d+8}},
	\\&\notag
	b = b(n) \leq \frac{\sqrt{\alpha\varkappa}}{A^{3/2}} \left[ \left( \frac{D \log n}n 
	\right)^{\frac1d} \vee \left( \frac{D M^2 
	\varkappa^2 \log n}n \right)^{\frac1{d+4}} \right].
\end{align}
The second one is the case of maximal admissible angle:
\begin{equation}
	\label{a3''}
	\tag{A4.2}
	b = \varkappa, \quad
	M = M(n) \leq \left(\frac{D^4 \alpha^{d+4}}{\varkappa^{d-4}}\right)^{\frac1{3d+4}} 
	n^{-\frac4{3d + 4}}.
\end{equation}

If \eqref{a3'} holds, we deal with \emph{almost} perpendicular noise.
Note that in this case the condition \eqref{a2} ensures that \( X \) is very close to the 
projection \( \pi_{\M^*}(Y) \) of \( Y \) onto \( \M^{*} \).
Here and further in this paper, for a closed set \( \M \) and a point \( x \), \( \pi_\M(x) \) stands 
for a Euclidean projection of \( x \) onto \( \M \).
Thus, estimating \( X_1, \dots, X_n \), we also estimate the projections of \( Y_1, \dots, Y_n 
\) onto \( \M^{*} \).
Also, we admit that the noise magnitude \( M \) may decrease as slow as \( n^{-2/(3d+8)} 
\).
We discuss this condition in details in Section \ref{sec_theoretical} after Theorem 
\ref{th1} and compare it with other papers to convince the reader that the assumption 
\( M \leq A n^{-2/(3d+8)} \) is mild.
In fact, to the best of our knowledge, only in \citep{gppvw12a, gppvw12b} the authors 
impose weaker assumptions on the noise magnitude.
At the first glance, the condition \eqref{a3'} looks very similar to the case of 
orthogonal noise \( b = 0 \).
However, our theoretical study reveals a surprising effect: the existing lower bounds for manifold 
estimation in the case of perpendicular noise are different from the rates we prove for 
the case of almost perpendicular noise satisfying \eqref{a3'}.
We provide the detailed discussion in Section \ref{sec_theoretical} below.

Finally, if \eqref{a3''} holds, the noise is not constrained to be orthogonal.
However, in this case, we must impose more restrictive condition on the noise 
magnitude than in \eqref{a3'}.
Nevertheless, under the condition \eqref{a3''}, we show that the result of \cite{al18}, 
Theorem 2.7, where the authors also consider bounded noise, can be improved if one 
additionally assumes that the log-density \( \log p(x) \) is Lipschitz.
A more detailed discussion is provided in Section \ref{sec_theoretical}.

\section{A Structure-adaptive Manifold Estimator (SAME)}
\label{sec_algorithm}

In this section we propose a novel manifold estimation procedure based on a 
nonparametric smoothing technique and structural adaptation idea.
One of the most popular methods in nonparametric estimation is weighted averaging:
\begin{equation}
	\label{nw}
	\widehat{X}_{i}\ind{loc} = \frac{ \sum\limits_{j=1}^{n} w_{ij}\ind{loc} Y_{j} }{ 
	\sum\limits_{j=1}^{n} w_{ij}\ind{loc} }, \quad 1 \leq i \leq n,
\end{equation}
and \( w_{ij}\ind{loc} \) are the localizing weights defined by
\[
	w_{ij}\ind{loc} = \K\left( \frac{\|Y_{i} - Y_{j}\|^2}{h^2} \right), \quad 1 \leq i, j \leq n,
\]
where \( \K(\cdot) \) is a smoothing kernel and the bandwidth \( h = h(n) \) is a tuning 
parameter.
In this paper, we consider the kernel \( \K(t) = e^{-t} \).
\begin{Rem}
Instead of $\K(t) = e^{-t}$, one can take any two times differentiable, monotonously 
decreasing on $\R_+$ function such that it and its first and second derivatives have 
either exponential decay or finite support.
We use $\K(t) = e^{-t}$ to avoid further complications of the proofs.
\end{Rem}
The estimate \eqref{nw} has an obvious limitation.
Consider a pair on indices \( (i, j) \) such that \( \|X_{i} - X_{j}\| < h \) and \( h = h(n) \) is of 
order \( (\log n/n)^{1/d} \), which is known to be the optimal choice in the presence of 
small noise (see \citep[Proposition 5.1]{al18} and \cite[Theorem 6]{al19}).
If the noise magnitude \( M \) is much larger than \( (\log n/n)^{1/d} \) (which is the case 
we also consider), then $M > h$ and the weights \( \smash{w_{ij}\ind{loc}} \) carry wrong 
information about the neighborhood of \( X_{i} \), 
i.e. \( \smash{w_{ij}\ind{loc}} \) can be very small even if the distance \( \|X_{i} - X_{j}\| \) is smaller 
than \( h \).
This leads to a large variance of the estimate \eqref{nw} when \( h \) is of order 
\( (\log n/n)^{1/d} \), and one has to increase the bandwidth \( h \), inevitably making 
the bias of the estimate larger.

The argument in the previous paragraph leads to the conclusion that the weights 
\( \smash{w_{ij}\ind{loc}} \) must be adjusted.
Let us fix any \( i \) from 1 to n.
\textquotedblleft Ideal\textquotedblright \, localizing weights \( w_{ij} \) are such that 
they take into account only those indices \( j \), for which the norm \( \|X_{i} - X_{j}\| \) 
does not exceed the bandwidth \( h \) too much.
Of course, we do not have access to compute the norms \( \|X_{i} - X_{j}\| \) for all pairs 
but assume for a second that the projector \( \bpi(X_{i}) \) onto the tangent space \( \T_{X_{i}}\M^* \) was known.
Then, instead of the weights \( \smash{w_{ij}\ind{loc}} \), one would rather use the ones of the form
\[
	w_{ij}(\bpi(X_{i})) = \K\left( \frac{\|\bpi(X_{i})(Y_{i} - Y_{j})\|^2}{h^2} \right), \quad 1 \leq j \leq 
	n,
\]
to remove a large orthogonal component of the noise.
The norm $\|\bpi(X_{i})(Y_{i} - Y_{j})\|$ turns out to be closer to $\|X_i - X_j\|$ than $\|Y_i 
- Y_j\|$, especially if the ambient dimension is large.
Thus, instead of the ball \( \{ Y \colon \| Y - Y_{i} \| \leq h \} \) around \( Y_{i} \), we 
consider a cylinder 
\( \{ Y \colon \|\bpi_{i}(Y_{i} - Y)\| \leq h\} \), where \( \bpi_{i} \) is a  projector, which is assumed to be 
close to \( \bpi(X_{i}) \).
One just has to ensure that the cylinder does not intersect \( \M^* \) several times.
For this purpose, we introduce the weights
\begin{equation}
\label{aw}
	w_{ij}(\bpi_{i}) = \K\left( \frac{\|\bpi_{i}(Y_{i} - Y_{j})\|^2}{h^2} \right) \1\left( \| Y_{i} - Y_{j} \| \leq 
	\tau \right), \quad 1 \leq  j \leq n,
\end{equation}
with a constant \( \tau < \varkappa \).

The adjusted weights \eqref{aw} require a ``good'' guess \( \bpi_{i} \) of the projector \( \bpi(X_{i}) \).
The question is how to find this guess.
We use the following strategy.
We start with poor estimates \( \smash{\widehat{\bpi}_1\mathstrut\ind{0}, \dots, \widehat{\bpi}\mathstrut_{n} \ind{0}} \) of 
\( \bpi(X_1), \dots, \bpi(X_n) \) and take a large bandwidth \( h_0 \).
Then we compute the weighted average estimates \( \smash{\widehat{X}\mathstrut_1\ind 1, \dots, \widehat 
X\mathstrut_n\ind 1} \) with the adjusted weights \eqref{aw} and the bandwidth \( h_0 \).
These estimates can be then used to construct estimates \( \smash{\widehat{\bpi}\mathstrut_{1}\ind{1}, 
\dots, \widehat{\bpi}\mathstrut_n\ind{1}} \) of \( \bpi(X_1), \dots, \bpi(X_n) \), which are better than \( 
\smash{\widehat{\bpi}\mathstrut_{1}\ind{0}, \dots, \widehat{\bpi}\mathstrut_n\ind{0}} \).
After that, we repeat the described steps with a bandwidth \( h_1 < h_0 \).
This leads us to an iterative procedure, which is given by Algorithm \ref{algorithm}.

\begin{algorithm}[t]
	\caption{Structure-adaptive manifold estimator (SAME)}
	\label{algorithm}
	\begin{algorithmic}[1]
		\State The sample of noisy observations $\Y_n = (Y_1, \dots, Y_n)$, the initial 
		guesses \( \widehat{\bpi}\mathstrut_{1}\ind{0}, \dots, \widehat{\bpi}\mathstrut_{n}\ind{0} \) of 
		\( \bpi(X_1), \dots, \bpi(X_n) \), the number of iterations \( K + 1 \), an initial 
		bandwidth \( h_0 \), the threshold \( \tau \) and constants \(a > 1\) and \(\gamma 
		> 0\) are given.
		\For{ \( k \) from \( 0 \) to \( K \)}
		\State Compute the weights \( w_{ij}\ind{k} \) according to the formula
		\[
			w_{ij}\ind{k} = \K \left( \frac{\| \widehat{\bpi}_{i}\ind{k} (Y_{i} - Y_{j}) \|^2}{h_k^2} 
			\right) \1 \left( \|Y_{i} - Y_{j} \| \leq \tau \right), \quad 1 \leq i, j \leq n \, .
		\]
		\State Compute the estimates
		\begin{equation}
			\label{weighted_avg}
			\widehat{X}_{i}\ind{k} 
			= 
			\sum\limits_{j=1}^n w_{ij}\ind{k} Y_{j} \Big/ \biggl( \sum\limits_{j=1}^n 
			w_{ij}\ind{k} \biggr), \quad 1 \leq i \leq n \, .
		\end{equation}
		\State If \( k < K \), for each \( i \) from 1 to n, define a set \( \J_{i}\ind{k} = \{ j :   
		\|\widehat{X}_{j}\ind{k} - \widehat{X}_{i}\ind{k} \| \leq \gamma h_k \} \) and 
		compute the matrices
		\[
			\widehat{\bsigma}_{i} \ind{k} = \sum\limits_{j\in\J_{i}\ind{k}} 
			(\widehat{X}_{j}\ind{k} - \widehat{X}_{i}\ind{k})(\widehat{X}_{j}\ind{k} - \widehat{X}_{i}\ind{k})^T, 
			\quad 1 \leq i \leq n \, .
		\]
		\State If \( k < K \), for each \( i \) from 1 to n, define \( \smash{\widehat{\bpi}\mathstrut_{i}\ind{k+1}} \) as a projector onto a linear span of eigenvectors of \( \smash{\widehat{\bsigma}\mathstrut_{i}\ind{k}} \), corresponding to the largest \( d \) eigenvalues.
		\State If \( k < K \), set \( h_{k+1} = a^{-1} h_k \).
		\EndFor
		\Return the estimates \( \widehat{X}_1 = \widehat{X}\mathstrut_1\ind K, \dots, 
		\widehat{X}_n = \widehat{X}\mathstrut_n\ind K \).
	\end{algorithmic}
\end{algorithm}
Let us discuss the role of the parameter $\gamma$ in Algorithm \ref{algorithm}.
After the computation of the estimates $\smash{\widehat X\mathstrut_1\ind k, \dots, \widehat X\mathstrut_n\ind k}$, 
our goal is to use them to update the projectors $\smash{\widehat\bpi\mathstrut_1\ind k, \dots, 
\widehat\bpi\mathstrut_n\ind k}$.
However, our theoretical analysis (see Lemma \ref{main_lem} in the proof of Theorem 
\ref{th1} below) reveals that the weighted average \eqref{weighted_avg} removes well an 
orthogonal component of the noise but causes a shift of $\smash{\widehat X\mathstrut_1\ind k, \dots, 
\widehat X\mathstrut_n\ind k}$ in tangent direction.
An illustration is given in Figure \ref{fig_shift}.
Even if the noise is nearly orthogonal, the tangent component of $\smash{\widehat X\mathstrut_i\ind k - 
X_i}$ is $O(h_k)$ while the orthogonal component is only $O(h_k^2 / \varkappa)$.
This means that if we take such $\smash{\widehat X\mathstrut_j\ind k}$'s that $X_j$ is close to $X_i$ we 
can get a good estimate of the projector $\smash{\widehat\bpi\mathstrut_i\ind k}$.
However, even if $X_i$ and $X_j$ are close, the distance between $\smash{\widehat X\mathstrut_i\ind k}$ 
and $\smash{\widehat X\mathstrut_j\ind k}$ can be as large as $O(h_K)$, because of the shift in the 
tangent direction.
We take it into account and introduce an auxiliary parameter $\gamma$ to construct a 
set of indices $\smash{\J\mathstrut_i\ind k}$ which includes only those $j$'s that $X_j$ is close to $X_i$.

\begin{figure}[t]
	\noindent\centering
	\begin{tikzpicture}
	[inner sep=0.5mm,
	dot/.style={circle,draw=black,fill=black,thick}]
	\draw (0,0) .. controls (1,-1) and (1, -5) .. (10, 1);
	\draw[dashed] (4.55, -1.8) -- (3.55, 0.2);
	\draw[dashed] (4.55, -1.8) -- (5.55, -3.8);
	\draw[thick] (4.55, -1.8) -- (4, 0);
	\node (m) at (9, 0.2) {};
	\node[below] at (m.south) {$\mathcal M^*$};
	\node (xi) at (4.6, -1.9) [dot] {};
	\node[below] at (xi.south) {$\boldsymbol X_{\boldsymbol i}$};
	\node (yi) at (4, 0) [dot] {};
	\node[right] at (yi.east) {$Y_i = X_i + \eps_i$};
	\filldraw[fill=blue!20!white, opacity=0.3] (5.15, -4) --  (3.15, 0) -- (3.95, 0.4) 
	-- (5.95, -3.6) -- cycle;
	\draw [<->] (3.1, 0.1) -- (3.9, 0.5);
	\node (cap1) at (3.4, 0.5) {};
	\node[left] at (cap1.north) {$2 Mb/\varkappa$};
	\draw [<->] (3.05, -0.05) -- (5.05, -4.05);
	\node (cap2) at (3.65, -1.25) {};
	\node[left] at (cap2.west) {$2 M$};
	\filldraw[fill=green!20!white, opacity=0.3] (2.3, -2.3) --  (2.9, -3.5) -- (6.9, -1.5) 
	-- (6.3, -0.3) -- cycle;
	\draw [<->] (2.2, -2.35) -- (2.8, -3.55);
	\node (cap3) at (2.5, -2.95) {};
	\node[left] at (cap3.west) {$O(h_k^2/\varkappa)$};
	\draw [<->] (2.95, -3.6) -- (6.95, -1.6);
	\node (cap4) at (6.4, -2.1) {};
	\node[below] at (cap4.south) {$O(h_k)$};
	\node (xhat) at (5.9, -1) [dot] {};
	\node[left] at (xhat.south) {$\widehat X_i\ind k$};
	\end{tikzpicture}
	\caption{An illustration of how the weighted average estimate \eqref{weighted_avg} 
	induces a shift in a tangent direction. Even if the noise is nearly orthogonal 
	and the observation lies in a thin blue cylinder around the point of interest, the weighted average estimate falls into the green rectangle stretched in the tangent directions with high probability. Nevertheless, the averaging successfully removes the orthogonal component of the noise.}
	\label{fig_shift}
\end{figure}

The computational complexity of Algorithm \ref{algorithm} is $O(n^2 D^2 K + n D^3 K)$.
This includes $O(n^2 D^2)$ operations to update the weights $w_{ij}\ind k$, $1 \leq i, j 
\leq n$, and the estimates $\widehat X_i\ind k$ and $\smash{\widehat \bsigma\mathstrut_i\ind k}$, $1 \leq i 
\leq n$, on each 
iteration and $O(n D^3)$ operations to update the projectors $\smash{\widehat\bpi\mathstrut_i\ind k}$, $1 
\leq i \leq n$, on each iteration.
SAME requires slightly more time than, for instance, the manifold blurring mean shift 
algorithm (\cite[MBMS]{wcp10}, see the pseudocode in Appendix \ref{sec_mbms} 
below).
The complexity of MBMS is $O(n^2 D + n (D+\mathsf{k}) (D \wedge \mathsf{k})^2)$ per iteration.
Here $\mathsf{k}$ is the number of neighbors used by MBMS to perform local PCA.
SAME needs more operations to update the weights $w_{ij}\ind k$, $1 \leq i, j 
\leq n$, because of multiplication of the projectors $\smash{\widehat\bpi\mathstrut_i\ind k}$, $1 \leq i \leq 
n$, by the vectors $(Y_j - Y_i)$, $1 \leq i, j \leq n$.
If the parameter $k$ in MBMS is greater than $D$, then SAME and MBMS require the 
same time to perform PCA-type procedures.

\section{Numerical Experiments}
\label{sec_numerical}

In this section, we carry out simulations to illustrate the performance of SAME.
For convenience, theoretical results  were obtained for manifolds without a boundary 
(which is a common assumption in the manifold learning literature) but we use some 
well-known surfaces with boundary in the experiments.
The source code of all the numerical experiments described in this section is 
available on GitHub (\href{https://github.com/npuchkin/SAME/tree/master}{link}).

\subsection{Manifold Denoising and Dimension Reduction}

\begin{figure}[t]
	\noindent
	\centering
	\includegraphics[width=\linewidth]{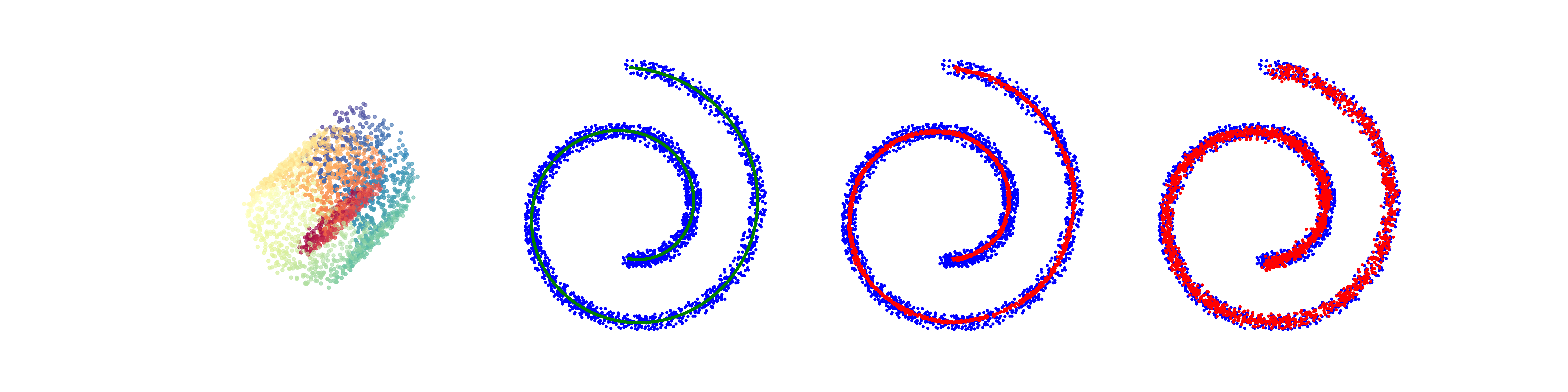}
	\vskip -16pt
	\includegraphics[width=\linewidth]{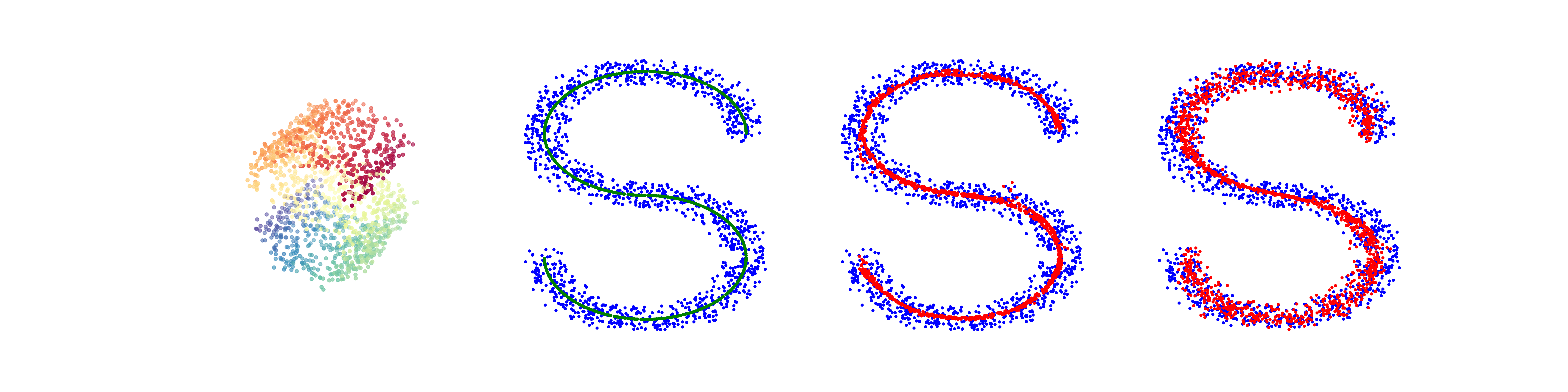}
	\vskip -16pt	
	\caption{Perfomance of SAME and MBMS on the Swiss Roll data set (top) and 
		on the S-shape data set (bottom). Column 1: noisy observations lying near a 
		two-dimensional manifold. Column 2: noisy observations (blue) and the true 
		manifold (green). Column 3: noisy observations (blue) and the projections onto 
		the 
		manifold (red), recovered by SAME. Column 4: noisy observations (blue) and 
		the 
		projections onto the manifold (red), recovered by MBMS.}
	\label{same_swiss_roll}
\end{figure}

In this section, we present the performance of SAME on two widely known artificial 
data sets: Swiss Roll and S-shape.
First, we show how our estimator denoises the manifold.
We start with the description of the experiment with the Swiss Roll.
We sampled \( n = 2500 \) points on a two-dimensional manifold in \( \R^3 \) and then 
embedded the surface into $\R^{20}$ adding $17$ dummy coordinates.
After that, we added a uniform noise with a magnitude \(0.75\) to each coordinate (thus, 
the noise magnitude $M$ was equal to $0.75 \cdot \sqrt{20}$).
In our algorithm, we initialized \( \smash{\widehat\bpi\mathstrut\ind 0_i = \bid_{20}} \) for all \( i \) from \( 1 \) 
to \( n \) and made $6$ iterations with \( h_k^2 = h_0^2 \cdot 1.25^{-k} \), \( 0 \leq k 
\leq 5 \), \( \tau = h_0 \), and \( \gamma = 4 \).
To choose the initial bandwidth $h_0$, we took $\alpha = 0.015$ and put $h_0$ equal to 
the distance to the $\lfloor\alpha n\rfloor$-th nearest neighbor of the first sample point.
The parameter $\gamma$ had minor influence on the behaviour of the algorithm and it 
was set to $4$ in all the experiments.
To quantify the performance of the algorithm, we used the mean squared error
\begin{equation}
	\label{mse}
	\frac1n \sum\limits_{i=1}^n \|\widehat X_i - X_i\|^2.
\end{equation}
The results are shown in Figure \ref{same_swiss_roll} (top) and in Table \ref{table_mse}.

\begin{table}[t]
	\noindent\centering
	\begin{tabular}{|c|c|c|}
		\hline
		Data set & MSE of SAME, $\times 10^2$ & MSE of MBMS, $\times 10^2$  \\
		\hline
		Swiss Roll & $\boldsymbol{67.4}$ & $69.3$ \\
		\hline
		S-shape & $\boldsymbol{3.9}$ & $4.7$ \\
		\hline
	\end{tabular}
	\caption{Mean squared errors (MSE, \eqref{mse}) of SAME and MBMS algorithms. 
		Best results 
		are boldfaced.}
	\label{table_mse}
\end{table}

We compare SAME with the manifold blurring mean shift algorithm \cite[MBMS]{wcp10}.
We provide a pseudocode of MBMS in Appendix \ref{sec_mbms} below to make the 
paper self-contained.
The parameters $\sigma$ and $k$ of MBMS as well as the number of iterations were 
chosen such that they minimized the mean squared error \eqref{mse} over a range of 
parameters.
The dimension $d$ was set to $2$.
The smallest MSE was achieved with $\sigma = 2.6 / \sqrt{2}$, $k = 150$ and only $1$ 
iteration.
The results are given in Figure \ref{same_swiss_roll} (top, right column).
We observed that the mean squared error of MBMS grew after $1$ or $2$ iterations.
In contrast to SAME, MBMS required much less iterations.
However, SAME recovered the surface better than MBMS.
The reason for that is hidden in the localizing weights $\smash{w_{ij}\ind k}$, $1 \leq i,j \leq n$, $1 
\leq k \leq K$.
Projection onto a tangent hyperplane removes a large orthogonal component of the noise 
and, hence, $\|\smash{\widehat\bpi\mathstrut_i\ind k (Y_j - Y_i)}\|$ better approximates the distance 
between $X_j$ and $X_i$ than $\|Y_j - Y_i\|$.
Besides, a weighted average estimate induces a shift in tangent direction (see the 
discussion after Algorithm \ref{algorithm} in Section \ref{sec_algorithm}) which may be 
harmful on early iterations.
As a consequence, $\smash{\|\widehat X\mathstrut_j\ind k - \widehat X\mathstrut_i\ind k\|}$ may also be a bad 
estimate for $\|X_j - X_i\|$.

The experiment with the S-shape manifold was carried in a similar way.
We took \( n = 1500 \) points on the manifold in $\R^3$, embedded the surface into 
$\R^{30}$, and added a uniform noise with a magnitude \(0.2 \) to each coordinate (thus, 
$M = 0.2 \cdot \sqrt{30}$).
Again, we initialized \( \smash{\widehat\bpi\mathstrut\ind 0_i = \bid_{30}} \) for all \( i \) from \( 1 \) to \( n \).
Then we put \( h_k^2 = h_0^2 \cdot 1.25^{-k} \), \( 0 \leq k \leq 13 \), i. e. the algorithm ran 
$14$ iterations.
The initial bandwidth $h_0$ was equal to the distance to the $\lfloor 0.1 n\rfloor$-th 
nearest neighbor of the first sample point. 
The parameters \(\tau\) and \(\gamma\) were equal to $h_0$ and $4$, respectively.
For MBMS, we took $\sigma = 0.45 / \sqrt{2}$, $k = 300$, $d=2$, and made $2$ 
iterations.
The tuning procedure of the parameters of MBMS was the same as in the example with 
the S-shape data set.
The result of this experiment is displayed in Figure \ref{same_swiss_roll} (bottom) and in 
Table \ref{table_mse}. 

\begin{figure}[t]
	\noindent
	\centering
	\includegraphics[width=\linewidth]{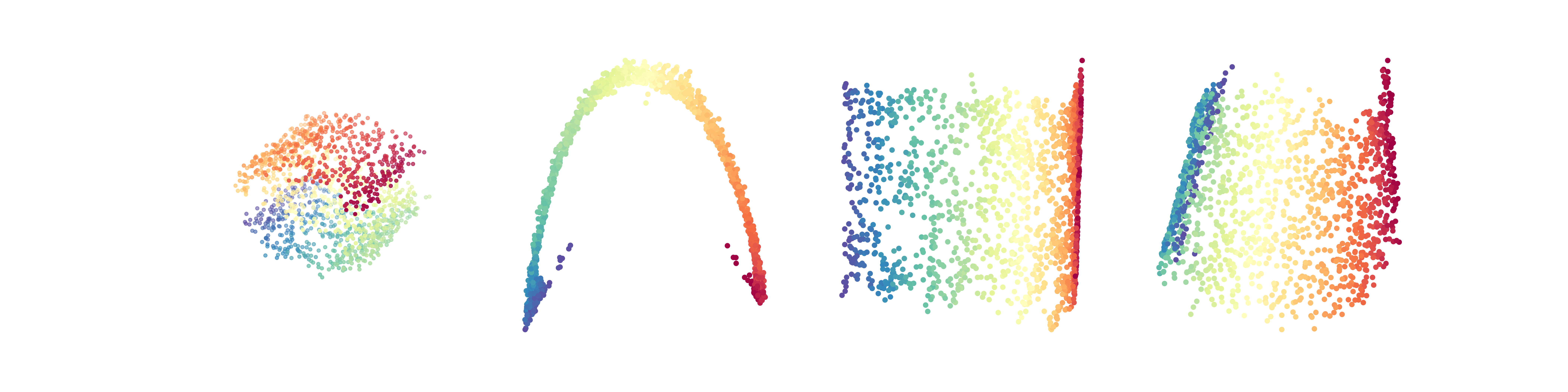}
	\vskip -26pt	
	\includegraphics[width=\linewidth]{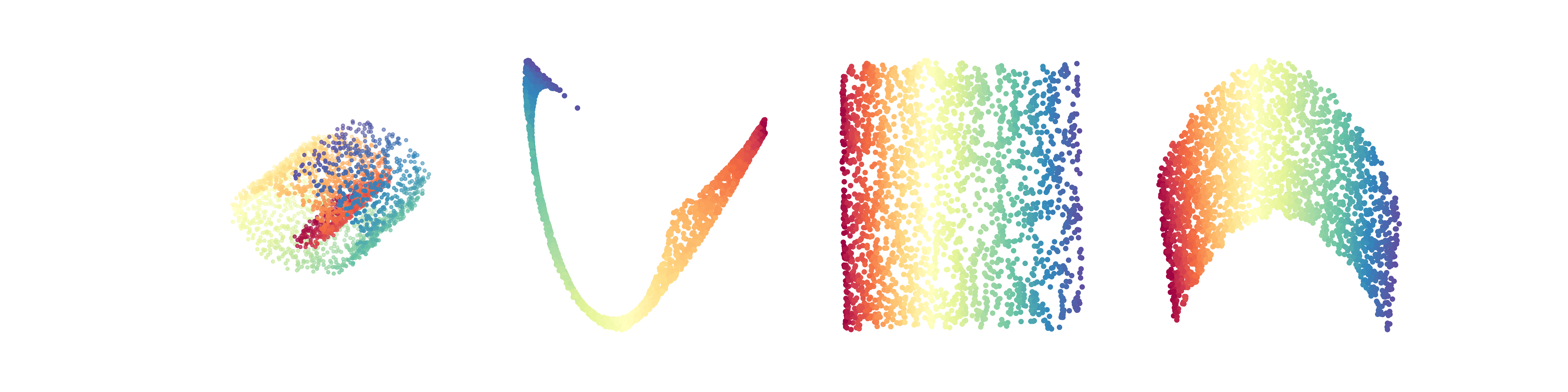}
	\vskip -16pt	
	\caption{The role of manifold denoising in a successful dimension reduction for the 
	S-shape data set (top) and Swiss Roll data set (Bottom). Column 1: noisy observations. Column 2:  application of MLLE to the data set without denoising. Column 3: application of MLLE to the data set with a preliminary denoising via SAME. Column 4: application of MLLE to the data set with a preliminary denoising via MBMS.}
	\label{dim_red_s-shape}
\end{figure}

Next, we show how the preliminary denoising step may improve a dimension reduction.
We consider the modified locally linear procedure \cite[MLLE for short]{zw06}, which is 
often used in applications due to its quality and computational efficiency.
MLLE takes high-dimensional vectors as an input and returns their low-dimensional 
representation.
In the case of S-shape and Swiss Roll data sets, one can easily find this map by 
straightening the curved surfaces into a plane.
In the noiseless case, MLLE solves this task.
However, as the other non-linear dimension reduction procedures based on Taylor's 
expansion, this algorithm deteriorates its performance in the presence of significant 
noise.
In Figure \ref{dim_red_s-shape} (center images) one can 
clearly observe that the MLLE procedure is not able to recognize a two-dimensional 
structure in the noisy data set.
Instead of a rectangular-like shape, which would be a natural choice to represent the 
two-dimensional structure of the S-shape and Swiss Roll data sets, we have a curve.
However, if one first uses SAME for manifold denoising and only after that applies MLLE 
for dimension reduction, then one obtains the desired result: both surfaces are 
straightened into planes.
Of course, popular dimension reduction methods (e.g. Isomap \citep{tsl00}, LLE 
\citep{rs00}, MLLE \citep{zw06}, Laplacian eigenmaps \citep{bn03}, t-SNE 
\citep{vdmh08}) still perform well in the presence of small noise.
However, a researcher should consider an option of using preliminary manifold denoising 
before dimension reduction in the case of larger noise.

\subsection{Manifold Denoising and Semi-supervised Learning}

Manifold denoising can be a preprocessing step in semi-supervised learning.
In the problem of semi-supervised learning, a statistician is usually given small amount of 
labelled data and a lot of unlabelled data.
The goal is to recover the labels of the unlabelled points (transductive semi-supervised 
learning) or to propose a rule for prediction of label of a test point \( x \) (``true'' semi-supervised learning).
In semi-supervised learning, it is usually assumed that the unlabelled data carries useful 
information, which may be useful for prediction.
The most popular assumptions is that the data has a cluster structure or data points lie 
in a vicinity of a low-dimensional manifold.

In this section, we pursue the goal of recovering labels of unlabelled points.
We take two artificial data sets g241c and g241n, which are described in \citep{csz10}.
The data sets g241c and g241n have \( n = 1500 \) pairs \( (Y_i, Z_i) \), \( 1 \leq i \leq n \), where 
\( Z_i \in \{-1, 1\} \) is a binary label and \( Y \in \R^D \), \( D = 241 \) is a high-dimensional 
feature vector.
According to \citep{csz10}, the data sets g241c and g241n were generated in such a way 
that they have a cluster structure and do not have hidden manifold structure.
However, in \citep{hm06}, the authors report that preliminary manifold denoising step, 
applied to these data sets, improves the quality of classification.
In this section, we illustrate that preliminary denoising with SAME also improves 
classification error.

We split the data sets into 100 train points and 1400 test points.
We use k-nearest neighbors classifier as a baseline for two reasons.
First, k-NN method is popular and often used in practice.
Second, k-NN classifier is based on pairwise distances between feature vectors and 
should gain from the manifold denoising.
For each of the data sets we perform the following procedure.
First, we apply k-NN method without denoising.
Then we make manifold denoising using SAME and apply k-NN to the denoised data set.
In the case of g241c data set, we took \( d = 10 \), \( \tau = 22 \), \( \gamma = 4 \) and \( 
h_k = 20 \cdot 1.2^{-k} \), \( 0 \leq k \leq 1 \).
In the case of g241n data set, we took \( d = 6 \), \( \tau = 21 \), \( \gamma = 4 \) and \( h_k 
= 20 \cdot 1.2^{-k} \), \( 0 \leq k \leq 2 \).
The results are summarized in Table \ref{table1}.
We observe that preliminary denoising improves quality of prediction.

\begin{table}
	\noindent\centering
	{\renewcommand{\arraystretch}{2}
		\begin{tabular}{|p{6em}|p{6em}|p{6em}|p{6em}|p{6em}|}
			\hline
			Data set & Best number of neighbors \( k \) without denoising & k-NN error 
			without denoising (\%) & Best number of neighbors \( k \) after denoising & 
			denoised k-NN error (\%) \\
			\hline
			g241c & 21 & 31.3 & 15 & 27.8 \\
			\hline
			g241n & 18 & 27.1 & 12 & 25.3 \\
			\hline
	\end{tabular}}
	\caption{Error rates with and without manifold denoising via SAME for k-NN 
	method,  applied to artificial data sets g241c and g241n.}
	\label{table1}
\end{table}

\section{Theoretical Properties of SAME}
\label{sec_theoretical}
This section states the main results. Here and everywhere in this paper, for any matrix \( 
\boldsymbol{A} \), \( \|\boldsymbol{A}\| \) denotes its spectral norm.
The notation \( f(n) \asymp g(n) \) means \( f(n) \lesssim g(n) \lesssim f(n) \).

\begin{Th}
	\label{th1}
	Assume \eqref{a1}, \eqref{a1'}, \eqref{a2}, and \eqref{a3} .
	Let the initial guesses \( \smash{\widehat{\bpi}\mathstrut_1\ind0, \dots, \widehat{\bpi}\mathstrut_n\ind0} \) of 
	\( \bpi(X_1), \dots, \bpi(X_n) \) be such that on an event with probability at least \( 1 - 
	n^{-1} \) it holds
	\[
		\max\limits_{1 \leq i \leq n} \|\widehat{\bpi}_{i}\ind0 - \bpi(X_{i})\| \leq \frac{\Delta 
		h_0}\varkappa
	\]
	with a constant \( \Delta \), such that \( \Delta h_0 \leq \varkappa/4 \), 
	and \(h_0 = C_0/\log n\), where \( C_0 >0 \) is an absolute constant.
	Choose \( \tau = 2C_0/\sqrt{\log n} \) and set any \( a \in (1, 2] \).
	If \( n \) is larger than a constant \( N_{\Delta} \), depending on \( \Delta \), and \( 
	h_{K} \gtrsim \left((D\log n/n)^{1/d} \vee (DM^{2}\varkappa^2 \log n/n)^{1/(d+4)}\right) 
	\) (with a sufficiently large hidden constant, which is greater than \( 1 \)) then there 
	exists a choice of \( \gamma \), such that after \( K \) iterations Algorithm 
	\ref{algorithm} produces estimates \( \widehat{X}_1, \dots, \widehat{X}_n \), such 
	that, with probability at least \( 1 - (5K + 4)/n \), it holds
	\begin{align*}
		&
		\max\limits_{1\leq i \leq n} \|\widehat{X}_{i} - X_{i} \| \lesssim
		\frac{Mb \vee M h_K \vee (1 + \Phi_{M, b, h_K, \varkappa}) 
		h_K^2}\varkappa + \sqrt{\frac{D(h_{K}^{2}\vee M^{2}) \log n}{nh_{K}^d}},
		\\&
		\max\limits_{1\leq i \leq n} \|\widehat{\bpi}_{i}\ind K - \bpi(X_i) \| \lesssim 
		\Psi_{M, b, h_K, \varkappa} \left(\frac{h_K}\varkappa + h_K^{-1} \sqrt{\frac{D 
		(h_K^{4}/\varkappa^2 \vee M^{2}) \log n}{nh_K^d}}\right),	
	\end{align*}
	where
	\begin{align}
		\label{phi_psi}
		\Phi_{M, b, h_K, \varkappa}
		&\notag
		= \frac{M^3(1 + b/h_K)^2}{h_K^2 \varkappa}
		+ \frac{M^2(1 + b/h_K + \sqrt{\log h_K^{-1}})}{\varkappa h_K}
		+ \frac{M h_K^{2}}{\varkappa^3}
		\lesssim \alpha + o(1), \quad n \rightarrow \infty,
		\\
		\Psi_{M, b, h_K, \varkappa}
		&
		= \left(1 + \frac{M(1 + b/h_K) \vee (1 + \Phi_{M, b, h_K, 
		\varkappa}) h_K}\varkappa \right)^{d+1} (1 + \Phi_{M, b, h_K, \varkappa})
		\\&\notag
		\leq (1+\alpha)\left( 4^{d+1} + (2\sqrt{\alpha})^{d+1} \right).
	\end{align}
	In particular, if one chooses the parameter \( a \) and the number of iterations \( K \) 
	in such a way that \( h_K \asymp \left((D\varkappa^2 \log n/ 
	n)^{1/(d+2)}\right.\) \(\left.\vee (DM^{2}\varkappa^2 \log n/n)^{1/(d+4)}\right) \) then
	\[
		\max\limits_{1\leq i \leq n} \|\widehat{X}_{i} - X_{i} \| \lesssim 
		\frac{Mb}\varkappa + \frac1\varkappa \left( \frac{D\varkappa^2 \log 
		n}n\right)^{\frac2{d+2}} \vee \frac{M}\varkappa \left(\frac{DM^2\varkappa^2 \log 
		n}n \right)^{\frac1{d+4}}.
	\]
	If \( h_K \asymp \left((D\log n/n)^{1/d} \vee (DM^{2}\varkappa^2 \log 
	n/n)^{1/(d+4)}\right) \) then
	\[
	\max\limits_{1\leq i \leq n} \|\widehat{\bpi}_{i}\ind K - \bpi(X_i) \| \lesssim 
	\frac1\varkappa \left(\frac{D\log n}n \right)^{\frac1d} \vee
	\frac1\varkappa \left(\frac{DM^{2}\varkappa^2 \log 
		n}n\right)^{\frac1{d+4}}.
	\]
\end{Th}

Note that one has to take the number of iterations \( K \) of order \( \log n \) since the 
sequence of bandwidths \( h_1, \dots, h_K \) decreases exponentially.

In Theorem \ref{th1}, we assume that \( \smash{\widehat\bpi\mathstrut_1\ind{0}, 
\dots, \widehat\bpi\mathstrut_n\ind{0}} \) may depend on \( Y_1, \dots, Y_n \).
The natural question is how to construct the initial guesses \( \smash{\widehat\bpi\mathstrut_1\ind{0}, 
\dots, \widehat\bpi\mathstrut_n\ind{0}} \) of the projectors \( \bpi(X_1), \dots, \bpi(X_n) \).
We propose a strategy for initialization of our procedure.
One can use \citep[Proposition 5.1]{al18} to get the estimates 
\( \smash{\widehat\bpi\mathstrut_1\ind{0}, \dots, \widehat\bpi\mathstrut_n\ind{0}} \).
For each \( i \) from \( 1 \) to \( n \) introduce
\[
	\widehat\bsigma_i\ind 0 = \frac1{n-1} \sum\limits_{j \neq i} (Y_j - \overline Y_i)(Y_j - 
	\overline Y_i)^T \1(Y_j \in \B(Y_i, h_0)), 
\]
where \( \overline Y_i = \frac1{N_i} \sum_{j \neq i} Y_j \1(Y_j \in \B(Y_i, h_0)) \), \( N_i 
= |\{ j : Y_j \in \B(Y_i, h)\}| \).
Let \( \smash{\widehat \bpi\mathstrut_i\ind 0} \) be the projector onto the linear span of the \( d \) largest 
eigenvalues of \( \smash{\widehat\bsigma\mathstrut_i\ind 0} \). Then the following result holds.
\begin{Prop}[\cite{al18}, Proposition 5.1]
	Assume \eqref{a1}, \eqref{a1'}, \eqref{a2}.
	Set \( h_0 \gtrsim (\log n/ n)^{1/d} \) for large enough hidden constant.
	Let \( M/{h_0} \leq 1/4 \) and let \( h_0 = h_0(n) = o(1) \), \( n \rightarrow \infty \).
	Then for \( n \) large enough, with probability larger than \( 1 - n^{-1} \), it holds
	\[
		\max\limits_{1 \leq i \leq n} \| \widehat\bpi_i\ind 0 - \bpi(X_i) \| \lesssim 
		\frac{h_0}\varkappa + \frac{M}{h_0}.
	\]
\end{Prop}

\begin{Rem}
	In \citep{al18}, the authors take \( h_0 \asymp (\log n / n)^{1/d} \).
	Nevertheless, a careful reading of the proofs reveals that one can also take larger 
	values of \( h_0 \).
\end{Rem}

\begin{Rem}
	One can also use local PCA procedure from \citep{cw13} with a more sophisticated 
	choice of the neighborhood for initialization.	
\end{Rem}

The condition \eqref{a3} and the choice of \( h_K \) in Theorem \ref{th1} yield that \( M 
= M(n) \) can decrease almost as slow as \( h_K^{2/3} = h_K^{2/3}(n) \).
Thus, we admit the situation when the noise magnitude \( M \) is much larger than the 
smoothing parameter \( h_K \).
For instance, in \citep{al19}, the authors use local polynomial estimates and require \( M = 
O(h^2) \) and \( h = h(n) \asymp n^{-1/d} \). 
In \citep{al18}, the authors assume \( M \leq \lambda (\log n/n)^{1/d} \), and \( \lambda \) does 
not exceed a constant \( \lambda_{d, p_0, p_1} \), depending on \( d \), \( p_0 \) and \( p_1 \).
In \citep{fikln18}, the authors deal with Gaussian noise \( \N(0, \sigma^2 \bid_D) \) and get 
the accuracy of manifold estimation \( O(\sigma\sqrt{D}) \) using \( O(\sigma^{-d}) \) 
samples.
This means that \( \sigma = O(n^{-1/d}) \), which yields that
\[
	\max\limits_{1 \leq i \leq n} \|\eps_i\| \lesssim n^{-1/d} \sqrt{D\log n}
\]
with overwhelming probability.
A similar situation is observed in \citep{gppvw14}, where the authors also consider the 
Gaussian noise \( \N(0, \sigma^2 \bid_D) \) and, using the kernel density estimate with 
bandwidth \( h \), obtain the upper bound
\[
	O\left(\sigma^2\log\sigma^{-1}+ h^2 + \sqrt{\frac{\log n}{nh^D}} \right)
\]
on the Hausdorff distance between \( \M^* \) and their estimate.
In order to balance the first and the second terms, one must take \( \sigma = 
O(h/\sqrt{\log h^{-1}}) \), which means that
\[
	\max\limits_{1 \leq i \leq n} \|\eps_i\| \lesssim h\sqrt{\frac{D\log n}{\log h^{-1}}},
\]
while we allow \( \max_{1 \leq i \leq n} \|\eps_i\| \) be as large as \( h_K^{2/3} \).
Finally, in \citep{hm06} the authors require \( M = O(h) \).
So, we see that the condition \eqref{a3} is quite mild.

Theorem \ref{th1} claims that, despite the relatively large noise, our procedure 
constructs consistent estimates of the projections of the sample points onto the 
manifold \( \M^* \).
The accuracy of the projection estimation is a bit worse than the accuracy of manifold 
estimation, which we provide in Theorem \ref{th2} below.
The reason for that is the fact that the estimate \( \widehat{X}_{i} \) is significantly 
shifted with respect to \( X_{i} \) in a tangent direction, while the orthogonal component 
of \( (\widehat{X}_{i} - X_{i}) \) is small.
A similar phenomenon was already known in the problem of efficient dimension 
reduction.
For instance, in \citep{hjs01a, hjs01b} the authors managed to obtain the rate 
\( n^{-2/3} \) for the bias of the component, which is orthogonal to the efficient 
dimension reduction space, while the rate of the bias in the index estimation was only \( 
n^{-1/2} \).
Moreover, the term \( Mh_K \) in Theorem \ref{th1} appears because of the correlation 
between the weights \( w_{ij}\ind k \) and the sample points \( Y_j \).

We proceed with upper bounds on the estimation of the manifold \( \M^* \).
\begin{Th}
	\label{th2}
	Assume conditions of Theorem \ref{th1}.
	Consider the piecewise linear manifold estimate
	\[
		\widehat{\M} = \left\{ \widehat{X}_{i} + h_K \widehat{\bpi}_{i}\ind K u : 1 \leq 
		i \leq n, \, u \in \B(0, 1) \subset \R^D \right\},
	\]
	where \( \smash{\widehat{\bpi}\mathstrut_{i}\ind K} \) is a projector onto \( d \)-dimensional space 
	obtained on the K-th iteration of Algorithm \ref{algorithm}.
	Then, as long as \( 
	h_{K} \gtrsim \left((D\log n/n)^{1/d} \vee (DM^{2}\varkappa^2 \log n/n)^{1/(d+4)}\right) 
	\) (with a sufficiently large hidden constant, which is greater than \( 1 \)), on an event 
	with probability at least 
	\( 1 - (5K + 5)/n \), it holds
	\[
		d_H(\widehat{\M}, \M^{*}) \lesssim \left( \frac{(1 + \Phi_{M, b, h_K, \varkappa} + 
		\Psi_{M, b, h_K, \varkappa})h_{K}^{2}}\varkappa \vee \frac{M^2 
		b^2}{\varkappa^3} \right) + \sqrt{\frac{D(h_{K}^4/\varkappa^2\vee M^{2}) \log 
		h_{K}^{-1}}{nh_{K}^d}},
	\]
	where $\Phi_{M, b, h_K, \varkappa}$ and $\Psi_{M, b, h_K, \varkappa}$ are defined in 
	\eqref{phi_psi}.
	In particular, if \( a \) and \( K \) are chosen such that \( h_K \asymp \left((D \log 
	n/n)^{1/d} \vee (D M^{2} \varkappa^2 /n \log n)^{1/(d+4)}\right) \), then
	\[
		d_H(\widehat{\M}, \M^{*}) \lesssim \frac{M^2 
		b^2}{\varkappa^3}  \vee \varkappa^{-1}\left( \frac{D \log 
		n}{n}\right)^{\frac{2}{d}} \vee \varkappa^{-1} \left(\frac{DM^{2}\varkappa^2\log 
		n}n \right)^{\frac2{d+4}}.
	\]
\end{Th}

\begin{Rem}
	For manifold reconstruction, one can also use a different technique, based on \citep[Theorem 4.1]{al18} and tangential Delaunay complexes. We emphasize that the manifold estimate $\widehat\M$, used in Theorem \ref{th2}, is not a manifold itself but a union of $d$-dimensional discs.
\end{Rem}

Let us elaborate on the result of Theorem \ref{th2}.
First, let us discuss the case of bounded non-orthogonal noise, that is, the situation 
when \eqref{a3''} holds.
The model with bounded noise was considered in \citep{al18}, where the authors 
assumed that $\M^*$ satisfies \eqref{a1} and the density of $X$ fulfils
\[
	0 < p_0 \leq p(x) \leq p_1, \quad \forall \, x \in \M^*,
\]
for some constants $p_0, p_1$.
Note that this is a slightly more general setup, since we additionally assume that the 
log-density is Lipschitz.
Under these assumptions, \cite{al18} proved (Theorem 2.7) the following upper bound on 
the Hausdorff distance using the tangential Delaunay complex (TDC):
\[
	d_H(\widehat\M_{TDC}, \M^*) \lesssim \left(\frac{\log n}n\right)^{2/d} + M^2 
	\left(\frac{\log n}n\right)^{-2/d},
\]
provided that \( M \leq \lambda_{d, p_0, p_1} (\log n/n)^{1/d} \), where the constant \( 
\lambda_{d, p_0, p_1} \) depends on \( d \), \( p_0 \) and \( p_1 \).
To the best of our knowledge, the situation, when \eqref{a1}, \eqref{a1'}, \eqref{a2}, and 
\eqref{a3''} hold, was not studied in the manifold learning literature.
One can observe that both TDC and SAME achieve the rate \( O\left(\log n/n\right)^{2/d} 
\) in the case of extremely small noise \( M \lesssim (\log n/n)^{2/d}\).
However, if \( (\log n/n)^{2/d} \lesssim M \lesssim n^{-4/(3d+4)} \) then the rate of 
convergence of SAME in the case of the density $p(x)$ satisfying \eqref{a1'} improves 
over the known rates of TDC in the case of bounded away from $0$ and $\infty$ density 
$p(x)$.

Now, let us discuss the case of almost orthogonal noise, i.e. when \eqref{a3} holds.
This model is completely new in the manifold learning literature.
The most similar one considered in the prior work is the model with perpendicular noise 
studied in \citep{gppvw12a, al19}, so we find it useful to compare this more restrictive 
model with our upper bounds for the case of almost orthogonal noise.
In \citep{gppvw12a}, the authors obtain the rates \( O(\log n/n)^{2/(d+2)} \) assuming that, 
given \( X \), the noise \( \eps \) has a uniform distribution on \( \B(X, M) \cap 
(\T_X\M^*)^\perp \).
In their work, the authors do not assume that \( M \) tends to zero as \( n \) tends to infinity, 
however, they put a far more restrictive assumption on the noise distribution than we do.
In \citep[Theorem 6]{al19}, the authors use local polynomial estimate \( \widehat\M_{LP} \) 
to prove the upper bound 
\[
	d_H(\widehat\M_{LP}, \M^*) \lesssim \left( \frac{\log n}n \right)^{k/d} \vee M
\]
for the case when \( \M^* \) is a \( \C^k \)-manifold with dimension \( d \) and reach at least 
\( \varkappa \) without a boundary.
If \( \M^* \) is a \( \C^2 \)-manifold, this rate is minimax optimal for the case of extremely 
small noise \( M \lesssim (\log n/n)^{2/d} \) but it can be improved when the noise 
magnitude exceeds \( (\log n/n)^{2/d} \).

Theorem \ref{th2} shows that our procedure achieves the classical nonparametric rate, 
where the bias and the variance terms correspond to the best one can hope for when 
deals with the locally linear estimator.
In the case of small noise (\( M \lesssim (\log n/n)^{2/d} \)), the result of Theorem \ref{th2} 
matches the lower bound obtained in \citep{kz15} and the upper bound from \citep{al19}.
It is not surprising, because when \( \M^* \) is a \( \C^2 \)-manifold, the local polynomial 
estimate considered in \citep{al19} becomes a piecewise linear estimate, based on local 
PCA, and achieves the optimal rate in the case of small noise.
Our algorithm acts in a similar manner and the only significant difference is hidden in the 
weights.
However, if the noise is very small, there is no need to adjust the weights, so local PCA 
and SAME behave comparably in this regime.

The same concerns the projector estimates \( \smash{\widehat\bpi\mathstrut_1\ind K, \dots, 
\widehat\bpi\mathstrut_n\ind K} \) from Theorem \ref{th1}.
In the case of small noise, we recover the minimax rate \( (\log n / n)^{1/d} \) obtained in 
\citep{al18, al19}.
However, as the magnitude of the noise grows, our procedure shows superior 
performance, compared to the estimates in \citep{al18, al19}.

The result of Theorem \ref{th2} cannot be improved for the case of general additive 
noise, which fulfils the assumption \eqref{a2} with \( b \gtrsim \left((\log n/n)^{1/d} \vee 
\right. \) \( \left.(M^{2}\varkappa^2 \log n/n)^{1/(d+4)}\right) \).
We justify this discussion by the following theorem.

\begin{Th}
	\label{lower_bound}
	Suppose that the sample \( \Y_n = \{Y_1, \dots, Y_n\} \) is generated according to the 
	model \eqref{model}, where \( \M^* \in \mclass_\varkappa^d \), the density \( p(x) \) 
	of \( X \) 
	fulfils \eqref{a1'} (with sufficiently large $p_1, L$ and sufficiently small $p_0$) and 
	the noise \( \eps \)  satisfies \eqref{a2}.
	Then, for any estimate $\widehat \M$, it holds that
	\begin{equation}
	\label{lb_1}
		\sup\limits_{\M^* \in \mclass_\varkappa^d} \E_{\M^*} 
		d_H(\widehat\M, \M^*)
		\gtrsim \frac{M^2 b^2}{\varkappa^3}.
	\end{equation}
	Moreover, if, in addition, \( n \) is sufficiently large, \( M\varkappa \gtrsim (\log 
	n/n)^{2/d} \), and the parameter $b$ in \eqref{a2} is such that
	\[
		b \gtrsim \left((\log n/n)^{1/d} \vee (M^{2}\varkappa^2 \log n/n)^{1/(d+4)}\right),
	\]
	with a  large enough hidden constant,
	then, for any estimate $\widehat \M$, it holds that
	\begin{equation}
		\label{lb_2}
		\sup\limits_{\M^* \in \mclass_\varkappa^d} \E_{\M^*} 
		d_H(\widehat\M, \M^*)
		\gtrsim \varkappa^{-1} \left(\frac{M^2\varkappa^2 \log n}n \right)^{\frac 2{d+4}}.
	\end{equation}
\end{Th}
Theorem \ref{lower_bound} studies the case \( M \gtrsim (\log n/n)^{2/d} \).
In \citep{kz15}, the authors proved the minimax lower bound
\[
	\inf\limits_{\widehat\M} \sup\limits_{\M^* \in \mclass_\varkappa^d} \E_{\M^*} 
	d_H(\widehat\M, \M^*)
	\gtrsim \left( \frac{\log n}n \right)^{2/d}
\]
for the noiseless case, which is also tight for \( M \lesssim (\log n/n)^{2/d} \).
Theorem \ref{lower_bound}, together with \citep[Theorem 1]{kz15} yields 
that SAME is minimax optimal in the model with almost orthogonal noise.
The lower bounds \eqref{lb_1} and \eqref{lb_2} are completely new and are different from 
the currently known results on manifold estimation from \citep{gppvw12a} and \citep{al19}, 
where the authors studied a perpendicular noise fulfilling \eqref{a2} with \( b = 0 \).
In \citep{gppvw12a}, the authors focused on the case of uniform noise and 
proved the lower bound \( (M/n)^{2/(d+2)} \) when \( \M^* \) fulfils \eqref{a1} and 
\( p_0 \leq p(x) \leq p_1 \) for all \( x \in \M^* \).
In \citep{al19}, the authors went further and proved that
\[
	\inf\limits_{\widehat\M} \sup\limits_{\M^*} \E 
	d_H(\widehat\M, \M^*) \gtrsim (M/n)^{k/(d+k)} \vee n^{-k/d},
\]
where \( \M^* \) runs over a class of compact, connected \(\C^k\)-manifolds of dimension 
\( d \) without a boundary and \( \reach{\M^*} \geq \varkappa \).
Theorem \ref{lower_bound} reveals a surprising effect: if one allows small deviations of 
\( \eps_i \)'s in tangent directions then the problem of manifold estimation becomes 
harder and this fact is reflected in the minimax rates of convergence.
Namely, if \( b \gtrsim (\log n/n)^{1/d} \vee (M^{2}\varkappa^2 \log 
n/n)^{1/(d+4)} \), the minimax rate is \( M^2 b^2 / \varkappa^3 \vee \left(M^2 \log n/n 
\right)^{2/(d+4)} \vee \left(\log n/n \right)^{2/d} \) which is different from the best 
known lower bound \( \left(M/n\right)^{2/(d+2)} \vee \left(\log n/n \right)^{2/d} \) for the 
case of perpendicular noise.

\section{Proofs}
\label{proofs}
This section collects the proof of the main results.
\subsection{Proof of Theorem \ref{th1}}
\label{sec_proof1}

The proof of Theorem \ref{th1} is given in several steps.
First, we show that the adjusted weights \( w_{ij}(\bpi_{i}) \) are informative, i.e. 
significant weights correspond only to points \( X_{j}\), which are close to \( X_{i}\).
\begin{Lem}
	\label{geom}
	Assume \eqref{a1}, \eqref{a2}.
	Let \( \bpi_{i} \) be any projector, such that \( \|\bpi_{i} - \bpi(X_{i})\| \leq \frac{\Delta 
		h}\varkappa\).
	Assume that \( M \leq \varkappa/16 \), \( \Delta h\leq 
	\varkappa/4 \), and \( M(\Delta + b/h) \leq \varkappa/4 \).
	Then for any \( i \) and \( j \), such that \( \| Y_{i} - Y_{j} \| \leq 0.5\varkappa \), it holds
	\[
		\frac{1}{2} \|\bpi_{i}(Y_{i} - Y_{j})\| - \frac{2M(\Delta h + b)}\varkappa \leq 
		\|X_{i} - X_{j} \| \leq 2 \|\bpi_{i}(Y_{i} - Y_{j})\| + \frac{4M(\Delta h + b)}\varkappa.
	\]
\end{Lem}

Lemma \ref{geom} quantifies the informal statement \( \bpi_i(Y_j - Y_i) \approx X_j - X_i \), 
giving explicit error bound depending on the error in the guess of the projector.
Note that if \eqref{a3} holds, \( h \in [h_K, h_0] \) and \( h_0, h_K \) 
satisfy the conditions of Theorem \ref{th1} then \( M \leq \varkappa/16 \), \( M(\Delta + b/h) 
\leq \varkappa/4 \), \( \Delta h \leq \varkappa/4 \) if \( n \) is large enough.
The next step is to show that the cylinder \( \{y : \|\bpi_i(y - Y_i)\| \leq h, \| y - Y_i \| \leq \tau 
\} \) contains enough sample points. 
For this purpose, we prove the regularity of the design points in the following sense.
\begin{Lem}
	\label{regularity}
	Assume \eqref{a1}---\eqref{a2}.
	Fix any \( i \) from \( 1 \) to \( n \) and let \( \bpi_{i} \) be any projector, such that \( 
	\|\bpi_{i} - \bpi(X_{i})\| \leq \frac{\Delta h}\varkappa\), where \( (\log n/n)^{1/d} \lesssim 
	h < h_0 \), \( \Delta h\leq \varkappa/4 \), \( M\leq \varkappa/16 \), and \( M 
	(\Delta + b/h) \leq \varkappa/4 \).
	Suppose that \( h_0 \) is chosen in a such way that \( h_0 \leq 0.5 \tau \), and \( n \) is 
	sufficiently large. 
	Then, on an event with probability at least \( 1 - n^{-2} \), it holds
	\begin{equation}
		\label{regularity_eq}
		\sum\limits_{j=1}^{n} w_{ij}(\bpi_{i}) \geq C' nh^d
	\end{equation}
	with an absolute constant \( C' > 0 \).
\end{Lem}
Roughly speaking, the cylinder \( \{y : \|\bpi_i(y - Y_i)\| \leq h, \| y - Y_i \| \leq \tau 
\} \) contains \( \sim n h^d \) sample points with high probability.
It is important, because the ball \( \B(Y_i, h) \) would contain \( \sim n h^D \) sample points if 
\( h < M \).
The sum of weights controls the variance of our estimates.
From this point of view, the choice of cylindric neighborhoods instead of the balls yields 
much better rates.

Now, we are ready to make the main step in the proof of Theorem \ref{th1}.
\begin{Lem}
	\label{main_lem}
	Assume conditions of Theorem \ref{th1}.
	Let \( \bpi_1, \dots, \bpi_n \) be any (possibly random) projectors, such that \( 
	\|\bpi_{i} - \bpi(X_{i})\| \leq \frac{\Delta h}\varkappa\) almost surely, \( (\log n/n)^{1/d} 
	\lesssim h \leq h_0 \), \( \Delta h\leq 
	\varkappa/4 \), \( M \leq \varkappa/16 \), and \( M (\Delta + b/h) \leq 
	\varkappa/4 \).
	Let \( w_{ij}(\bpi_{i}) \), \( 1 \leq i, j \leq n \), be the localizing weights, computed 
	according to
	\[
	w_{ij}(\bpi_{i}) = \K\left( \frac{\|\bpi_{i}(Y_{i} - Y_{j})\|^2}{h^2} \right) \1\left( \| 
	Y_{i} - Y_{j} \| \leq 
	\tau \right), \quad 1 \leq  j \leq n,
	\]
	with a constant \( \tau < 0.5\varkappa \).
	Then, conditionally on \( \bpi_1, \dots, \bpi_n \), on an event with probability at least \( 
	1 - 2n^{-1} \), it simultaneously holds
	\begin{align*}
		&
		\max\limits_{1\leq i \leq n} \left\|\sum\limits_{j=1}^{n} w_{ij}(\bpi_{i})(Y_{j} - X_{i}) 
		\right\| 
		\lesssim \left(M\left(\Delta + \frac bh \right) \vee h \vee \frac{\Delta^2 
		h^2}\varkappa \right) 
		\frac{h^{d+1}}\varkappa
		\\&
		+ \Phi_{M, b, h, \varkappa, \Delta} 
		\frac{nh^{d+2}}\varkappa + \sqrt{D(h^{2}\vee M^{2}) nh^d \log n},
		\\&
		\max\limits_{1\leq i \leq n} \left\|\sum\limits_{j=1}^{n} w_{ij}(\bpi_{i})\big((\bid - 
		\bpi(X_{i}))(Y_{j} - X_{i}) \big)\right\| 
		\\&
		\lesssim \left(1 + \Phi_{M, b, h, \varkappa, \Delta} \right)
		\frac{nh^{d+2}}\varkappa + \sqrt{D(h^4/\varkappa^2 \vee M^{2}) nh^d \log n},
	\end{align*}
	where
	\[
		\Phi_{M, b, h, \varkappa, \Delta} =
		\frac{M^3(1 + \Delta + b/h)^2}{h^2 \varkappa}
		+ \frac{M^2(\Delta + b/h + \sqrt{\log h^{-1}})}{\varkappa h}
		+ \frac{(1 + \Delta^4) M h^{2}}{\varkappa^3},
	\]
	and the hidden constants do not depend on \( \Delta \).
\end{Lem}
The proof of Lemma \ref{main_lem} is moved to Appendix \ref{app_c}.
In Lemma \ref{main_lem}, the assumption \eqref{a3} comes into play.
The condition \eqref{a3} implies that \( \Phi_{M, b, h, \varkappa, \Delta} \leq \alpha +  o(1) 
\) as \( n\rightarrow \infty \).
This follows from the fact that, under \eqref{a3}, for any \( h  = h(n) \geq h_K \), it holds 
\( M^3 = o(h^2) \), \( n \rightarrow\infty \).
Indeed, we have
\begin{align*}
	\frac{M^3}{h^2}
	&
	\leq \frac{M^3}{h_K^2}
	= \frac{M^3}{(M^2/n)^{2/(d+4)}} \cdot \frac{(M^2/n)^{2/(d+4)}}{h_K^2}
	\\&
	= M^{(3d+8)/(d+4)}\cdot n^{2/(d+4)} \cdot \frac{(M^2/n)^{2/(d+4)}}{h_K^2}
	\\&
	\leq A \cdot \frac{(M^2/n)^{2/(d+4)}}{h_K^2}
	\leq \frac A{(\log n)^{2/(d+4)}}
	\rightarrow 0, \quad n \rightarrow \infty.
\end{align*}
Moreover, under \eqref{a3}, it holds
\[
	\frac{M^3 b^2}{\varkappa h^4}
	\leq \frac{M^3 b^2}{\varkappa h_K^4}
	\leq \frac{M^3 b^2}{\varkappa\left(\frac{D\log n}n \right)^{\frac4d} \vee
	\varkappa \left(\frac{DM^{2}\varkappa^2 \log n}n\right)^{\frac4{d+4}}}
	\leq \alpha.
\]
Therefore, \( \Phi_{M, b, h, \varkappa, \Delta} \leq \alpha + o(1) \) as \( n\rightarrow\infty \).
Similarly, if \eqref{a3''} holds, we have \( M^3 b^2  = M^3 \varkappa^2 = o(h_K^4) \), 
which also yields \( \Phi_{M, b, h, \varkappa, \Delta} \rightarrow 0 \).

We need one more auxiliary result.
Lemma \ref{geom}, \ref{regularity} and \ref{main_lem} imply that if we have good 
guesses of the projectors \( \bpi(X_1), \dots, \bpi(X_n) \) then we can get good estimates 
of \( X_1, \dots, X_n \) even if the noise magnitude \( M \) is quite large.
Now we have to show that these estimates of \( X_1, \dots, X_n \) can be used to construct 
good estimates of \( \bpi(X_1), \dots, \bpi(X_n) \) and then we can carry out the proof by 
induction.
Our next lemma is devoted to this problem.

Let us work on the event, on which \eqref{regularity_eq} holds with \( h = h_k \).
Note that
\[
	\|\widehat{X}_{i}\ind{k} - X_{i} \|
	= \left\| \frac{\sum\limits_{j=1}^{n} w_{ij}\ind{k} Y_{j}}{\sum\limits_{j=1}^{n} w_{ij}\ind{k}} - 
	X_{i} \right\|
	= \frac{\left\| \sum\limits_{j=1}^{n} w_{ij}\ind{k} (Y_{j} - X_{i}) 
	\right\|}{\sum\limits_{j=1}^{n} w_{ij}\ind k}
\]
and
\begin{align*}
	d(\widehat{X}_{i}\ind{k}, \{X_i\} \oplus \T_{X_{i}}\M^{*})
	&
	= \left\|\widehat{X}_{i}\ind{k} - X_{i} - \bpi(X_{i})(\widehat{X}_{i}\ind{k} - X_{i}) \right\|
	\\&
	= \frac{ \left\| \sum\limits_{j=1}^{n} w_{ij}\ind{k} \big(Y_{j} - X_{i} - 
	\bpi(X_{i})(Y_{j} - X_{i}) \big) \right\| }{\sum\limits_{j=1}^{n} w_{ij}\ind{k}} \, .
\end{align*}
Here we used the fact that the projection of a point \( x \) onto the tangent plane \( \{X_i\} 
\oplus \T_{X_i}\M^* \) is given by
\[
	\pi_{\{X_i\} \oplus \T_{X_i}\M^*}(x) = X_i + \bpi(X_i)(x - X_i).
\]
Then Lemma \ref{main_lem} and Lemma \ref{regularity} immediately yield that, if we 
have
\[
	\max\limits_{1\leq i\leq n}\| \widehat\bpi_{i}\ind{k} - \bpi(X_{i}) \| \leq \Delta h_{k} / \varkappa
\]
on the \( k \)-th iteration with probability at least \( 1 - (5k+1)/n \), then
\begin{align*}
	&
	\max\limits_{1\leq i \leq n} \|\widehat{X}_{i}\ind{k} - X_{i}\| \lesssim  
	\left(\frac{M(\Delta h_k + b) \vee (1 + \Phi_{M, b, h_k, \varkappa, \Delta}) 
	h_k^2}\varkappa \vee \frac{\Delta^2h_k^3}{\varkappa^2}\right)
	+ \sqrt{\frac{D(h_{k}^{2}\vee M^{2}) \log h_{k}^{-1}}{nh_{k}^d}},
	\\&\notag
	\max\limits_{1\leq i \leq n} d(\widehat{X}_{i}\ind{k}, \{X_i\} \oplus  \T_{X_{i}}\M^{*}) 
	\lesssim \frac{(1 + \Phi_{M, b, h_k, \varkappa, \Delta})h_{k}^{2}}\varkappa + 
	\sqrt{\frac{D(h_{k}^4/\varkappa^2\vee M^{2}) \log h_{k}^{-1}}{nh_{k}^d}}
\end{align*}
with probability at least \( 1 - (5k+1)/n - 3/n = 1 - (5k+4)/n \).
It only remains to check that the projector estimates \( \smash{\widehat{\bpi}\mathstrut_1\ind{k+1}, 
\dots, \widehat{\bpi}\mathstrut_n\ind{k+1}} \) also satisfy
\[
	\max\limits_{1\leq i\leq n} \|\widehat\bpi_{i}\ind{k+1} - \bpi(X_{i})\| \lesssim 
	\frac{h_{k+1}}\varkappa
\]
with high probability.
The precise statement is given in the following lemma.

\begin{Lem}
	\label{proj}
	Assume conditions of Theorem \ref{th1}.
	Let \( \Omega_{k} \) be an event, such that on this event it holds
	\begin{align}
		\label{beta}
		&\notag
		\max\limits_{1\leq i \leq n} \|\widehat{X}_{i}\ind{k} - X_{i}\|
		\leq \beta_1 \left( h_k + \sqrt{\frac{D(h_{k}^{2}\vee M^{2}) \log 
		h_{k}^{-1}}{nh_{k}^d}} \right)
		\leq 2\beta_1 h_k ,
		\\&
		\max\limits_{1\leq i \leq n} d(\widehat{X}_{i}\ind{k}, \{X_i\} \oplus \T_{X_{i}}\M^{*}) 
		\leq \beta_2 \left( \frac{h_{k}^{2}}\varkappa + 
		\sqrt{\frac{D(h_{k}^4/\varkappa^2\vee M^{2}) \log h_{k}^{-1}}{nh_{k}^d}} \right)
		\leq \frac{2\beta_2 h_k^2}\varkappa.
	\end{align}
	Then there exists \( \gamma \asymp 1 + \beta_1 \) such that, with probability at least 
	\( \p(\Omega_k) - 2n^{-1} \), it holds
	\begin{align*}
		\max\limits_{1\leq i \leq n} \|\widehat{\bpi}_{i}\ind{k+1} - \bpi(X_{i})\| 
		\lesssim \frac{\gamma(\gamma + 4\beta_1)^d \beta_2 h_{k}}\varkappa
		\lesssim \frac{(1 + \beta_1)^{d+1} \beta_2 h_{k}}\varkappa.
	\end{align*}
\end{Lem}
The proof of Lemma \ref{proj} can be found in Appendix \ref{app_d}.
From the derivations before Lemma \ref{proj}, the event \( \Omega_k \) from Lemma 
\ref{proj} has probability at least \( 1 - (5k+4)/n \).
Note that, if
\( h_{k} \gtrsim \left((D\log n/n)^{1/d} \vee (DM^{2}\varkappa^2\log n/n)^{1/(d+4)}\right) \) 
with a hidden constant greater than \( 1 \), as given in the conditions of Theorem 
\ref{th1}, then the bias terms in \eqref{beta} are dominating, i.e.
\[
	h_k + \sqrt{\frac{D(h_{k}^{2}\vee M^{2}) \log h_{k}^{-1}}{nh_{k}^d}}
	\leq 2 h_{k}
	= 2a h_{k+1}
\]
and
\[
	\frac{h_{k}^{2}}\varkappa + \sqrt{\frac{D(h_{k}^4/\varkappa^2\vee M^{2}) \log 
	h_{k}^{-1}}{nh_{k}^d}}
	\leq \frac{2 h_{k}^{2}}\varkappa.
\]
Due to the discussion before Lemma \ref{proj}, we can take
\begin{align*}
	&\notag
	\beta_1 = \left(\frac{M(\Delta + b/h_k) \vee (1 + \Phi_{M, b, h_k, \varkappa, \Delta}) 
	h_k}\varkappa \vee \frac{\Delta^2 h_k^2}{\varkappa^2}\right),
	\\&
	\beta_2 = 1 + \Phi_{M, b, h_k, \varkappa, \Delta}.
\end{align*}
Then Lemma \ref{proj} yields
\begin{equation}	
	\label{proj_ineq}
	\max\limits_{1\leq i \leq n} \|\widehat{\bpi}_{i}\ind{k+1} - \bpi(X_{i})\|
	\lesssim a(1 + \beta_1)^{d+1} \beta_2 \frac{h_{k+1}}\varkappa.
\end{equation}

The proof of Theorem \ref{th1} goes by induction.
Let \( C \) be the hidden constant in \eqref{proj_ineq}.
Assume that on the \( k \)-th iteration
\[
	\max\limits_{1 \leq i \leq n} \|\widehat \bpi_i \ind k - \bpi(X_i)\| 
	\leq \frac{\Delta \vee (1 + \alpha) Ca(4^{d+1} + (2\sqrt{\alpha})^{d+1})}\varkappa h_k
\]
with probability at least \( 1 - (5k + 1)/n \).
Here \( \alpha \) is the constant from \eqref{a3}.
Lemma \ref{regularity}, Lemma \ref{main_lem} and Lemma \ref{proj} imply that, with 
probability at least \( 1 - (5(k+1) + 1)/n \), it holds
\[
	\max\limits_{1\leq i \leq n} \|\widehat{\bpi}_{i}\ind{k+1} - \bpi(X_{i})\|
	\leq aC(1 + \beta_1)^{d+1} \beta_2 \frac{h_{k+1}}\varkappa.
\]
We have to check that \( (1 + \beta_1)^{d+1} \beta_2 \leq (1+\alpha)\left( 4^{d+1} + 
(2\sqrt{\alpha})^{d+1}\right) \).
From the definition of \( \beta_1 \) and \( \beta_2 \), we have\\
\begin{align*}
	\beta_2
	&
	= 1 +  \Phi_{M, b, h, \varkappa, \Delta} \leq 1 + \alpha + o(1), \quad n 
	\rightarrow \infty,
	\\
	\beta_1
	&
	\leq \frac{M((Ca(1 + \alpha)(4^{d+1} + (2\sqrt{\alpha})^{d+1}) \vee \Delta) + b/h_K) 
	}\varkappa
	\\&
	\vee \frac{(1 + \Phi_{M, h, \varkappa, \Delta}) h_0}\varkappa \vee 
	\frac{(Ca(1 + \alpha)(4^{d+1} + (2\sqrt{\alpha})^{d+1}) \vee\Delta)^2 
	h_0^2}{\varkappa^2}
	\\&
	= \frac{Mb}{\varkappa h_K} + o(1),
	\quad n \rightarrow \infty
\end{align*}
Show that \eqref{a3} yields \( Mb/(\varkappa h_K) \leq \sqrt{\alpha} + o(1) \), as 
\( n \rightarrow \infty \).
Then, if \( n \) is sufficiently large, i.e. \( n \geq N_\Delta \), this will imply \( \beta_1 \leq 
\sqrt{\alpha} + 1 \), \( \beta_2 \leq 2(1 + \alpha) \).
Due to \eqref{a3},
\[
	\frac{M^2 b^3}{\varkappa h_K^4}
	\leq \frac{M^2 b^3}{\varkappa\left(\frac{D\log n}n \right)^{\frac4d} \vee
	\varkappa \left(\frac{DM^{2}\varkappa^2 \log n}n\right)^{\frac4{d+4}}}
	\leq \alpha.
\]
If \( M \leq h_K^2/\varkappa \) then
\[
	\frac{Mb}{\varkappa h_K} \leq \frac{h_K b}{\varkappa^2} \leq \frac{h_K}{\varkappa} = 
	o(1), \quad n\rightarrow\infty.
\]
Otherwise, we have
\[
	\left(\frac{Mb}{\varkappa h_K}\right)^2 \leq \frac{M\varkappa}{h_K^2} \cdot  
	\frac{M^2 b^2}{\varkappa^2 h_K^2} = \frac{M^3 b^2}{\varkappa h_K^4} \leq \alpha. 
\]
Thus, \( Mb/(\varkappa h_K) \leq \sqrt{\alpha} + o(1) \), \( n \rightarrow \infty \), and we have 
\( \beta_1 \leq \sqrt{\alpha} + 1 \), \( \beta_2 \leq 2(1 + \alpha) \) for \( n \geq N_\Delta \).
This yields
\begin{align*}
	&
	(1 + \beta_1)^{d+1} \beta_2
	\leq 2(1 + \alpha) \left( 2 + \sqrt{\alpha} \right)^{d+1}
	\\&
	\leq 2(1 + \alpha) \cdot 2^d \left( 2^{d+1} + \alpha^{(d+1)/2}\right)
	= (1 + \alpha) \left( 4^{d+1} + (2\sqrt{\alpha})^{d+1} \right).
\end{align*}
Thus, on the next iteration we have
\[
	\max\limits_{1\leq i \leq n} \|\widehat{\bpi}_{i}\ind{k+1} - \bpi(X_{i})\|
	\leq \frac{(1 + \alpha) Ca \left( 4^{d+1} + (2\sqrt{\alpha})^{d+1} \right)}\varkappa 
	h_{k+1}
\]
with probability at least \( 1 - (5(k+1) + 1)/n \).
The confidence level \( 1 - (5K + 4)/n \) in the claim of Theorem \ref{th1} follows from the 
fact that we do not recompute the projectors on the final step.

\subsection{Proof of Theorem \ref{th2}}

Fix any \( x \in \widehat \M \)
By definition of \( \widehat \M \), there exist \( i \) and \( u \in \B(0, 1) \), such that
\[
	x = \widehat{X}_{i} + h_K \widehat{\bpi}_{i}\ind K u.
\]
Lemma \ref{regularity}, Lemma \ref{main_lem}, Lemma \ref{proj}, and the union bound 
imply that, with probability at least \( 1 - (5K+4)/n \),
\[
	\max\limits_{1\leq i \leq n} d( \widehat{X}_{i}, \{X_i\}\oplus\T_{X_{i}}\M^{*} ) \lesssim 
	\frac{(1 + \Phi_{M, b, h_K, \varkappa})h_{K}^{2}}\varkappa + 
	\sqrt{\frac{D(h_{K}^4/\varkappa^2\vee M^{2}) \log h_{K}^{-1}}{nh_{K}^d}}
\]
and
\[
	\max\limits_{1\leq i\leq n} \|\widehat{\bpi}_{i}\ind K - \bpi(X_{i})\| \lesssim 
	\Psi_{M, b, h_K, \varkappa} \frac{h_K}\varkappa,
\]
where $\Phi_{M, b, h_K, \varkappa}$ and $\Psi_{M, b, h_K, \varkappa}$ are defined in 
\eqref{phi_psi}.
Recall that \( \proj{\M}x \) denotes the projection of \( x \) onto a closed set \( \M \).
Using the result Lemma \ref{proj}, we immediately obtain
\begin{align}
	\label{621}
	&\notag
	d(x, \{X_i\} \oplus \T_{X_{i}}\M^{*})
	= \inf\limits_{v\in\R^D} \|\widehat{X}_{i} + h_K \widehat{\bpi}_{i}\ind K u - X_{i} - 
	\bpi(X_{i}) v\|
	\\&
	\leq d(\widehat{X}_{i}, \{X_i\} \oplus \T_{X_{i}}\M^{*})
	\\&\notag
	+ \inf\limits_{v\in\R^D} \left\|\proj{\{X_i\} \oplus \T_{X_{i}}\M^{*}}{\widehat{X}_{i}} + 
	h_K \widehat{\bpi}_{i}\ind K u - X_{i} - \bpi(X_{i}) v\right\|
\end{align}
Since the vector \( \pi_{\{X_i\} \oplus \T_{X_{i}}\M^{*}}(\widehat{X}_{i}) - X_{i} \) belongs to 
\( \T_{X_{i}}\M^{*} \), we have
\[
	\pi_{\{X_i\} \oplus \T_{X_{i}}\M^{*}}(\widehat{X}_{i}) - X_{i} 
	= \bpi(X_i)\left( \pi_{\{X_i\} \oplus \T_{X_{i}}\M^{*}}(\widehat{X}_{i})  - X_{i} \right).
\]
Then, substituting $v + X_i - \pi_{\{X_i\} \oplus \T_{X_{i}}\M^{*}}(\widehat{X}_{i})$ by 
$\widetilde v$, we obtain that 
the last expression in \eqref{621} is equal to
\begin{align*}
	&
	d(\widehat{X}_{i}, \{X_i\} \oplus \T_{X_{i}}\M^{*})
	\\&
	+ \inf\limits_{v\in\R^D} \left\| h_K \widehat{\bpi}_{i}\ind K u - \bpi(X_{i}) \left(v + X_i - 
	\proj{\{X_i\} \oplus \T_{X_{i}}\M^{*}}{\widehat{X}_{i}} \right)\right\|
	\\&
	= d(\widehat{X}_{i}, \{X_i\} \oplus \T_{X_{i}}\M^{*}) + \inf\limits_{\widetilde v\in\R^D} 
	\|h_K \widehat{\bpi}_{i}\ind K u - \bpi(X_{i}) \widetilde v\|
	\\&
	\leq d(\widehat{X}_{i}, \{X_i\} \oplus \T_{X_{i}}\M^{*}) + \|h_K \widehat{\bpi}_{i}\ind K 
	u - h_K \bpi(X_{i}) u\|
	\\&
	\leq d(\widehat{X}_{i}, \{X_i\} \oplus \T_{X_{i}}\M^{*}) + h_K \|\widehat{\bpi}_{i}\ind K - 
	\bpi(X_{i})\|
	\\&
	\lesssim \frac{(1 + \Phi_{M, b, h_K, \varkappa} + \Psi_{M, b, h_K, 
	\varkappa})h_{K}^{2}}\varkappa + \sqrt{\frac{D(h_{K}^4/\varkappa^2\vee M^{2}) \log 
	h_{K}^{-1}}{nh_{K}^d}}.
\end{align*}
Next, note that, \( \|x - \widehat X_i\| \leq h_K \) and, due to Theorem \ref{th1}, we have
\begin{align*}
	\|\proj{\{X_i\} \oplus \T_{X_{i}}\M^{*}}x - X_i\|
	&
	\leq \|x - X_i\|
	\leq \|x - \widehat X_i\| + \|\widehat X_i - X_i\|
	\\&
	\lesssim \frac{Mb \vee M h_K \vee (1 + \Phi_{M, b, h_K, \varkappa}) 
	h_K^2}\varkappa + \sqrt{\frac{D(h_{K}^{2}\vee M^{2}) \log h_{K}^{-1}}{nh_{K}^d}}
	\\&
	\lesssim \left(\frac{Mb}\varkappa \vee h_K \right) + \sqrt{\frac{D(h_K^2 \vee M^2) \log 
	h_{K}^{-1}}{nh_K^d}}
	\lesssim \frac{Mb}\varkappa \vee h_K .
\end{align*}
The inequalities in the last line follow from the fact that $M < \varkappa$, $\Phi_{M, b, 
h_K, \varkappa} \lesssim \alpha + o(1)$, $n \rightarrow \infty$, and \( h_{K} \geq 
\left((D\log n/n)^{1/d} \vee (DM^{2}\varkappa^2 \log n/n)^{1/(d+4)}\right) \) due to the 
conditions of Theorem \ref{th1}.
Since \( \M^* \) is a \( \C^2 \)-manifold with a reach at least \( \varkappa \), it holds that
\[
	d(\proj{\{X_i\} \oplus \T_{X_{i}}\M^{*}}x, \M^*)
	\lesssim \frac{\|\proj{\{X_i\} \oplus \T_{X_{i}}\M^{*}}x - X_i\|^2}\varkappa
	\lesssim \frac{h_K^2}\varkappa \vee \frac{M^2 b^2}{\varkappa^3}.
\]	
Finally, we obtain
\begin{align*}
	d(x, \M^{*}) 
	&
	\leq d(x, \{X_i\} \oplus \T_{X_{i}}\M^{*}) + d(\proj{\{X_i\} \oplus \T_{X_{i}}\M^{*}}x, \M^*)
	\\&
	\lesssim \left( \frac{(1 + \Phi_{M, b, h_K, \varkappa} + \Psi_{M, b, h_K, 
	\varkappa})h_{K}^{2}}\varkappa \vee \frac{M^2 b^2}{\varkappa^3} \right) + 
	\sqrt{\frac{D(h_{K}^4/\varkappa^2\vee M^{2}) \log h_{K}^{-1}}{nh_{K}^d}}.
\end{align*}
Thus, \( \widehat{\M} \subseteq \M^{*} \oplus \B(0, r) \) with
\[
	r \lesssim \left( \frac{(1 + \Phi_{M, b, h_K, \varkappa} + \Psi_{M, b, h_K, 
	\varkappa})h_{K}^{2}}\varkappa \vee \frac{M^2 b^2}{\varkappa^3} \right) + 
	\sqrt{\frac{D(h_{K}^4/\varkappa^2\vee M^{2}) \log h_{K}^{-1}}{nh_{K}^d}}.
\]

It remains to prove that \( \M^{*} \subseteq \widehat{\M} \oplus \B(0, r) \) with the same \( 
r \).
Fix \( x\in\M^{*} \).
Note that there exist constants \( c_1 \) and \( r_0 \), such that
\begin{align*}
	&
	\p_{X}\left( X \in \B(x, r) \right) \geq p_0 \text{Vol}(\B(x, r) \cap \M^*) \geq  c_1 p_0 r^d, 
	\quad \forall r < r_0
	\\&
	\p_{X}\left( X \notin \B(x, r) \right) \leq 1 - c_1 p_0 r^d \leq e^{-c_1 p_0 r^d}, 
	\quad \forall r < r_0
\end{align*}
Let \( \N_\eps(\M^{*}) \) stand for an \( \eps \)-net of \( \M^{*} \).
It is known (see, for example, \cite{gppvw12a}, Lemma 3) that \( |\N_\eps(\M^{*})| 
\lesssim \eps^{-d} \).
Then
\begin{align*}
	&
	\p\left( \exists \, x \in \M^{*} : \forall \, i \quad X_{i} \notin \B(x, 2\eps) \right)
	\\&
	\leq \p\left( \exists \, x \in \N_\eps(\M^{*}) : \forall \, i \quad X_{i} \notin \B(x, \eps) \right)
	\\&
	\leq \sum\limits_{x \in \N_\eps(\M^{*})} \p\left( \forall \, i \quad X_{i} \notin \B(x, \eps) 
	\right)
	\lesssim \eps^{-d} e^{-c_1 p_0 n\eps^d}.
\end{align*}
This implies that with probability at least \( 1 - 1/n \)
\[
	\sup\limits_{x\in\M^{*}} \min\limits_{1 \leq i \leq n} \|x - X_{i}\| \lesssim \left( 
	\frac{\log n}n \right)^{\frac1d} \, .
\]
According to \citep[Theorem 4.18]{f59}, on the same event we have
\[
	\sup\limits_{x\in\M^{*}} \min\limits_{1 \leq i \leq n} d(x, \{X_i\} \oplus \T_{X_{i}}\M^{*}) 
	\lesssim \varkappa^{-1} \left( \frac{\log n}n \right)^{\frac2d}.
\]
Now, fix any \( x \in \M^{*} \).
Without loss of generality, assume that \( \min_{1 \leq i \leq n} d(x, \{X_i\} \oplus 
\T_{X_{i}}\M^{*}) \) is attained with \( i=1 \).
Let \( \proj{\{X_1\} \oplus \T_{X_1}\M^{*}}x \) be the projection of \( x \) onto the tangent 
plane \( \{X_1\} \oplus \T_{X_1}\M^{*} \).
It is clear that 
\[
	\left\| \proj{\{X_1\} \oplus\T_{X_1}\M^{*}}x - X_1 \right\| \lesssim \left( 
	\frac{\log n}n \right)^{\frac1d} \leq \frac{h_K}2.
\]
Then there exists \( u \in \B(0,1) \), such that
\begin{align*}
	\left\| \widehat{X}_1 + h_K \widehat{\bpi}_{1}\ind K u - \proj{\{X_1\} \oplus 
	\T_{X_1}\M^{*}}x \right\| 
	&
	\lesssim \left( \frac{(1 + \Phi_{M, b, h_K, \varkappa} + \Psi_{M, b, h_K, 
	\varkappa})h_{K}^{2}}\varkappa \vee \frac{M^2 b^2}{\varkappa^3} \right)
	\\&\quad
	+ \sqrt{\frac{D(h_{K}^4/\varkappa^2\vee M^{2}) \log h_{K}^{-1}}{nh_{K}^d}}.
\end{align*}
Thus,
\begin{align*}
	d(x, \widehat\M)
	&
	\lesssim \varkappa^{-1} \left( \frac{\log n}n \right)^{\frac2d} + \left( \frac{(1 + \Phi_{M, 
	b, h_K, \varkappa} + \Psi_{M, b, h_K, \varkappa})h_{K}^{2}}\varkappa \vee \frac{M^2 
	b^2}{\varkappa^3} \right) + 
	\sqrt{\frac{D(h_{K}^4/\varkappa^2\vee M^{2}) \log h_{K}^{-1}}{nh_{K}^d}}
	\\&
	\lesssim \left( \frac{(1 + \Phi_{M, b, h_K, \varkappa} + \Psi_{M, b, h_K, 
	\varkappa})h_{K}^{2}}\varkappa \vee \frac{M^2 b^2}{\varkappa^3} \right) + 
	\sqrt{\frac{D(h_{K}^4/\varkappa^2\vee M^{2}) \log h_{K}^{-1}}{nh_{K}^d}}.
\end{align*}

\subsection{Proof of Theorem \ref{lower_bound}}
\label{proof_lower_bound}

For the sake of convenience, the proof of Theorem \ref{lower_bound} is divided into two 
steps.
First, we use the following lemma to obtain the lower bound \eqref{lb_1}.
\begin{Lem}
	\label{lem_lower_bound_1}
	Suppose that the sample \( \Y_n = \{Y_1, \dots, Y_n\} \) is generated according to the 
	model \eqref{model}, where \( \M^* \in \mclass_\varkappa^d \), $\M^* \subseteq \B(0, 
	R)$ with $R = \sqrt{\varkappa^2 + M^2 b^2/\varkappa^4}$, \( \eps \)  satisfies 
	\eqref{a2}, and the density \( p(x) \) of \( X \) fulfils \eqref{a1'} with $L > 0$, $p_0 \leq 
	((d+1) \omega_{d+1}  R^d)^{-1}$, $p_1 \geq ((d+1) \omega_{d+1} R^d)^{-1}$, where 
	$\omega_{d + 1}$ is the volume of the unit Euclidean ball in $\R^{d+1}$.
	Then, for any estimate $\widehat \M$, it holds that
	\begin{equation*}
		\sup\limits_{\M^* \in \mclass_\varkappa^d} \E_{\M^*} 
		d_H(\widehat\M, \M^*)
		\geq \frac{M^2 b^2}{6 \varkappa^3}.
	\end{equation*}
\end{Lem}
The proof of Lemma \ref{lem_lower_bound_1} is moved to Appendix 
\ref{app_lower_bounds_1}.
The construction used in the proof of Lemma \ref{lem_lower_bound_1} is extremely 
simple.
We show that if $(\eps \cond X)$ is supported on a tangent space $\T_X \M_0$, where 
$\M_0$ is a $d$-dimensional sphere of radius $\varkappa$, then a statistician cannot 
distinguish between $\M_0$ and a sphere $\M_1$ with greater radius.

Second, we prove the lower bound \eqref{lb_2} using Lemma \ref{lem_lower_bound_2} 
below.
The proof of Lemma \ref{lem_lower_bound_2} is based on application of
\cite[Theorem 2.5]{t09} to a family of smooth manifolds with small bumps at different 
points.
\begin{Lem}
	\label{lem_lower_bound_2}
	Suppose that the sample \( \Y_n = \{Y_1, \dots, Y_n\} \) is generated according to the 
	model \eqref{model}, where \( \M^* \in \mclass_\varkappa^d \), the density \( p(x) \) 
	of \( X \) fulfils \eqref{a1'} with sufficiently large $p_1$, sufficiently small $p_0 > 0$, 
	and $L \geq 4p_1 / 3$.
	Let the noise \( \eps \)  satisfy \eqref{a2} with
	\[
		b \gtrsim \left((\log n/n)^{1/d} \vee (M^{2}\varkappa^2 \log n/n)^{1/(d+4)}\right),
	\]
	where the hidden constant is large enough.
	Then, for any estimate $\widehat \M$, if \( n \) is sufficiently large and \( M\varkappa 
	\gtrsim (\log n/n)^{2/d} \), it holds that
	\begin{equation*}
		\sup\limits_{\M^* \in \mclass_\varkappa^d} \E_{\M^*} 
		d_H(\widehat\M, \M^*)
		\gtrsim \varkappa^{-1} \left(\frac{M^2\varkappa^2 \log n}n \right)^{\frac 2{d+4}}.
	\end{equation*}
\end{Lem}
The proof of Lemma \ref{lem_lower_bound_2} is moved to Appendix 
\ref{app_lower_bounds_2}.
The claim of Theorem \ref{lower_bound} follows from Lemma \ref{lem_lower_bound_1} 
and Lemma \ref{lem_lower_bound_2}.


\acks{The authors are grateful to the action editor and three anonymous reviewers for valuable suggestions and remarks. This work was partly supported by the German Ministry for Education and 
Research as BIFOLD. It was partly carried out within the framework of the HSE University 
Basic Research Program. Results of Section 5 have been obtained under support of the 
RSF grant No. 19-71-30020. Nikita Puchkin is a Young Russian Mathematics award 
winner and would like to thank its sponsors and jury.}



\appendix

\section{Proof of Lemma \ref{geom}}

We have
\begin{align*}
	&
	\| X_j - X_i - \bpi_i(Y_j - Y_i) \|
	\\&
	\leq \| X_j - X_i - \bpi_i(X_j - X_i)\| + \|\bpi_i(\eps_j - \eps_i)\|
	\\&
	\leq \| X_j - X_i - \bpi(X_i)(X_j - X_i)\|
	\\&
	+ \|\bpi_i - \bpi(X_i)\| \|X_j - X_i\| + \|\bpi_i(\eps_j - \eps_i)\| \, .
\end{align*}
According to \citep[Theorem 4.18]{f59}, 
\[
	 \| X_j - X_i - \bpi(X_i)(X_j - X_i)\| \leq \frac{\|X_j - X_i\|^2}{2\varkappa} \, .
\]
Consider the term $\|\bpi_i(\eps_j - \eps_i)\|$.
It holds
\begin{align}
	&\notag
	\|\bpi_i(\eps_j - \eps_i)\|
	\leq \|(\bpi_i - \bpi(X_i)(\eps_j - \eps_i)\| + \|\bpi(X_i)(\eps_j - \eps_i)\|
	\\&\label{l1e1}
	\leq \frac{2M\Delta h}\varkappa + \|\bpi(X_i)\eps_i\| + \|\bpi(X_j)\eps_j\| + \|(\bpi(X_i) 
	- 
	\bpi(X_j))\eps_j\|
	\\&\notag
	\leq \frac{2M\Delta h}\varkappa + \frac{2M b}\varkappa + \frac{3M\|X_j - 
	X_i\|}\varkappa,
\end{align}
where in the last inequality we used the condition \eqref{a2} and applied Proposition 
\ref{aux}.
Taking into account that $\|\bpi_i - \bpi(X_i)\| \leq \frac{\Delta h}\varkappa$, we conclude
\begin{align}
	\label{aa1}
	&
	\| X_j - X_i - \bpi_i(Y_j - Y_i) \|
	\leq \frac{\|X_j - X_i\|^2}{2\varkappa} + \frac{\Delta h \|X_j - X_i\|}\varkappa
	\\&\notag
	+ \frac{2M(\Delta h + b)}\varkappa + \frac{3M\|X_i - X_j\|}\varkappa.
\end{align}
Using the triangle inequality
\[
	\|X_i - X_j\| - \|\bpi_i(Y_i - Y_j)\| \leq \| X_j - X_i - \bpi_i(Y_j - Y_i) \|
\]
and solving the quadratic inequality
\begin{align*}
	&
	\|X_i - X_j\| - \|\bpi_i(Y_i - Y_j)\|
	\leq \frac{\|X_j - X_i\|^2}{2\varkappa}
	\\&
	+ \frac{(3M + \Delta h) \|X_j - X_i\|}\varkappa
	+ \frac{2M(\Delta h + b)}\varkappa
\end{align*}
with respect to $\|X_i - X_j\|$, we obtain
\[
	\| X_i - X_j \|
	\leq \frac{ \|\bpi_i(Y_i - Y_j)\| + 2M(\Delta h + b) /\varkappa}{ 1 - (\Delta h + 
	3M)/\varkappa}
	\leq 2\|\bpi_i(Y_i - Y_j)\| + \frac{4M(\Delta h + b)}\varkappa.
\]
Here we used the fact that, due to condition of the lemma,
\[
	\Delta h + 3M \leq \frac{\varkappa}4 + \frac{3\varkappa}{16} < \frac{\varkappa}2.
\]

On the other hand, from \eqref{aa1} we have
\begin{align*}
	&
	\|X_i - X_j\| - \|\bpi_i(Y_i - Y_j)\|
	\geq -\frac{\|X_j - X_i\|^2}{2\varkappa}
	\\&\notag
	- \frac{2M(\Delta h + b)}\varkappa - \frac{(3M + \Delta h)\|X_i - X_j\|}\varkappa.
\end{align*}
If
\[
	\frac{\|X_i - X_j\|^2}{2\varkappa} + \frac{3M + \Delta h}\varkappa \|X_i - 
	X_j \| + \frac{2M(\Delta h + b)}\varkappa \leq \frac{ \|\bpi_i(Y_j - Y_i)\| }2,
\]
then $\|X_i - X_j\| \geq 0.5 \|\bpi_i(Y_j - Y_i)\|$.
Otherwise, it holds
\begin{align*}
	&
	\|X_i - X_j\|
	\geq -(3M+\Delta h)
	\\&
	+ \frac1\varkappa \sqrt{ \left( \frac{\Delta h + 
	3M}\varkappa\right)^2 + \frac{\|\bpi_i(Y_i - Y_j)\| - 4M(\Delta h + b) 
	/\varkappa}\varkappa }.
\end{align*}
Introduce a  function $g(t) = \sqrt{a^2 + t} - a$, $a > 0$, $t \geq -a^2$.
The function $g(t)$ is concave, increasing and $g(0) = 0$.
Therefore, for any $t_0$ and any $t \in [0, t_0]$ it holds
\[
	 g(t) \geq g(t_0) \frac t{t_0}.
\]
Taking $a = (\Delta h + 3M)/\varkappa$ and $t_0 = 1 - 4M(\Delta h + b)/\varkappa^2$, we 
immediately obtain
\begin{align*}
	\|X_i - X_j\|
	&\geq \left( -\frac{3M+\Delta h}\varkappa + \sqrt{ \left( \frac{3M+\Delta h}\varkappa 
	\right)^2 + \frac{\varkappa - 4M\Delta h/\varkappa}\varkappa } \right) 
	\\&
	\cdot \left( \|\bpi_i(Y_i - Y_j)\| - \frac{4M(\Delta h + b)}\varkappa \right).
\end{align*}
Now it is easy to see that, if $M \leq \varkappa/16$, $\Delta h \leq 
\varkappa/4$, and $M(\Delta + b/h) \leq \varkappa/4$ then $3M + \Delta h < 
\varkappa/2 < \varkappa$ and
\[
	1 - \frac{4M(\Delta h + b)}{\varkappa^2}
	\geq \frac34 \frac{(3M + \Delta h)^2}{\varkappa^2}.
\]
The last inequality yields
\[
	-\frac{3M+\Delta h}\varkappa + \sqrt{ \left( \frac{3M+\Delta h}\varkappa \right)^2 + 
	\frac{\varkappa - 4M\Delta h/\varkappa}\varkappa }
	\geq \frac12,
\]
which, in its turn, implies
\[
	\|X_i - X_j\| \geq \frac{\|\bpi_i(Y_i - Y_j)\|}2 - 2M(\Delta h + b).
\]

\section{Proof of Lemma \ref{regularity}}

Show that for any $\bpi_i$, such that $\|\bpi_i - \bpi(X_i)\| \leq \Delta h/\varkappa$, it 
holds
\[
	\E\ind{-i} w_{i1}(\bpi_i) \geq C_1 h^d.
\]
Here and further in this paper, $\E\ind{-i}(\cdot) \equiv \E(\cdot \cond (X_i, Y_i))$.
Due to Lemma \ref{geom}, we have
\[
	\|\bpi_i(Y_i - Y_1)\| \leq 2 \|X_i - X_1\| + \frac{4M(\Delta h + b)}\varkappa,
\]
which yields
\[
	\|\bpi_i(Y_i - Y_1)\|^2 \leq 8 \|X_i - X_1\|^2 + \frac{32M^2(\Delta h + b)^2}\varkappa,
\]
\begin{align*}
	&
	\E\ind{-i} w_{i1}(\bpi_i)
	= \E\ind{-i} e^{-\frac{\|\bpi_i(Y_i - Y_1)\|^2}{h^2}} \1\left( \|Y_i - Y_1\| \leq \tau \right)
	\\&
	\geq \E\ind{-i} e^{-32M^2(\Delta + \alpha)^2/\varkappa^2} e^{-\frac{8\|X_i - 
	X_1\|^2}{h^2}} \1\left( \|X_i - X_1\| \leq 
	\tau - 2M \right)
	\\&
	\geq e^{-2} \int\limits_{\M^*\cap\B(X_i, \tau-2M)} e^{-\frac{8\|X_i - 
	x\|^2}{h^2}} p(x)dW(x)
	\\&
	\geq e^{-2}p_0 \int\limits_{\B(X_i, h)} e^{-\frac{8\|X_i - x\|^2}{h^2}} dW(x)
	\\&
	= \frac{p_0}{e^2} \int\limits_{\Exp^{-1}(\B(X_i, h))} e^{-\frac{8\|\Exp_{X_i}(p) - 
	\Exp_{X_i}(0)\|^2}{h^2}} \sqrt{\det g(v)} dv.
\end{align*}
Here we used the fact that, due to conditions of the lemma, $M(\Delta + b/h) \leq 
\varkappa/4$ and, also, $\tau - 2M \geq 0.5 \tau \geq h_0 \geq h$, if $h_0$ is chosen 
sufficiently small.
Next, due to \cite[Lemma 1]{al19}, 
\[
	\| \Exp_{X_i}(v) - \Exp_{X_i}(0) - v \| \leq C_{\perp} \|v\|^2 \leq C_\perp\varkappa \|v\|.
\]
It also holds $\det g(v) \geq \frac12$ for any $v \in \Exp^{-1}(\B(X_i, h))$.
Then there exists a constant $C'$, depending on $d$, $p_0$ and $\varkappa$, such that
\[
	\E\ind{-i} w_{i1}(\bpi_i) \geq 2C' h^d.
\]
Now, consider the sum
\[
	\sum\limits_{j=1}^n w_{ij}(\bpi_i)
\]
Given $Y_i$, the weights $w_{ij}(\bpi_i)$ are conditionally independent and identically 
distributed.
The Bernstein's inequality implies that
\[
	\p\ind{-i} \left( \sum\limits_{j=1}^n w_{ij}(\bpi_i) \leq C'h^d \right) \leq e^{-\frac{C'^2 
	n h^{2d}}{2\sigma^2 + 2C'h^d/3}} \leq e^{-C'' n h^d} \, ,
\]
and $e^{-C'' n h^d} \leq n^{-1}$ if $h \gtrsim \left(\frac{\log n}n \right)^{1/d}$ (with a 
sufficiently large hidden constant).
Therefore, with probability at least $1 - n^{-2}$, it holds
\[
	\sum\limits_{j=1}^n w_{ij}(\bpi_i) \geq C'h^d.
\]

\section{Proof of Lemma \ref{main_lem}}

\label{app_c}

Fix any $i$ from $1$ to $n$ and denote $\E\ind{-i}(\cdot) \equiv \E\left( \cdot | (X_i, Y_i) 
\right)$ and $\p\ind{-i}(\cdot) \equiv \p\left(\cdot | (X_i, Y_i) \right)$.
Also, let $\pset_i(\Delta h/\varkappa)$ be a set of projectors $\bpi$ onto 
$d$-dimensional space, such that $\|\bpi - \bpi(X_i)|| \leq \Delta h/\varkappa$.
First, we study the supremum of the empirical process
\[
	 \sup\limits_{\bpi_i\in\pset_i(\Delta h/\varkappa)} \left\| \sum\limits_{j=1}^n 
	 w_{ij}(\bpi_i) (Y_j - X_i) \right\|
	 = \sup\limits_{\substack{u \in \B(0, 1)\\\bpi_i\in\pset_i(\Delta h/\varkappa)}} 
	 \sum\limits_{j=1}^n 
	 w_{ij}(\bpi_i) u^T(Y_j - X_i) \, .
\]
The rest of the proof can be summarized as follows.
First, we fix $u \in \B(0, 1)$ and $\bpi_i\in\pset_i(\Delta h/\varkappa)$ and bound the 
supremum of 
the expectation
\[
	\sup\limits_{\substack{u \in \B(0, 1)\\\bpi_i\in\pset_i(\Delta h/\varkappa)}}  \E\ind{-i} 
	\sum\limits_{j=1}^n w_{ij}(\bpi_i) u^T(Y_j - X_i) \, .
\]
Then we provide uniform bounds on
\[
	\E\ind{-i} \sup\limits_{\substack{u \in \B(0, 1)\\\bpi_i\in\pset_i(\Delta h/\varkappa)}} 
	\left( 
	\sum\limits_{j=1}^n w_{ij}(\bpi_i) u^T(Y_j - X_i) - \E\ind{-i} \sum\limits_{j=1}^n 
	w_{ij}(\bpi_i) u^T(Y_j - X_i) \right).
\]
Finally, we derive large deviation results for
\[
	\sup\limits_{\substack{u \in \B(0, 1)\\\bpi_i\in\pset_i(\Delta h/\varkappa)}} 
	\sum\limits_{j=1}^n 
	w_{ij}(\bpi_i) u^T(Y_j - X_i) - \E\ind{-i} \sup\limits_{\substack{u \in \B(0, 
	1)\\\bpi_i\in\pset_i(\Delta h/\varkappa)}} \sum\limits_{j=1}^n w_{ij}(\bpi_i) u^T(Y_j - 
	X_i).
\]

As it was said earlier, we start with bounds on the expectation.
The rigorous result is given in the next proposition.
\begin{Prop}
	\label{lem_expectation}
	Under conditions of Theorem \ref{th1} and Lemma \ref{main_lem}, for any $u \in 
	\B(0, 1)$ and $\bpi_i\in\pset_i(\Delta h/\varkappa)$, it holds
	\begin{align}
		\label{expbound}
		\tag{a}
		&
		\E\ind{-i} \sum\limits_{j=1}^n w_{ij}(\bpi_i) u^T(X_j - X_i) \lesssim \left(M(\Delta + 
		b/h) \vee h \vee \frac{\Delta^2 h^2}\varkappa \right) \frac{h^{d+1}}\varkappa, 
		\\&
		\label{expbound2}
		\tag{b}
		\E\ind{-i} \sum\limits_{j=1}^n w_{ij}(\bpi_i) u^T(\bid - \bpi(X_i))(X_j - X_i) \lesssim 
		\frac{nh^{d + 2}}\varkappa, 
		\\&
		\label{expbound3}
		\tag{c}
		\E\ind{-i} \sum\limits_{j=1}^n w_{ij}(\bpi_i) u^T\eps_j
		\lesssim \Phi_{M, b, h, \varkappa, \Delta} \frac{nh^{d+2}}\varkappa,
	\end{align}
	where
	\[
		\Phi_{M, b, h, \varkappa, \Delta} =
		\frac{M^3(1 + \Delta + b/h)^2}{h^2 \varkappa}
		+ \frac{M^2(\Delta + b/h + \sqrt{\log h^{-1}})}{\varkappa h}
		+ \frac{(1 + \Delta^4) M h^{2}}{\varkappa^3},
	\]
	and the hidden constants do not depend on $\Delta$.
\end{Prop}
The proof of Proposition \ref{lem_expectation} relies on Taylor's expansion but it is quite 
technical.
Therefore, it is moved to Appendix \ref{app_e}.
We continue with a uniform bound on the expectation
\[
	\E \sup\limits_{\substack{u \in \B(0, 1),\\ \| \bpi_i - \bpi(X_i) \| \leq \Delta h/\varkappa}} 
	\sum\limits_{j=1}^n \left( w_{ij}(\bpi_i) u^T Y_j - \E w_{ij}(\bpi_i) u^T Y_j \right).
\]
Introduce the class of functions
\begin{align*}
	\F_i = \big\{
	f(y) = \K\left(\frac{\| \bpi_i(Y_i - y) \|^2}{h^2} \right) \1( \|Y_i - y \| \leq \tau) u^T(y - X_i) :
	\\
	\| \bpi_i - \bpi(X_i) \| \leq \Delta h/\varkappa, \, Y_i \in \B(X_i, M), \, u \in \B(0, 1) \big\}.
\end{align*}
We use the same trick as in \citep[Section 4]{gk06}.
Note that the class
\begin{align}
	\label{f1class}
	\F_i\ind1 = \big\{ 
	&
	f_1(y) = \| \bpi_i(Y_i - y ) \| :
	\\&\notag
	\| \bpi_i - \bpi(X_i) \| \leq \Delta h/\varkappa, \, Y_i \in \B(X_i, M)\big\}
\end{align}
is VC subgraph, because the stripe $\{ y : \| \bpi(Y_i - y ) \| \leq t \}$ is an intersection of a 
finite number of halfspaces.
According to \citep[Theorem 2.6.18 (viii)]{vdvw96}, the class
\begin{align*}
	\widetilde\F_i\ind1 = \bigg\{ 
	&
	f_1(y) = \K\left( \frac{\| \bpi_i(Y_i - y ) \|^2}{h^2} \right) :
	\\&
	\| \bpi_i - \bpi(X_i) \| \leq \Delta h/\varkappa, \, Y_i \in \B(X_i, M)\bigg\}
\end{align*}
is also VC subgraph, since $\K(\cdot)$ monotonously decreases.
The class of balls
\begin{align}
	\label{f2class}
	\F_i\ind2 = \big\{ 
	&
	f_2(y) = \1( \|Y_i - y \| \leq \tau) : Y_i \in \B(X_i, M)\big\}
\end{align}
and the class of hyperplanes
\[
	\F_i\ind3 = \big\{ f_3(y) = u^T(y - X_i) : u \in \B(0, 1) \big\}
\]
are VC subgraph.
The functions from the classes $\widetilde\F_i\ind1, \F_i\ind2$ and $\F_i\ind3$ are 
bounded by $1, 1$ and $R+M$ respectively.
Then there exist constants $\A$ and $\nu$, depending only on the VC characteristics of 
the classes $\widetilde\F_i\ind1, \F_i\ind2$ and $\F_i\ind3$, such that
\[
	\N(\F_i, L_2(\p\ind{-i}_n), \delta) \leq \left( \frac\A\delta \right)^\nu,
\] 
where $\N(\F_i, L_2(\p\ind{-i}_n), \delta)$ is the $\delta$-covering number of $\F_i$ with 
respect to the $L_2(\p\ind{-i}_n)$ metric.
Theorem 6 in \cite{rxh11} (see also \cite{s98}) implies that we can take $\nu \lesssim 
Dd$ and $\A$ to be an absolute constant, which does not depend on $D, d$ or 
$\varkappa$.
Corollary 2.2 from \cite{gg02} implies
\begin{equation}
	\label{vcbound}
	\E\ind{-i} \sup\limits_{f \in \F_i} \sum\limits_{j=1}^n \left( f(Y_j) - \E\ind{-i} f(Y_j) \right)
	\leq \Rcal \sigma \left( \sqrt{D n \log \frac\A\sigma} \vee D \log 
	\frac\A\sigma \right),
\end{equation}
with an absolute constant $\Rcal$ and $\sigma^2 \geq \sup_{f\in \F} \text{Var} f(Y_1)$.
Lemma \ref{geom} yields
\[
	\|X_{i} - X_{j} \|^2 \leq 8\|\bpi_{i}(Y_{i} - Y_{j})\|^2 + 
	\frac{32M^2(\Delta+b/h)^2h^2}{\varkappa^2}.
\]
Using this, we can derive
\begin{align*}
	&
	\E\ind{-i} e^{-\frac{2\|\bpi_i(Y_j - Y_i)\|^2}{h^2}} (u^T(Y_j - X_i))^2
	\\&
	\leq \E\ind{-i} e^{\frac{8M^2(\Delta+\alpha)^2}{\varkappa^2}} e^{-\frac{\|X_j - 
	X_i\|^2}{4h^2}} \|Y_j - 
	X_i\|^2
	\\&
	\leq 2 e^{1/2} \E\ind{-i} e^{-\frac{\|X_j - X_i\|^2}{4h^2}} \|X_j - X_i\|^2
	+ 2 e^{1/2} \E\ind{-i} e^{-\frac{\|X_j - X_i\|^2}{2h^2}} \|\eps_j\|^2
	\\&
	\leq 2 e^{1/2}  \E\ind{-i} e^{-\frac{\|X_j - X_i\|^2}{4h^2}} \|X_j - X_i\|^2
	+ 2 e^{1/2}  \E\ind{-i} e^{-\frac{\|X_j - X_i\|^2}{4h^2}} M^2.
\end{align*}
Here we used the fact that $4M(\Delta + \alpha) \leq \varkappa$.
Next, due to Lemma \ref{p1a4}, there exist absolute constants $B_1$ and $B_2$, such 
that
\begin{align*}
	&
	\E\ind{-i} e^{-\frac{\|X_j - X_i\|^2}{4h^2}} \|X_j - X_i\|^2 \leq B_1 h^{d+2} \, ,
	\\&
	\E\ind{-i} e^{-\frac{\|X_j - X_i\|^2}{4h^2}} \leq B_2 h^d \, .
\end{align*}
Therefore, we can take $\sigma^2 = B^2 (h^2\vee M^2) h^d$ with an absolute constant 
$B$.
Thus, there exists a constant $C_{R, M, d}$, depending on $R, M$ and $d$ only (but not 
on $\Delta$), such that
\begin{align*}
	&
	\E\ind{-i} \sup\limits_{f \in \F_i} \sum\limits_{j=1}^n \left( f(Y_j) - \E\ind{-i} f(Y_j) \right)
	\\&
	\leq \Rcal B \sqrt{D(h^2\vee M^2) n h^d \log\frac\A{B^2 (h^2\vee M^2)h^d}}.
\end{align*}
Finally, we use the Talagrand's concentration inequality \citep{t96} and obtain bounds on 
large deviations of
\[
	\sup\limits_{\substack{u \in \B(0, 1),\\ \|\bpi_i - \bpi(X_i)\| \leq \Delta h/\varkappa}} 
	\sum\limits_{j=1}^n w_{ij}(\bpi_i) u^T (Y_j - X_i).
\]
More precisely, we use the version of Talagrand's inequality from \citep{b02}, where a  
deviation bound with nice constants was derived.
Denote
\[
	Z_i = \sup\limits_{f \in \F_i} \sum\limits_{j=1}^n \left( f(Y_j) - \E f(Y_j) \right).
\]
Then \citep[Theorem 2.3]{b02} claims that, on an event with probability $1 - n^{-2}$, 
it holds that
\[
	Z_i \leq \E Z_i + \sqrt{4 v \log n} + \frac{2 \log n}3,
\]
with $v = n\sigma^2 + 2 \E Z_i$ and the same $\sigma$ as in \eqref{vcbound}.
This, together with \eqref{expbound} and \eqref{vcbound}, yields
\begin{align*}
	&
	\sup\limits_{\substack{u \in \B(0, 1),\\ \|\bpi_i - \bpi(X_i)\| \leq \Delta h/\varkappa}} 
	\sum\limits_{j=1}^n w_{ij}(\bpi_i) u^T(Y_j - X_i) 
	\lesssim \left(M(\Delta + b/h) \vee h \vee \frac{\Delta^2 h^2}\varkappa \right) 
	\frac{h^{d+1}}\varkappa
	\\&
	+ \Phi_{M, b, h, \varkappa, \Delta} 
	\frac{nh^{d+2}}\varkappa + \sqrt{D (h^{2}\vee M^{2}) nh^d \log n}
\end{align*}
on an event with probability at least $1 - n^{-2}$.
The union bound implies that, on an event with probability at least $1 - n^{-1}$, it holds that
\begin{align*}
	&
	\max\limits_{1 \leq i \leq n} \left\| \sum\limits_{j=1}^n w_{ij}(\bpi_i)(Y_j - X_i)\right\|
	\lesssim \left(M(\Delta + b/h) \vee h \vee \frac{\Delta^2 h^2}\varkappa \right) 
	\frac{h^{d+1}}\varkappa
	\\&
	+ \Phi_{M, b, h, \varkappa, \Delta} 
	\frac{nh^{d+2}}\varkappa + \sqrt{D(h^{2}\vee M^{2}) nh^d \log n}.
\end{align*}
The bound
\begin{align*}
	&
	\max\limits_{1\leq i \leq n} \left\|\sum\limits_{j=1}^n w_{ij}(\bpi_i)\big((\bid - 
	\bpi(X_i))(Y_j 
	- X_i) \big)\right\| 
	\\&
	\lesssim (1 + \Phi_{M,b, h, \varkappa,\Delta})nh^{d+ 2}/\varkappa + 
	\sqrt{D(h^4/\varkappa^2 \vee M^2) nh^d \log n}
\end{align*}
is proven in a completely similar way.
Proposition \ref{lem_expectation} yields
\[
	\sup\limits_{\substack{u \in \B(0, 1),\\ \|\bpi_i - \bpi(X_i)\| \leq \Delta h/\varkappa}} 
	\E\ind{-i} \sum\limits_{j=1}^n w_{ij}(\bpi_i) u^T (\bid - \bpi(X_i))(Y_j - X_i) \lesssim 
	(1 + \Phi_{M,b, h, \varkappa,\Delta}) \frac{n h^{d+2}}\varkappa.
\]
Again, applying the VC subgraph argument and using \citep[Corollary 2.2]{gg02}, we 
obtain
\begin{align*}
	&
	\E\ind{-i} \sup\limits_{\substack{u \in \B(0, 1),\\ \|\bpi_i - \bpi(X_i)\| \leq \Delta 
	h/\varkappa}} \sum\limits_{j=1}^n w_{ij}(\bpi_i) u^T (\bid - \bpi(X_i))(Y_j - X_i)
	\\&
	\lesssim (1 + \Phi_{M,b, h, \varkappa,\Delta})\frac{n h^{d+2}}\varkappa + \Rcal 
	\sigma' 
	\left( \sqrt{D n \log \frac\A{\sigma'}} \vee D \log \frac\A{\sigma'} \right).
\end{align*}
The only difference is that we can take $(\sigma')^2 \asymp 
h^d(h^4/\varkappa^2 \vee M^2)$ in this case.
The reason for that is \cite[Theorem 4.18]{f59}, which implies
\[
	\left\| (\bid - \bpi(X_i))(X_j - X_i) \right\| \leq \frac{\|X_j - X_i\|^2}{2\varkappa},
\]
and then
\begin{align*}
	&
	\E\ind{-i} w_{ij}^2(\bpi_i) \left( u^T(\bid - \bpi(X_i))(Y_j - X_i) \right)^2
	\\&
	\leq \E\ind{-i} w_{ij}^2(\bpi_i) \| (\bid - \bpi(X_i))(Y_j - X_i)\|^2
	\\&
	\leq 2 \E\ind{-i} w_{ij}^2(\bpi_i) \| (\bid - \bpi(X_i))(X_j - X_i)\|^2
	\\&
	+ 2 \E\ind{-i} w_{ij}^2(\bpi_i) \| \eps_j \|^2
	\lesssim n h^d \left( h^4/\varkappa^{2}  \vee M^2 \right).
\end{align*}
Thus,
\begin{align*}
	&
	\E\ind{-i} \sup\limits_{\substack{u \in \B(0, 1),\\ \|\bpi_i - \bpi(X_i)\| \leq \Delta 
	h/\varkappa}} \sum\limits_{j=1}^n w_{ij}(\bpi_i) u^T (\bid - \bpi(X_i))(Y_j - X_i)
	\\&
	\lesssim (1 + \Phi_{M,b, h, \varkappa,\Delta})\frac{n h^{d+2}}\varkappa + \sqrt{D nh^d 
	(h^4/\varkappa^2 \vee M^2) \log h^{-1}}.
\end{align*}
Finally, applying the Talagrand's concentration inequality \citep[Theorem 2.3]{b02}), we 
obtain that, for a fixed $i$, with probability at least $1 - n^{-2}$, 
\begin{align*}
	&
	\left\| \sum\limits_{j=1}^n w_{ij}(\bpi_i)(\bid - \bpi(X_i))(Y_j - X_i)\right\|
	\\&
	\lesssim \left((1 + \Phi_{M,b, h, \varkappa,\Delta})nh^{d+2}/\varkappa + 
	\sqrt{D(h^4/\varkappa^2 \vee M^2)nh^d \log n} \right).
\end{align*}
Applying the union bound, we conclude
\begin{align*}
	&
	\max\limits_{1 \leq i \leq n} \left\| \sum\limits_{j=1}^n w_{ij}(\bpi_i)(\bid - \bpi(X_i))(Y_j 
	- X_i)\right\|
	\\&
	\lesssim (1 + \Phi_{M,b, h, \varkappa,\Delta})nh^{d+2}/\varkappa + 
	\sqrt{D(h^4/\varkappa^2 \vee M^2)nh^d \log n}.
\end{align*}

\section{Proof of Lemma \ref{proj}}
\label{app_d}

Throughout the proof of Lemma \ref{proj}, we work on the event $\Omega_k$, on which 
\eqref{beta} holds.
Consider
\[
	\widehat \bsigma_i\ind k = \sum\limits_{j=1}^n (\widehat X_j\ind k - \widehat X_i\ind 
	k)(\widehat X_j\ind k - \widehat X_i\ind k)^T \1\left( \|\widehat X_j\ind k - \widehat 
	X_i\ind k\| \leq \gamma h_k \right) \, .
\]
Denote $v_{ij} = \1\left( \|\widehat X_j\ind k - \widehat X_i\ind k\| \leq \gamma h_k \right)$.
Let
\begin{align*}
	Z_{ij}&= \proj{\{X_i\}\oplus\T_{X_i}\M^*}{X_j}, \quad 1 \leq j \leq n,
	\\
	\widehat Z_{ij}\ind k &= \proj{\{X_i\}\oplus\T_{X_i}\M^*}{\widehat X_j\ind k}, \quad 1 
	\leq j \leq n,
\end{align*}
and introduce a matrix
\[
	\widehat\bxi_i\ind k = \sum\limits_{j=1}^n v_{ij} (\widehat Z_{ij}\ind k - \widehat 
	Z_{ii}\ind k)(\widehat Z_{ij}\ind k - \widehat Z_{ii}\ind k)^T \, .
\]

From the conditions of Lemma \ref{proj}, we have
\[
	\max\limits_{1 \leq j \leq n} \|\widehat Z_{ij}\ind k - \widehat X_j\ind k\| \leq \beta_2 
	\left( \frac{h_k^2}\varkappa + \sqrt{\frac{h_k^4/\varkappa^2 \vee M^2}{nh_k^d} D \log 
	n} \right) \leq \frac{2\beta_2 h_k^2}\varkappa.
\]
This yields
\begin{align*}
	&
	\|\widehat \bsigma_i\ind k - \widehat \bxi_i\ind k \|
	\\&
	= \sup\limits_{u \in \B(0,1)} \left| \sum\limits_{j=1}^n v_{ij} \left[ 
	(u^T(\widehat X_j\ind k - \widehat X_i\ind k))^2 - (u^T(\widehat Z_{ij}\ind k - \widehat 
	Z_{ii}\ind k))^2 \right] \right|
	\\&
	\leq \sum\limits_{j=1}^n v_{ij} \left( \|\widehat X_j\ind k - \widehat X_i\ind k \| + 
	\|\widehat Z_{ij}\ind k - \widehat Z_{ii}\ind k \| \right)
	\\&
	\qquad\cdot\left( \|\widehat X_j\ind k - \widehat Z_{ij}\ind k\| + \|\widehat X_i\ind k - 
	\widehat Z_{ii}\ind k\| \right).
\end{align*}
Since $T_{X_i}\M^*$ is a convex set, then
\begin{align*}
	\|\widehat Z_{ij}\ind k - \widehat Z_{ii}\ind k\|
	&
	= \left\| \proj{\{X_i\}\oplus\T_{X_i}\M^*}{\widehat X_j\ind k} - 
	\proj{\{X_i\}\oplus\T_{X_i}\M^*}{\widehat X_i\ind k} \right\|
	\\&
	\leq \|\widehat X_j\ind k - \widehat X_i\ind k\|.
\end{align*}
Thus,
\[
	\|\widehat \bsigma_i\ind k - \widehat \bxi_i\ind k \|
	\leq \frac{8\beta_2 h_k^2}\varkappa \sum\limits_{j=1}^n v_{ij} \|\widehat X_j\ind k - 
	\widehat X_i\ind k\|
	\leq \frac{8\gamma \beta_2 h_k^3}\varkappa \sum\limits_{j=1}^n v_{ij}.
\]
Next, we are going to prove that, with probability at least $1 - n^{-2}$,
\begin{align}
	\label{e1}
	\sum\limits_{j=1}^n v_{ij}
	&
	\leq \sum\limits_{j=1}^n \1 \left( \|X_j - X_i\| \leq (\gamma + 
	4\beta_1) h_k \right)
	\\&\notag
	\leq 2C' n(\gamma + 4\beta_1)^d h_k^d.
\end{align}
The first inequality follows from the fact that $\|\widehat X_i\ind k - \widehat X_j\ind k\| 
\leq \gamma h_k$ implies $\|X_i - X_j\| \leq (\gamma + 4\beta_1) h_k$.
Next, we have
\begin{align*}
	&
	\p\ind{-i} \left( \|X_j - X_i\| \leq (\gamma + 4\beta_1) h_k \right)
	= \int\limits_{\M^*\cap \|x- X_i\| \leq (\gamma + 4\beta_1) h_k} p(x)dW(x)
	\\&
	\leq p_1 \int\limits_{\|v\| \leq (\gamma + 4\beta_1) h_k} \sqrt{\det g(v)} dv.
\end{align*}
Using the inequality $|\sqrt{\det g(v)} - 1| \lesssim d \|v\|^2/\varkappa^2$ (see 
\citep[Equation 2.1]{tsay19}), we have that $\sqrt{\det g(v)} \lesssim 1$, provided that 
$\|v\| \leq (\gamma + 4\beta_1) h_k$.
Then
\begin{align*}
	\p\ind{-i} \left( \|X_j - X_i\| \leq (\gamma + 4\beta_1) h_k \right)
	\lesssim \text{Vol}(\B(0,  (\gamma + 4\beta_1) h_k))
	\lesssim (\gamma + 4\beta_1)^d h_k^d.
\end{align*}
Thus, there exists a constant $C$, such that
\[
	\p\ind{-i} \left( \|X_j - X_i\| \leq (\gamma + 4\beta_1) h_k \right) \leq C (\gamma + 
	4\beta_1)^d h_k^d \, .
\]
Denote $C' = C \vee 16/3$.
The Bernstein's inequality yields
\begin{align*}
	&
	\p\ind{-i}\left( \sum\limits_{j=1}^n \1 \left( \|X_j - X_i\| \leq 2\gamma h_k \right) > 2C' 
	nh_k^d \right)
	\leq e^{-\frac{(C' nh_k^d)^2}{2 \cdot (C' nh_k^d) + 2/3 \cdot (C' nh_k^d)}}
	\\&\notag
	= e^{-\frac{3C'nh_k^d}8}
	\leq e^{-2n h_k^d}
	\leq e^{-2n h_K^d}
	\leq \frac 1{n^2},
\end{align*}
and then \eqref{e1} holds.
From now on, we are working on the event, on which \eqref{e1} holds.
On this event, we have
\begin{equation}
	\label{cov_op_norm}
	\|\widehat \bsigma_i\ind k - \widehat \bxi_i\ind k \| \leq 8\gamma (\gamma + 
	4\beta_1)^d \beta_2 C' nh_k^{d+3} \, .
\end{equation}

Consider the matrix $\widehat\bxi_i\ind k$.
According to Lemma \ref{lem_spectral_gap}, we have the following guarantee on the 
spectral gap of $\widehat\bxi_i\ind k$:
\begin{align*}
	&
	\lambda_d(\widehat\bxi_i) - \lambda_{d+1}(\widehat\bxi_i)
	\\&
	\geq \frac c4 \left(1 - \frac{2}{c (\gamma - 4\beta_1)^d} - 
	\sqrt{\frac{6C'}{c^2(\gamma - 4\beta_1)^d}} \right) (\gamma - 4\beta_1)^{d+2} 
	n h_k^{d+2}
	\\&
	- 9C' n^{-2/d} (\gamma + 4\beta_1)^{d+2} nh_k^{d+2}
	- 16 C'\beta_1^2 (\gamma + 4\beta_1)^d n h_k^{d+2}
	\\&
	- \frac{C' (\gamma + 4\beta_1)^{d+4} nh_k^{d+4}}{\varkappa^2}.
\end{align*}
with probability at least $1 - n^{-2}$.
Take $\gamma$ satisfying the inequalities
\begin{align}
	\label{gamma}
	&\notag
	(\gamma - 4\beta_1)^d \geq \frac 8c,
	\\&\notag
	(\gamma - 4\beta_1)^d \geq \frac{96C'}{c^2},
	\\&
	\frac{c}{32}(\gamma - 4\beta_1)^{d+2} \geq C' n^{-2/d} (\gamma + 4\beta_1)^{d+2},
	\\&\notag
	\frac{c}{32}(\gamma - 4\beta_1)^{d+2} \geq 16 C'\beta_1^2 (\gamma + 4\beta_1)^d,
	\\&\notag
	\frac{c}{32}(\gamma - 4\beta_1)^{d+2} \geq \frac{C' (\gamma + 4\beta_1)^{d+4}
	h_0^2}{\varkappa^2}.
\end{align}
Note that such $\gamma$ always exists if $n^{-2/d}$ and $h_0$ are 
sufficiently small.
Then 
\begin{equation}
	\label{spectral_gap}
	\lambda_d(\widehat\bxi_i) - \lambda_{d+1}(\widehat\bxi_i) \geq \frac c8 nh_k^{d+2} 
	- \frac{3c}{32} nh_k^{d+2} = \frac c{32} nh_k^{d+2} \, .
\end{equation}
The Davis-Kahan $\sin\theta$ theorem \citep{dk70} and the inequalities 
\eqref{cov_op_norm}, \eqref{spectral_gap} imply that for a fixed $i$ from $1$ to 
$n$ with probability at least $1 - 2n^{-2}$ it holds
\[
	\| \widehat\bpi_i\ind{k+1} - \bpi(X_i) \| \leq \frac{ 256\gamma (\gamma + 4\beta_1)^d 
	\beta_2 C' nh_k^{d+3} }{ c nh_k^{d+2} } = \widetilde C h_k
\]
with $\widetilde C = (256\gamma (\gamma + 4\beta_1)^d \beta_2 C')/c$.
Applying the union bound, we have that
\[
	\max \limits_{1\leq i \leq n}\| \widehat\bpi_i\ind{k+1} - \bpi(X_i) \| \leq \widetilde C h_k
\]
with probability at least $1 - 2n^{-1}$, and the proof of Lemma \ref{proj} is finished.

\section{Proof of Proposition \ref{lem_expectation}}
\label{app_e}

The proof of Proposition \ref{lem_expectation} is divided into three parts for the sake of 
convenience.
On each step we prove one of the inequalities \eqref{expbound}, \eqref{expbound2}, 
\eqref{expbound3}.

\subsection{Proof of Proposition \ref{lem_expectation}a}

First, consider the expression
\[
	\E\ind{-i} w_{ij}(\bpi_i)u^T(X_j - X_i).
\]
Let $r_d = 4h \sqrt{(d+2) \log h^{-1}}$.
Then
\begin{align*}
	\E\ind{-i} w_{ij}(\bpi_i)u^T(X_j - X_i)
	&
	= \E\ind{-i} w_{ij}(\bpi_i)u^T(X_j - X_i) \1\left( X_j \in \B(X_i, r_d) \right)
	\\&
	+ \E\ind{-i} w_{ij}(\bpi_i)u^T(X_j - X_i) \1\left( X_j \notin \B(X_i, 
	r_d) \right).
\end{align*}
Due to Lemma \ref{geom},
\begin{equation}
	\label{geom2}
	\|X_i - X_j\|^2 \leq 8\|\bpi_i(Y_i - Y_j)\|^2 + \frac{32M^2(\Delta+b/h)^2 
	h^2}{\varkappa^2},
\end{equation}
and, if $X_j \notin \B(X_i, r_d)$, we conclude
\begin{align*}
	&
	w_{ij}(\bpi_i)
	\leq e^{-\frac{\|X_i - X_j\|^2}{8h^2} + \frac{4M^2(\Delta+b/h)^2}{\varkappa^2}}
	\\&
	\leq e^{\frac{4M^2(\Delta+b/h)^2}{\varkappa^2} -\frac{\|X_i - X_j\|^2 + r_d^2}{16h^2}} 
	\leq e^{\frac14 -\frac{\|X_i - X_j\|^2}{16h^2}}h^{d+2}.
\end{align*}
Here we used the fact that, due to the conditions of Theorem \ref{th1}, $M(\Delta + 
b/h) \leq \varkappa/4$.
Using the equality
\[
	\max\limits_{t > 0} t e^{-\frac{t^2}{16h^2}} = 2h \sqrt2 e^{-1/2},
\]
we conclude
\begin{equation}
	\label{p1a5}
	\E\ind{-i} w_{ij}(\bpi_i)u^T(X_j - X_i) \1\left( X_j \notin \B(X_i, 
	r_d) \right) \lesssim h^{d+3}.
\end{equation}

We see that outside the ball $\B(X_i, r_d)$, the weights $w_{ij}(\bpi_i)$ become very 
small.
It remains to consider the event $\left\{ X_j \in \B(X_i, r_d) \right\}$.
We assume that $h_0$ is sufficiently small, so it holds $r_d \leq 2h_0 \sqrt{2(d+2) \log 
h_0^{-1}} \leq \tau/4$.
On this event $\|Y_i - Y_j \| \leq 2M + r_d < \tau$, which yields
\begin{align*}
	&
	\E\ind{-i} w_{ij}(\bpi_i)u^T(X_j - X_i) \1\left( X_j \in \B(X_i, r_d) \right) 
	\\&
	= \E\ind{-i} e^{-\frac{\|\bpi_i(Y_j - Y_i)\|^2}{h^2}} u^T(X_j - X_i) \1\left( X_j \in \B(X_i, 
	r_d) \right).
\end{align*}
Using the Taylor's expansion, one has
\begin{align}
	\label{d6}
	&
	e^{-\frac{\|\bpi_i(Y_j - Y_i)\|^2}{h^2}}
	= e^{-\frac{\|X_j - X_i\|^2}{h^2}}
	\\&\notag
	+ e^{-\frac{\|\xi\|^2}{h^2}}\left( \frac{ \|\bpi_i(Y_j - Y_i)\|^2 - \|X_j - X_i\|^2 }{h^2} \right),
\end{align}
where $\xi = \theta(X_j - X_i) + (1 - \theta)\bpi_i(Y_i - Y_j)$ for some $\theta \in (0,1)$.

Consider the expectation
\[
	\E\ind{-i} e^{-\frac{\|\bpi_i(X_j - X_i)\|^2}{h^2}} u^T(X_j - X_i) \1\left( X_j \in \B(X_i, 
	r_d) \right).
\]
According to Lemma \ref{p1a1}, it does not exceed
\[
	\E\ind{-i} e^{-\frac{\|\bpi_i(X_j - X_i)\|^2}{h^2}} u^T(X_j - X_i) \1\left( X_j \in \B(X_i, 
	r_d) \right) \lesssim \frac{d h^{d+2}}\varkappa.
\]

Next, consider the second term in \eqref{d6}.
Note that
\begin{align*}
	\|\xi\|
	&
	= \| \theta\bpi_i(Y_j - Y_i) + (1 - \theta)(X_j - X_i) \|
	\\&
	\geq \|X_i - X_j\| - \|\bpi_i(Y_j - Y_i) - (X_j - X_i) \| \, .
\end{align*}
From the proof of Lemma \ref{geom} (see Equation \ref{aa1}) we know that
\begin{align*}
	&
	\| X_j - X_i - \bpi_i(Y_j - Y_i) \|
	\leq \frac{\|X_j - X_i\|^2}{2\varkappa} + \frac{\Delta h \|X_j - X_i\|}\varkappa
	\\&\notag
	+ \frac{2M(\Delta h + b)}\varkappa + \frac{3M\|X_i - X_j\|}\varkappa.
\end{align*}
Then
\begin{align*}
	\|\xi\|
	&
	\geq \left(1 -  \frac{3M + \Delta h}\varkappa \right)\|X_i - X_j\| - \frac{\|X_i - 
	X_j\|^2}{2\varkappa} - \frac{2M(\Delta h + b)}\varkappa
	\\&
	\geq \left(1 -  \frac{3M + \Delta h + r_d}\varkappa \right) \|X_i - X_j\| 
	- \frac{2M(\Delta h + b)}\varkappa
	\\&
	\geq \left(1 -  \left(\frac{3}{16} + \frac14 + \frac18 \right) \right) \|X_i - X_j\| 
	- \frac{2M(\Delta h + b)}\varkappa
	\\&
	> \frac14 \|X_i - X_j\| - \frac{2M(\Delta h + b)}\varkappa.
\end{align*}
This yields
\begin{equation}
	\label{p1a6}
	\|X_i - X_j\|^2
	\leq 16\|\xi\|^2 + 128M^2(\Delta + b/h)^2h^2/\varkappa^2
	\leq 16\|\xi\|^2 + 8h^2.
\end{equation}
Then
\begin{align}
	\label{p1a2}
	&\notag
	\E\ind{-i} e^{-\frac{\|\xi\|^2}{h^2}} \left( \frac{ \|\bpi_i(Y_j - Y_i)\|^2 - \|X_j - X_i\|^2 
	}{h^2} \right)
	\\&
	\quad \cdot u^T(X_j - X_i) \1\left( X_j \in \B(X_i, r_d) \right)
	\\&\notag
	\leq e^{1/2} \E\ind{-i} e^{-\frac{\|X_i - X_j\|^2}{16h^2}} \left| \frac{ \|\bpi_i(Y_j - Y_i)\|^2 
	- \|X_j - X_i\|^2 }{h^2} \right|
	\\&\notag
	\quad \cdot \|X_j - X_i\| \1\left( X_j \in \B(X_i, r_d) \right)
\end{align}
Note that
\begin{align*}
	&
	\|\bpi_i(Y_j - Y_i)\|^2 - \|X_j - X_i\|^2
	= \|\bpi_i(\eps_j - \eps_i)\|^2
	\\&
	+ 2(\eps_j - \eps_i)^T\bpi_i(X_j - X_i) - \|(\bid - \bpi_i)(X_j - X_i)\|^2.
\end{align*}
Due to \citep[Theorem 4.18]{f59},
\begin{align*}
	\|(\bid - \bpi_i)(X_j - X_i)\|
	&
	\leq \|\bpi_i - \bpi(X_i)\| \|X_j - X_i\| + \|(\bid - \bpi(X_i))(X_j - X_i)\|
	\\&
	\leq \frac{\Delta h \|X_j - X_i\|}\varkappa + \frac{\|X_j - X_i\|^2}{2\varkappa}.
\end{align*}
Using \eqref{l1e1}, we obtain
\begin{align}
	\label{p1a3}
	&\notag
	\left| \|\bpi_i(Y_j - Y_i)\|^2 - \|X_j - X_i\|^2 \right|
	\\&\notag
	\leq \left(\frac{2M(\Delta h + b)}\varkappa + \frac{3M\|X_j - 
	X_i\|}\varkappa \right)^2
	\\&
	+ 2\|X_j - X_i\| \left(\frac{2M(\Delta h + b)}\varkappa + \frac{3M\|X_j - 
	X_i\|}\varkappa\right)
	\\&\notag
	+ \left(\frac{\Delta h \|X_j - X_i\|}\varkappa + \frac{\|X_j - X_i\|^2}{2\varkappa} 
	\right)^2.
\end{align}

Next, it is useful to control the expectations of the form
\[
	\E\ind{-i} e^{-\frac{\|X_i - X_j\|^2}{16h^2}} \|X_j - X_i\|^q \1\left( X_j \in \B(X_i, r_d) 
	\right).
\]
Lemma \ref{p1a4} yields
\[
	\E\ind{-i} e^{-\frac{\|X_i - X_j\|^2}{16h^2}} \|X_j - X_i\|^q \1\left( X_j \in \B(X_i, r_d) 
	\right) \lesssim h^{q+d}.
\]
The inequalities \eqref{p1a2} and \eqref{p1a3} and Lemma \ref{p1a4} yield that, up to a 
multiplicative constant, the left-hand side of \eqref{p1a2} is bounded by
\begin{align*}
	&
	\frac{M^2(\Delta + b/h)^2 h^{d+1}}{\varkappa^2} + \frac{M^2 
	h^{d+1}}{\varkappa^2} + \frac{M(\Delta + b/h) h^{d+1}}\varkappa + \frac{M 
	h^{d+3}}\varkappa + \frac{\Delta^2 h^{d+3}}{\varkappa^2} + 
	\frac{h^{d+3}}{\varkappa^2}
	\\&
	\lesssim \frac{M(\Delta + b/h + 1) h^{d+1}}\varkappa + \frac{(\Delta^2 + 1) 
	h^{d+3}}{\varkappa^2}.
\end{align*}
This and Lemma \ref{p1a1} imply
\begin{equation}
	\label{d4}
	\E\ind{-i} w_{ij}(\bpi_i)u^T(X_j - X_i) \lesssim \left(M(\Delta + b/h) \vee h \vee 
	\frac{\Delta^2 h^2}\varkappa \right) \frac{h^{d+1}}\varkappa
\end{equation}
with a hidden constant, which does not depend on $\Delta$.

\bigskip

\subsection{Proof of Proposition \ref{lem_expectation}b}

Consider the expectation
\[
	\E\ind{-i} \sum\limits_{j=1}^n w_{ij}(\bpi_i) u^T(\bid - \bpi(X_i))(X_j - X_i).
\]
Since for each $j \neq i$ the summand has the same conditional distribution with 
respect to $(X_i, Y_i)$, it is enough to prove that
\[
	\E\ind{-i} w_{ij}(\bpi_i) u^T(\bid - \bpi(X_i))(X_j - X_i) \lesssim \frac{h^{d+2}}\varkappa
\]
for any distinct $j$.

Again, we use the decomposition
\begin{align*}
	&
	\E\ind{-i} \sum\limits_{j=1}^n w_{ij}(\bpi_i) u^T(\bid - \bpi(X_i))(X_j - X_i)
	\\&
	= \E\ind{-i} \sum\limits_{j=1}^n w_{ij}(\bpi_i) u^T(\bid - \bpi(X_i))(X_j - X_i) \1(X_j \in 
	\B(X_i, r_d))
	\\&
	+  \E\ind{-i} \sum\limits_{j=1}^n w_{ij}(\bpi_i) u^T(\bid - \bpi(X_i))(X_j - X_i) \1(X_j 
	\notin \B(X_i, r_d)).
\end{align*}
From \eqref{p1a5}, the second term is of order $h^{d+3} \ll h^{d+2}/\varkappa$.	
On the event $\{ X_j \in \B(X_i, r_d) \}$, we can use \citep[Theorem 4.18]{f59}:
\begin{align*}
	&
	\E\ind{-i} w_{ij}(\bpi_i) u^T(\bid - \bpi(X_i))(X_j - X_i) \1(X_j \in 
	\B(X_i, r_d))
	\\&
	\leq \E\ind{-i} w_{ij}(\bpi_i) \left\| (\bid - \bpi(X_i))(X_j - X_i) \right\| \1(X_j \in 
	\B(X_i, r_d))
	\\&
	\leq \frac1{2\varkappa} \E\ind{-i} w_{ij}(\bpi_i) \|X_j - X_i\|^2 \1(X_j \in 
	\B(X_i, r_d)).
\end{align*}
Using \eqref{geom2}, we obtain
\[
	w_{ij}(\bpi_i)
	= e^{-\frac{\|\bpi_i(Y_j - Y_i)\|^2}{h^2}}
	\leq e^{-\frac{\|X_j - X_i\|^2}{8h^2} + \frac{4M^2(\Delta+b/h)^2}{h^2}}
	\leq e^{\frac12 -\frac{\|X_j - X_i\|^2}{8h^2}}.
\]
The assertion of Proposition \ref{lem_expectation}b now follows from Lemma \ref{p1a4}.

\bigskip

\subsection{Proof of Proposition \ref{lem_expectation}c}

To complete the proof of Proposition \ref{lem_expectation}, it remains to bound the 
expectation
\[
	\E\ind{-i} w_{ij}(\bpi_i)u^T\eps_j.
\]
Outside the ball $\B(X_i, r_d)$, we have
\[
	\E\ind{-i} w_{ij}(\bpi_i)u^T\eps_j \1(X_j \notin \B(X_i, r_d) \lesssim Mh^{d+2},
\]
so it remains to control the expectation
\[
	\E\ind{-i} w_{ij}(\bpi_i)u^T\eps_j \1(X_j \in \B(X_i, r_d).
\]

Again, use the Taylor's expansion:
\begin{align*}
	e^{-\frac{\|\bpi_i(Y_i - Y_j)\|^2}{h^2}}
	&
	= e^{-\frac{\|X_i - X_j\|^2}{h^2}}
	+ e^{-\frac{\|X_i - X_j\|^2}{h^2}} \left( \frac{ \|\bpi_i(Y_i - Y_j)\|^2 - \|X_i - X_j\|^2}{h^2} 
	\right)
	\\&
	+ e^{-\frac{\|\zeta\|^2}{h^2}} \frac{ \left( \|\bpi_i(Y_i - Y_j)\|^2 - \|X_i - X_j\|^2 
	\right)^2}{2h^4},
\end{align*}
where $\zeta = \vartheta\bpi_i(Y_j - Y_i) + (1 - \vartheta)(X_j - X_i)$ for some $\vartheta 
\in (0, 1)$.
On the event $\{ \|X_i - X_j \| \leq r_d \}$ it holds $\|Y_i - Y_j\| \leq 2M + r_d \leq \tau$.
This yields
\begin{align}
	\label{p1c5}
	&
	\E\ind{-i} e^{-\frac{\|X_i - X_j\|^2}{h^2}} \1(\|Y_i - Y_j\| \leq \tau) \1(\|X_i - X_j\| \leq 
	r_d) u^T\eps_j
	\\&\notag
	= \E\ind{-i} e^{-\frac{\|X_i - X_j\|^2}{h^2}} \1(\|X_i - X_j\| \leq 
	r_d) u^T\eps_j = 0.
\end{align}
Now, consider the term
\[
	\E\ind{-i} e^{-\frac{\|X_i - X_j\|^2}{h^2}} \left( \frac{ \|\bpi_i(Y_i - Y_j)\|^2 - \|X_i - 
	X_j\|^2}{h^2} \right)\1(\|X_i - X_j\| \leq r_d)  u^T\eps_j.
\]
It is equal to
\begin{align*}
	 &
	 \E\ind{-i} e^{-\frac{\|X_i - X_j\|^2}{h^2}} \1(\|X_i - X_j\| \leq r_d) u^T\eps_j
 	\\&
 	\cdot \left( \frac{ \|\bpi_i(\eps_i - \eps_j)\|^2 + 2(X_i - X_j)^T\bpi_i(\eps_i - \eps_j) - 
 	\|(\bid - \bpi_i)(X_i - X_j)\|^2}{h^2} \right).
\end{align*}
First, note that
\begin{equation}
	\label{p1c1}
	\E\ind{-i} e^{-\frac{\|X_i - X_j\|^2}{h^2}} \frac{\|(\bid - \bpi_i)(X_i - X_j)\|^2}{h^2} 
	\1(\|X_i - X_j\| \leq r_d) u^T\eps_j = 0.
\end{equation}
Next, \eqref{l1e1} implies
\[
	\|\bpi_i(\eps_i - \eps_j)\|^2
	\leq \left(\frac{2M(\Delta h + b)}\varkappa + \frac{3M\|X_j - X_i\|}\varkappa \right)^2.
\]
Then, using the inequality $\|\eps_j\| \leq M$ (due to \eqref{a2}) and Lemma \ref{p1a4}, 
we obtain
\begin{align}
	\label{p1c2}
	&
	\E\ind{-i} e^{-\frac{\|X_i - X_j\|^2}{h^2}} \1(\|X_i - X_j\| \leq r_d) u^T\eps_j  \frac{ 
	\|\bpi_i(\eps_i - \eps_j)\|^2}{h^2}
	\\&\notag
	\lesssim \frac{M^3(\Delta + b/h)^2 h^d}{\varkappa^2}.
\end{align}
Finally, consider the expectation
\[
	\E\ind{-i} e^{-\frac{\|X_i - X_j\|^2}{h^2}} \cdot \frac{ 2(X_i - X_j)^T\bpi_i(\eps_i - 
	\eps_j) }{h^2} \1(\|X_i - X_j\| \leq r_d) u^T\eps_j.
\]
Denote
\[
	v_{ij} = 2 \E\left( \bpi_i(\eps_i - \eps_j) u^T\eps_j \cond X_i, X_j ,Y_i \right).
\]
According to \eqref{l1e1}, the norm of the vector $v_{ij}$ is bounded by
\[
	\|v_{ij}\|
	\leq 2M \left\| \bpi_i(\eps_i - \eps_j) \right\|
	\leq  2M \left(\frac{2M(\Delta h + b)}\varkappa + \frac{3M\|X_j - X_i\|}\varkappa 
	\right).
\]
On the event $\{ \|X_i - X_j \| \leq r_d\}$, we have
\[
	\|v_{ij}\| \leq 2M \left(\frac{2M(\Delta h + b)}\varkappa + \frac{3M r_d}\varkappa 
	\right).
\]
Applying Lemma \ref{p1a1} with the vector $u = v_{ij}/ \|v_{ij}\|$, we obtain
\begin{align}
	\label{p1c3}
	&\notag
	\E\ind{-i} e^{-\frac{\|X_i - X_j\|^2}{h^2}} \cdot \frac{ 2(X_i - X_j)^T\bpi_i(\eps_i - 
	\eps_j) }{h^2} \1(\|X_i - X_j\| \leq r_d) u^T\eps_j
	\\&
	\lesssim \frac{Mh^{d}}\varkappa \left(\frac{M(\Delta h + b)}\varkappa + \frac{M 
	r_d}\varkappa \right)
	\\&\notag
	\lesssim
	\frac{Mh^{d}}\varkappa \left(\frac{M(\Delta h + b)}\varkappa + \frac{M 
	h\sqrt{\log h^{-1}}}\varkappa \right).
\end{align}
Taking \eqref{p1c1}, \eqref{p1c2} and \eqref{p1c3} together, one obtains
\begin{align}
	\label{p1c4}
	&
	\E\ind{-i} e^{-\frac{\|X_i - X_j\|^2}{h^2}} \left( \frac{ \|\bpi_i(Y_i - Y_j)\|^2 - \|X_i - 
	X_j\|^2}{h^2} u^T\eps_j  \right)
	\\&\notag
	\lesssim \frac{M^3(\Delta + b/h)^2 h^d}{\varkappa^2} + \frac{Mh^{d}}\varkappa 
	\left(\frac{M(\Delta h + b)}\varkappa + \frac{M h\sqrt{\log h^{-1}}}\varkappa 
	\right).
\end{align}

To complete the proof of Proposition \ref{lem_expectation}, it remains to bound
\[
	\E\ind{-i} e^{-\frac{\|\zeta\|^2}{h^2}} \cdot \frac{ \left( \|\bpi_i(Y_i - Y_j)\|^2 - \|X_i - 
	X_j\|^2 \right)^2}{2h^4} \1(\|X_i - X_j\| \leq r_d) u^T\eps_j,
\]
where $\zeta = \vartheta\bpi_i(Y_j - Y_i) + (1 - \vartheta)(X_j - X_i)$ for some $\vartheta 
\in (0, 1)$.
The same argument, as in the analysis of the vector $\xi$ (see the proof of Proposition 
1a, Inequality \ref{p1a6}), yields
\[
	\|X_i - X_j\|^2
	\leq 16\|\zeta\|^2 + 128M^2(\Delta + b/h)^2h^2/\varkappa^2
	\leq 16\|\zeta\|^2 + 8h^2
\]
and then
\[
	e^{-\frac{\|\zeta\|^2}{h^2}} \leq e^{1/2} e^{-\frac{\|X_i - X_j\|^2}{16h^2}}.
\]
Due to \eqref{p1a3},
\begin{align*}
	&
	\left( \|\bpi_i(Y_i - Y_j)\|^2 - \|X_i - X_j\|^2 \right)^2
	\leq \left[  \left(\frac{2M(\Delta h+ b)}\varkappa + \frac{3M\|X_j - 
	X_i\|}\varkappa \right)^2
	\right.
	\\&
	+ 2\|X_j - X_i\| \left(\frac{2M(\Delta h + b)}\varkappa + \frac{3M\|X_j - 
	X_i\|}\varkappa\right)
	\\&\notag
	\left.
	+ \left(\frac{\Delta h \|X_j - X_i\|}\varkappa + \frac{\|X_j - X_i\|^2}{2\varkappa} 
	\right)^2 \right]^2
\end{align*}
Applying the inequality $(a + b + c)^2 \leq 3a^2 + 3b^2 + c^2$ and Lemma \ref{p1a4}, we 
obtain
\begin{align}
	\label{p1c6}
	&
	\E\ind{-i} e^{-\frac{\|\zeta\|^2}{h^2}} \cdot \frac{ \left( \|\bpi_i(Y_i - Y_j)\|^2 - \|X_i - 
	X_j\|^2 \right)^2}{2h^4} \1(\|X_i - X_j\| \leq r_d) u^T\eps_j
	\\&\notag
	\lesssim M h^{d-4} \left( \frac{M^4(\Delta + b/h)^4 h^4}{\varkappa^4} + 
	\frac{M^4 h^4}{\varkappa^4} + \frac{M^2 (\Delta + b/h)^2 h^4}{\varkappa^2} + 
	\frac{M^2 h^4}{\varkappa^2} + \frac{\Delta^4 h^8}{\varkappa^4} + 
	\frac{h^8}{\varkappa^4} \right)
	\\&\notag
	\lesssim \frac{M^5(1 + \Delta + b/h)^4 h^d}{\varkappa^4} + 
	\frac{M^3 (1 + \Delta + b/h)^2 h^d}{\varkappa^2} + \frac{(1 + \Delta^4) 
	M h^{d+4}}{\varkappa^4}.
\end{align}
The assertion of Proposition \ref{lem_expectation}c follows from Inequalities 
\ref{p1c5}, \ref{p1c4}, \ref{p1c6} and the fact that, due to conditions of Theorem 
\ref{th1},
\[
	\frac{M(1 + \Delta + b/h)}{\varkappa} \leq \frac1{16} + \frac1{4} = \frac{5}{16}.
\]

\section{Proofs Related to Theorem \ref{lower_bound}}
\label{app_lower_bounds}

\subsection{Proof of Lemma \ref{lem_lower_bound_1}}
\label{app_lower_bounds_1}

It is enough to consider the case $D = d + 1$.
For any $x_0 \in \R^{d+1}$ and $r > 0$, denote a sphere of radius $r$ centered at $x_0$ 
by $\partial \B(x, r) = \{x \in \R^{d+1} : \|x - x_0\| = r\}$.
Consider $\M_0 = \partial B(0, \varkappa)$.
Let a random element $X$ have a uniform distribution on $\M_0$.
Clearly, $\M_0$ satisfies \eqref{a1} and the distribution of $X$ satisfies \eqref{a1'} with 
$L = 0$, $p_0 = p_1 = ((d+1) \omega_{d+1}  \varkappa^d)^{-1}$, where $\omega_{d + 1}$ 
is the volume of the Euclidean ball in $\R^{d+1}$ with radius $1$.
Given $X$, let $\eps$ have a uniform distribution on $\T_X \M_0 \cap \partial \B(X, 
Mb/\varkappa)$.
Then
\[
\E(\eps \cond X) = 0,
\quad
\|\eps\| = \frac{Mb}{\varkappa}
\quad
\p( \cdot \cond X)\text{-almost surely},
\]
and the assumption \eqref{a2} is fulfilled.
Consider $Y = X + \eps$.
Since $X$ is orthogonal to $\eps$ by the construction, we have
\[
\|Y\|^2 = \|X\|^2 + \|\eps\|^2 = \varkappa^2 + \frac{M^2 b^2}{\varkappa^2}
\quad \text{almost surely}.
\]
Consequently, the random element $Y$ is supported on the sphere $\partial \B(0, R)$ 
with $R = \sqrt{\varkappa^2 + M^2 b^2 / \varkappa^2}$.
By the spherical symmetry construction, $Y$ has a uniform distribution on $\partial \B(0, 
R)$.

Let $X'$ have a uniform distribution on $\M_1 = \partial \B(0, R)$ and, for any $X'$, let 
$(\eps' \cond X') = 0$ $\p(\cdot \cond X')$-almost surely.
Then $X' + \eps'$ and $X + \eps$ have the same distribution.
However,
\[
	d_H(\M_1, \M_0)
	= R - \varkappa
	= \varkappa \left( \sqrt{1 + \frac{M^2 b^2}{\varkappa^4}} - 1 \right)
	> \frac{M^2 b^2}{3\varkappa^3}.
\]
Here we used the fact that, due to the concavity of $\sqrt{1 + x} - 1$, it holds that
\[
	\sqrt{1 + x} - 1 \geq (\sqrt 2 - 1) x > \frac x3, \quad \forall x \in [0,1].
\]
Thus, for any estimate $\widehat\M$, we have
\[
	\sup\limits_{\M^* \in \mclass_\varkappa^d} d_H(\widehat\M, \M)
	\geq \frac12 d_H(\M_1, \M_0)
	> \frac{M^2 b^2}{6\varkappa^3}.
\]

\subsection{Proof of Lemma \ref{lem_lower_bound_2}}
\label{app_lower_bounds_2}

Without loss of generality, we assume \( D = d+1 \).
We write a \( (d+1) \)-dimensional vector as \( (u, v) \), where \( u \in \R^{d} \), \( v \in \R \).
Let $\Z^{(0)} \subset \R^D$ be a $d$-dimensional $\C^\infty$-manifold without a 
boundary with reach greater than $1$ such that
\[
	\{ (z, 0) : z \in \R^d, \|z\| \leq 1/2 \} \subset \Z^{(0)}.
\]
In \citep{al19}, the authors claim that such a manifold can be constructed by flattering 
smoothly a unit sphere in $\R^{D}$.
Let $\Z = 4\varkappa \Z^{(0)}$.
Then $\Z$ is a $\C^\infty$-manifold without a boundary.
Moreover, its reach is at least $4\varkappa$ and
\[
	\{ (z, 0) : z \in \R^d, \|z\| \leq 2\varkappa \} \subset \Z.
\]
We construct manifolds in the following way.
Let \( \psi : \R^d \rightarrow \R \) be a smooth function, such that \( \max_u \psi(u) = 
\psi(0) 
= 1 \), \( \psi(u) = 0 \) for any \( u \notin \B(0,1) \) and \( \sup_u \| \nabla^2 
\psi(u) \| \leq \Lambda \) for an absolute constant \( \Lambda \).
Let \( (z_1, 0), \dots, (z_N, 0) \), where \( z_1, \dots, z_N \in \R^d \), be a \( 2h \)-packing of 
\( d \)-dimensional ball \( \Z \cap B(0, \varkappa/2) \), \( N = 
\left(4h/\varkappa\right)^{-d} \), \( h < \varkappa/4 \).
For any \( j \in \{1, \dots, N\} \), introduce a manifold
\[
\M_j = \left\{ \begin{pmatrix} z\\0 \end{pmatrix} + \frac{h^2}{\varkappa L} \psi\left( \frac{z - z_j}h \right)e_{d+1} : 
(z, 0) \in \Z \cap \B(0, \varkappa) \right\}
\cup \left( \Z \backslash \B(0, \varkappa) \right),
\]
where the vector \( e_i \) is the \( i \)-th vector of the canonical basis in \( \R^{d+1} \) with 
the 
components \( \smash{e_i\ind j = \1(i=j)} \).
Let \( \M_0 \) be equal to \( \Z \).
Notice that \( \M_j \), \( j \in \{1, \dots, N\} \) differs from \( \M_0 \) only on the set 
\( \B((z_j, 0), h) \), and for any \( k \neq j \) the balls \( \B((z_j, 0), h) \) and \( \B((z_k, 0), h) 
\) do 
not intersect.
In other words, we consider a family of manifolds with a small bump in one of the points 
\( (z_1, 0), \dots, (z_N, 0) \).

Show that the family of manifolds \( \mclass^\circ_\varkappa = \left\{ \M_j : 1 \leq j \leq N 
\right\} \) with
\[
	h = c_0 \left( \frac{M^2\varkappa^2 \log n}n \right)^{1/(d+4)},
\]
where \( c_0 \) is a constant to be chosen later, is contained in the class \( \mclass_\varkappa^d \) introduced in 
\eqref{a1}.
It is clear that \( \M_j \) is a compact, connected, smooth \( d \)-dimensional manifold 
without a boundary for any \( j \) from \( 1 \) to \( N \).
The most important part is to check that the reach of \( \M_j \) is not less than 
\( \varkappa \).
For this purpose, we use Theorem 4.18 from \cite{f59}, which states that \( \reach\M \geq 
\varkappa \) if and only if for any \( x, x' \in \M \) it holds \( d(x', \{x\} \oplus \T_x\M) \leq \|x 
- x'\|^2/(2\varkappa) \).
Fix arbitrary \( j \in \{1, \dots, N\} \) and introduce
\[
	f_j(z) = \begin{pmatrix} z\\0 \end{pmatrix} + \frac{h^2}{L\varkappa} \psi\left( \frac{z - z_j}h \right) e_{d+1}, \quad z \in \Z.
\]
Then for any \( x \in \M_j \cap \B(0, \varkappa) \) the exists unique \( (z, 0)\in\Z \), such 
that \( x = f_j(z) \).
By the construction, the inverse function to \( f_j(z) \), \( (z, 0) \in \Z \cap \B(0, \varkappa) 
\), 
is given by
\[
	f_j^{-1}(x) = \left(x\ind1, \dots, x\ind d\right)^T,
\]
where \( x \in \M_j \cap \B(0, \varkappa) \) and \( x\ind j \) is the \( j \)-th component of the vector \( x \).
Moreover, the unit normal to \( \M_j \) at the point \( x = f_j(z) \) is given by
\begin{equation}
\label{normal}
\nu_j(z) = C_h^{-1} \left(-\frac{h}{\varkappa L}\nabla\psi\left(\frac{z - 
	z_j}h\right)^T, 1 \right)^T,
\end{equation}
where
\[
C_h = \sqrt{1 + \left(\frac{h}{\varkappa L}\right)^2 \left\|\nabla\psi\left(\frac{z - z_j}h 
	\right) \right\|^2}.
\]
Fix arbitrary \( x = f_j(z), x_0 = f_j(z_0) \in \M_j \) and check that
\[
\left|\nu_j(z_0)^T(x - x_0) \right|
= \left|\nu_j(z_0)^T(f_j(z) - f_j(z_0)) \right|
\leq \frac{\|z - z_0\|^2}{2\varkappa}
\leq \frac{\|x - x_0\|^2}{2\varkappa} \, .
\]
The last inequality is obvious, since \( (z - z_0) \) is a subvector of \( (x - x_0) \).
It remains to check the second inequality.
It holds that
\begin{align*}
	&
	\left|\nu_j(z_0)^T(f_j(z) - f_j(z_0)) \right|
	\\&
	= C_h^{-1} \left| -\frac{h}{\varkappa \Lambda}\nabla\psi^T\left(\frac{z_0 - 
	z_j}h\right)(z - 
	z_0) + \frac{h^2}{\varkappa \Lambda} \left(\psi\left(\frac{z - z_j}h\right) - 
	\psi\left(\frac{z_0 - z_j}h\right) \right) \right|
	\\&
	\leq  \left| -\frac{h}{\varkappa \Lambda}\nabla\psi^T\left(\frac{z_0 - z_j}h\right)(z - 
	z_0) + \frac{h^2}{\varkappa \Lambda} \left(\psi\left(\frac{z - z_j}h\right) - 
	\psi\left(\frac{z_0 - z_j}h\right)\right) \right|
	\\&
	\leq \frac{h^2}{\varkappa \Lambda} \cdot \frac{\Lambda \|z - z_0\|^2}{2h^2} 
	= \frac{\|z - z_0\|^2}{2\varkappa}.
\end{align*}
Here we used Taylor's expansion of \( \psi \) up to the second order and the fact that 
\( \|\nabla^2\psi\| \leq \Lambda \).

Now, we are going to describe distributions of \( X \) and \( \eps \) in the model 
\eqref{model}.
Let a random element \( Z \) have a uniform distribution on \( \Z \).
For any fixed \( j \in \{1, \dots, N\} \), we take \( X = f_j(Z) \), where
\[
	f_j(z) = \begin{pmatrix} z\\0 \end{pmatrix} + \frac{h^2}{\varkappa \Lambda} \psi\left(\frac{z - z_j}h\right) e_{d+1}, \quad  z \in \Z.
\]
Denote a volume of the set \( \Z \) by \( V_\Z \).
Then for any \( x \), such that \( f_j^{-1}(x) \notin \B(z_j, h) \), the density \( p_j(x) \) is just 
\( V_\Z^{-1} \).
Otherwise, the density \( p_j(x) \) of \( X \) is defined by the formula
\begin{align*}
	p_j(x)
	&
	= \frac1{V_\Z} \left( \det \nabla f_j(f_j^{-1}(x)) \right)^{-1}
	\\&\notag
	= \frac1{V_\Z} \left( \det\left( I + \frac{h}{\varkappa \Lambda} \nabla\psi\left( 
	\frac{f_j^{-1}(x) - z_j}h \right)e_{d+1}^T \right)\right)^{-1}.
\end{align*}
Since for any two vectors \( u, v \in \R^{d+1} \) it holds \( \det(I + u v^T) = 1 + u^T v \), we 
have
\begin{equation}
\label{qc}
	p_j(x)
	= \frac1{V_\Z} \left( 1 + \frac{h}{\varkappa \Lambda} e_{d+1}^T\nabla\psi\left( 
	\frac{f_j^{-1}(x) - z_j}h \right) \right)^{-1}.
\end{equation}
Note that \( \nabla\psi(0) = 0 \) by construction. Taking into account that 
\( \sup_u\|\nabla^2\psi(u)\| \leq \Lambda \), we conclude that
\[
	\left\|\nabla\psi\left( \frac{f_j^{-1}(x) - z_j}h \right) \right\|
	\leq \frac{\Lambda \|f_j^{-1}(x) - z_j\|}h \leq \Lambda \quad \forall \, x : f_j^{-1}(x) \in 
	\B(z_j, h).
\]
This and the fact that \( h < \varkappa/4 \) yield
\[
	p_0 = \frac4{5V_\Z}
	\leq \frac1{V_\Z\left(1 + \frac{h}\varkappa\right)}
	\leq p_j(x)
	\leq \frac1{V_\Z\left(1 - \frac{h}\varkappa\right)}
	\leq \frac4{3V_\Z} = p_1.
\]
Thus, the density of \( X \) is bounded from above and below by \( p_1 \) and \( p_0 \) 
respectively.

Show that \( p_j(x) \) has \( 4p_1/(3\varkappa) \)-Lipschitz derivative.
Differentiating \eqref{qc}, we obtain
\[
	\|\nabla p_j(x)\|
	= \frac{V_\Z  p_j^2(x)}{\varkappa \Lambda} \left\| \bid_{d, d+1} \nabla^2 \psi\left( 
	\frac{f_j^{-1}(x) - z_j}h \right) e_{d+1} \right\|
	\leq \frac{V_\Z  p_1^2}{\varkappa}
	\leq \frac{4p_1}{3\varkappa},
\]
where \( \bid_{d, d+1} \in \R^{d\times(d+1)} \) is the matrix of the first \( d \) rows of the 
identity matrix \( \bid_{d+1} \).
Thus, for each \( j \), the density \( q_j(x) \) fulfils \eqref{a1'}.

Next, we describe the conditional distribution of \( Y \) given \( X \).
We generate \( Y \) from the model
\begin{equation}
\label{y}
Y = X + \xi e_{d+1} \1(X \in \B(0, \varkappa)), \quad X \in \M_j,
\end{equation}
where \( \p(\xi = 0.5M - X\ind{d+1} | X) = \eta(X), \p(\xi = -0.5M - X\ind{d+1}) = 1 - 
\eta(X) \), \( \eta = \eta(X) = 1/2 + X\ind{d+1} / M \), and \( X\ind{d+1} \) is the \( (d+1) 
\)-th component of \( X \).
Note that \( Y \) belongs either to the set \( \M^+ = \{ (z, 0.5M) : (z, 0) \in \B(0, \varkappa) 
\subset \R^{d+1}\} \), or to the set \( \M^- = \{ (z, 0.5M) : (z, 0) \in \B(0, \varkappa) \subset 
\R^{d+1} \} \), or to the set \( \Z \backslash B(0, \B(0, \varkappa)) \).
It remains to check the condition \eqref{a2}.
Take
\[
	h = c_0 \left(\frac{M^2\varkappa^2 \log n}n \right)^{1/(d+4)},
\]
where \( c_0 \) is such that the condition \( h^2/(M\varkappa \Lambda) \leq 1/2 \) is fulfilled.
Such \( c_0 \) exists since \( M \gtrsim (\log n / n)^{2/d} \). 
First, note that the noise magnitude is not greater than \( 0.5M + h^2/(\varkappa 
\Lambda) \), which is less than \( M \). 
Second, using the expression \eqref{normal} of the unit normal to \( \M_j \) at the 
point \( x = f_j(z) \), we have
\begin{align*}
	&
	\| \xi \bpi(f_j(z)) e_{d+1}\|^2
	= |\xi|^2 - |\xi e_{d+1}^T \nu_j(z)|^2
	\\&
	= |\xi|^2 - \frac{|\xi|^2}{1 + h^2 \|\nabla\psi(z/h)\|^2/ (\varkappa \Lambda)^2}
	\\&
	= \frac{|\xi|^2 h^2 \|\nabla\psi(z/h)\|^2}{\varkappa^2 \Lambda^2 + h^2 
	\|\nabla\psi\left(\frac{z - z_j}h \right)\|^2}
	\leq \Lambda^2 |\xi|^2 \cdot \frac{h^2}{\varkappa^2 \Lambda^2}
	\leq \frac{M^2 h^2}{\varkappa^2},
\end{align*}
and \eqref{a2} holds with
\[
	b = c_0 \left(\frac{M^2\varkappa^2 \log n}n \right)^{1/(d+4)}.
\]
Here we used the fact that for any \( u \in \B(0, 1) \)
\[
\|\nabla\psi(u)\| = \|\nabla\psi(u) - \nabla\psi(0)\| \leq \max\limits_{u'\in\B(0,1)} 
\|\nabla^2\psi(u')\| \|u\| \leq L \, .
\]

We use \citep[Theorem 2.5]{t09} to prove the lower bound in Theorem \ref{lower_bound}.
Let \( P_j \), \( 0 \leq j \leq N \), be the probability measure, generated by \( Y = X + \eps \), 
where \( X \in \M_j \).
Then \( P_0 \) is a dominating measure, i.e. \( P_j \ll P_0 \) for all \( j \) from \( 1 \) to \( N \), 
and for 
any \( j \neq k \) we have
\[
d_H(\M_j, \M_k) \geq \frac{h^2}{\varkappa L}.
\]
Prove that, for sufficiently large \( n \), it holds
\begin{equation}
\label{cond_th25}
\frac1N \sum\limits_{j=1}^N \kl(P_j^{\otimes n}, P_0^{\otimes n}) \leq \alpha \log N,
\end{equation}
where \( \kl(P, Q) \) is the Kullback-Leibler divergence between \( P \) and \( Q \), \( \alpha 
\) is 
a constant from the interval \( (0, 1/8) \).
Then Theorem 2.5 in \cite{t09} yields
\[
\inf\limits_{\widehat \M} \sup\limits_{\M^* \in \mclass_\varkappa^d} \E 
d_H(\widehat\M, \M^*) \gtrsim \varkappa^{-1} \left(\frac{M^2\varkappa^2 \log n}n 
\right)^{1/(d+4)}.
\]

It remains to check \eqref{cond_th25}.
For any \( j \) from \( 1 \) to \( N \), it holds
\[
\kl(P_j^{\otimes n}, P_0^{\otimes n})  = n \int\limits \log \frac{dP_j(y)}{dP_0(y)} 
dP_j(y).
\]
The density of \( P_j \) with respect to the Hausdorff measure on \( (\Z\backslash B(0, 
\varkappa)) \cup \M^+ \cup \M^{-} \) is given by the formula
\[
	\log p_j(y) = 
	\begin{cases}
		\log\left(1 + \frac{2h^2}{M\varkappa \Lambda} \psi\left(\frac{\widetilde y - z_j}h 
		\right) \right) - \log(2V_\Z), \quad y = (\widetilde y, 0.5 M) \in \M^+,
		\\
		\log\left(1 - \frac{2h^2}{M\varkappa \Lambda} \psi\left(\frac{\widetilde y - z_j}h 
		\right) \right) - \log(2V_\Z), \quad y = (\widetilde y, -0.5 M) \in \M^-,
		\\
		- \log(V_\Z), \quad y \in \Z\backslash \B(0, \varkappa).
	\end{cases}
\]
The density of \( P_0 \) with respect to the same measure is just
\[
	\log p_0(y) = 
	\begin{cases}
		- \log(2V_\Z), \quad y = (\widetilde y, \pm M) \in \M^+ \cup \M^-,
		\\
		- \log(V_\Z), \quad y \in \Z\backslash \B(0, \varkappa).
	\end{cases}
\]
Then
\begin{align*}
	&
	\kl(P_j^{\otimes n}, P_0^{\otimes n}) 
	\\&
	= \frac{n}{2 V_\Z} \int\limits_{\|z - z_j\| \leq h} 
	\log\left(1 + \frac{2h^2}{M\varkappa \Lambda} \psi\left(\frac{z - z_j}h \right) \right) 
	\left(1 + \frac{2h^2}{M\varkappa \Lambda} \psi\left(\frac{z - z_j}h \right) \right) dz
	\\&
	+ \frac{n}{2 V_\Z} \int\limits_{\|z - z_j\| \leq h} 
	\log\left(1 - \frac{2h^2}{M\varkappa \Lambda} \psi\left(\frac{z - z_j}h \right) \right) 
	\left(1 - \frac{2h^2}{M\varkappa \Lambda} \psi\left(\frac{z - z_j}h \right) \right) dz.
\end{align*}
Using the inequality
\[
(1+t)\log(1 + t) \leq t + \frac{t^2}2, \quad t \in (-1, 1),
\]
we obtain that
\[
(1+t)\log(1 + t) + (1-t)\log(1 - t) \leq t^2, \quad t \in (-1, 1).
\]
Then, since \( 2h^2 < M\varkappa \Lambda \), we have
\[
	\kl(P_j^{\otimes n}, P_0^{\otimes n}) 
	\leq \frac{n}{2 V_\Z} \int\limits_{\|z\| \leq h} \frac{4h^4}{M^2\varkappa^2 \Lambda^2}
	\psi^2\left(\frac{z}h \right) dz
	= \frac{C_{\Lambda, \Z} n h^{d+4}}{M^2 \varkappa^2},
\]
where \( C_{\Lambda, \Z} > 0 \) is a constant, depending on \( \Z \) and \( \Lambda \).
Substituting \( h \) by  \( c_0 \left(\frac{M^2\varkappa^2 \log n}n \right)^{1/(d+4)} \), we 
obtain
\[
	\frac{C_{\Lambda, \Z} n h^{d+4}}{M^2 \varkappa^2} = c_0^{d+4} C_{\Lambda,\Z} 
	\log n.
\]
On the other hand,
\begin{align*}
	\log N
	&
	= \log \left(\frac{4h}{\varkappa}\right)^{-d} 
	= d\log \frac{\varkappa}4 - d \log c_0 + \frac{d}{d+4} \log 
	\frac{n}{M^2\varkappa^2 \log n}
	\\&
	\geq - d \log c_0 + \frac{d}{2(d+4)} \log n,
\end{align*}
where in the last inequality we assumed that \( n \) is large, so it holds \( M^2\varkappa^2 
\log n \leq n^{1/4} \) and \( \varkappa \geq 4 n^{-1/(4d + 16)} \).
Choose any constant \( c_0 > 0 \), satisfying the inequality
\[
8 c_0^{d+4} C_{L,\Z} \log n + d\log c_0 < \frac{d}{2(d+4)} \log n.
\]
Such constant always exists.
Thus, \eqref{cond_th25} is fulfilled, and \citep[Theorem 2.5]{t09} yields the claim of 
Theorem \ref{lower_bound}.

\section{Auxiliary Results}

This section contains some auxiliary results, which are used in the proofs.
The results below often use technique concerning integration over manifolds.
Therefore, we would like to start with a short background, which will help a reader follow 
the proofs.

Given $x \in \M$ and $v \in \T_x\M$, let $\gamma(t, x, v)$ be a geodesic starting at $x$, 
such that
\[
	\left.\frac{d\gamma(t, x, v)}{dt}\right|_{t=0} = v.
\]
The exponential map of $\M$ at the point $x$ $\Exp_x : \T_x\M \rightarrow \M$ is defined 
as $\Exp_x(p) = \gamma(1, x, v)$.
Note that $d_\M(\Exp_x(v), x) = \|v\|$ for $v \leq \varkappa/4$, where $d_\M(x, x')$ is the 
length of the shortest path on $\M$ between $x$ and $x'$.
We extensively use integration over manifolds.
In these cases, the exponential map is useful to apply the change of variables formula.
For a small open set $U$, $x \in U \subset \M$, it holds
\[
	\int\limits_{U} f(x) dW(x) = \int\limits_{\Exp_x^{-1}(U)} f(\Exp_x(p)) \sqrt{\det g(p)} dp,
\]
where $g(p)$ is the metric tensor.
Without going deep into details, we just mention that the metric tensor allows a nice 
decomposition (see, for instance, \citep[Equation 2.1]{tsay19}), which is enough for our 
purposes:
\[
	\left| \sqrt{\det g(v)} - 1 \right| \lesssim \frac{d \|v\|^2}{\varkappa^2}.
\]

\begin{Prop}
	\label{aux}
	Let $\M \in \mclass_\varkappa^d$ and $x, x' \in \M$, $\|x - x'\| \leq 2\varkappa$.
	Let $\bpi(x)$ and $\bpi(x')$ be the projectors onto the tangent spaces $\T_x\M$ 
	and $\T_{x'}\M$ respectively.
	Then the following inequalities hold:
	\begin{align*}
		&
		\tag{a}
		\|x - x'\| \leq d_{\M}(x, x') \leq 2\|x - x'\|,
		\\&
		\tag{b}
		\|\bpi(x) - \bpi(x')\| \leq \frac{\|x - x'\|}\varkappa.
	\end{align*}
\end{Prop}

\begin{proof}
	Proposition \ref{aux}a follows from \citep[Lemma 2.5]{blw18} and the inequality
	\[
		\arcsin \frac{2t}\pi \leq t, \quad \forall \, t \in (0, \pi/2).
	\]
	
	To prove Proposition \ref{aux}b, we use \citep[Theorem 2.5.1]{gvl13}:
	\[
		\|\bpi(x) - \bpi(x')\|
		= \|(\bid - \bpi(x)) \bpi(x')\|
		= \sin \angle(\T_x\M, \T_{x'}\M).
	\]
	Then the claim of the proposition follows from \citep[Corollary 3.6]{blw18}:
	\[
		\|\bpi(x) - \bpi(x')\|
		= \sin \angle(\T_x\M, \T_{x'}\M)
		\leq 2 \sin \frac{\angle(\T_x\M, \T_{x'}\M)}2
		\leq \frac{\|x - x'\|}\varkappa.
	\]
	
\end{proof}

\begin{Lem}
	\label{lem_spectral_gap}
	Fix any $i$ from $1$ to $n$.
	There are absolute constants $c$ and $C'$, such that, with probability at least $1 - 
	n^{-2}$, it holds
	\begin{align*}
	&
	\lambda_d(\widehat\bxi_i) - \lambda_{d+1}(\widehat\bxi_i)
	\\&
	\geq \frac c4 \left(1 - \frac{2}{c (\gamma - 4\beta_1)^d} - 
	\sqrt{\frac{6C'}{c^2(\gamma - 4\beta_1)^d}} \right) (\gamma - 4\beta_1)^{d+2} 
	n h_k^{d+2}
	\\&
	- 9C' n^{-2/d} (\gamma + 4\beta_1)^{d+2} nh_k^{d+2}
	- 16 C'\beta_1^2 (\gamma + 4\beta_1)^d n h_k^{d+2}
	\\&
	- \frac{C' (\gamma + 4\beta_1)^{d+4} nh_k^{d+4}}{\varkappa^2}.
	\end{align*}
\end{Lem}

\begin{proof}

Now, consider the matrix $\smash{\widehat\bxi\mathstrut_i\ind k}$.
It is clear that all the eigenvectors of $\smash{\widehat\bxi\mathstrut_i\ind k}$ belong to the linear space 
$\T_{X_i}\M^*$.
Thus, $\smash{\widehat\bxi\mathstrut_i\ind k}$ has at most $d$ non-zero eigenvalues.
In what follows, we show that the $d$-th largest eigenvalue of $\smash{\widehat\bxi\mathstrut_i\ind k}$ is 
non-zero and give a lower bound on the spectral gap $\smash{\lambda_d(\widehat\bxi\mathstrut_i\ind k) - 
\lambda_{d+1}(\widehat\bxi\mathstrut_i\ind k)}$.
It holds that
\[
\lambda_d(\widehat\bxi_i\ind k) - \lambda_{d+1}(\widehat\bxi_i\ind k)
= \min\limits_{u\in \T_{X_i}\M^*, \|u\|=1} \sum\limits_{j=1}^n v_{ij} (u^T(\widehat 
Z_{ij}\ind k - \widehat Z_{ii}\ind k))^2 \, .
\]
Using the inequality
\begin{align*}
&
\| \widehat Z_{ij}\ind k - X_j \|
\leq \| \widehat Z_{ij}\ind k - Z_{ij} \| + \| X_j - Z_{ij} \|
\\&
\leq \| \widehat X_j\ind k - X_j \| + \| X_j - Z_{ij} \|
\leq 2\beta_1 h_k + \frac{\|X_j - X_i\|^2}{2\varkappa},
\end{align*}
we obtain that for any $u$ it holds
\begin{align*}
(u^T(\widehat Z_{ij}\ind k - \widehat Z_{ii}\ind k))^2
\geq \frac12 (u^T(X_j - X_i))^2 - 8\beta_1^2 h_k^2 - \frac{\|X_j - 
	X_i\|^4}{2\varkappa^2},
\end{align*}
which yields
\begin{align}
\label{e2}
\lambda_d(\widehat\bxi_i\ind k) - \lambda_{d+1}(\widehat\bxi_i\ind k)
&\notag
\geq \frac12 \min\limits_{u\in \B(0,1) \cap \T_{X_i}\M^*} \sum\limits_{j=1}^n v_{ij} 
(u^T(X_j - X_i))^2
\\&
- 8\beta_1^2 h_k^2 \sum\limits_{j=1}^n v_{ij} - \sum\limits_{j=1}^n v_{ij} \frac{\|X_j - 
	X_i\|^4}{2\varkappa^2}
\\&\notag
\geq \frac12 \min\limits_{u\in \B(0,1) \cap \T_{X_i}\M^*} \sum\limits_{j=1}^n v_{ij} 
(u^T(X_j - X_i))^2
\\&\notag
- 16 C'\beta_1^2 (\gamma + 4\beta_1)^d n h_k^{d+2} - \frac{C' (\gamma + 
	4\beta_1)^{d+4} nh_k^{d+4}}{\varkappa^2} \, .
\end{align}
In the last inequality we used \eqref{e1} and the fact that $\|X_i - X_j\| \leq (\gamma + 
4\beta_1) h_k$, if $\|\widehat X_i\ind k - \widehat X_j\ind K \| \leq \gamma h_k$.

It remains to provide a lower bound for the sum
\[
\min\limits_{u\in \B(0,1) \cap \T_{X_i}\M^*} \sum\limits_{j=1}^n v_{ij} (u^T(X_j - X_i))^2.
\]
Let $\N_\eps$ stand for a $\eps$-net of the set $\B(0,1) \cap \T_{X_i}\M^*$.
It is known that $|\N_\eps| \leq (3/\eps)^d$.
Here and further in this proof we will assume $\eps = \eps_n = 3 n^{-1/d}$.
Then, for any $t > 0$, it holds
\begin{align}
\label{e3}
&\notag
\left\{ \min\limits_{u\in \B(0,1) \cap\T_{X_i}\M^*} \sum\limits_{j=1}^n v_{ij} (u^T(X_j - 
X_i))^2 < t \right\}
\\&
\subseteq \left\{ \min\limits_{u\in \N_\eps} \sum\limits_{j=1}^n v_{ij} (u^T(X_j - X_i))^2 
< 2t + 2\eps^2  \sum\limits_{j=1}^n v_{ij} \|X_i - X_j\|^2  \right\}
\\&\notag
\subseteq \bigcup\limits_{u\in \N_\eps} \left\{ \sum\limits_{j=1}^n v_{ij} (u^T(X_j - 
X_i))^2 < 2t + 4 C' \eps^2 (\gamma + 4\beta_1)^{d+2} nh_k^{d+2} \right\} \, .
\end{align}
Fix any $u\in\N_\eps$ and consider
\[
\p\ind{-i}\left( \sum\limits_{j=1}^n v_{ij} (u^T(X_j - X_i))^2 < 2t + 4 C' \eps^2 (\gamma 
+ 4\beta_1)^{d+2} nh_k^{d+2} \right) \, .
\]
Note that
\begin{align*}
&
\sum\limits_{j=1}^n v_{ij} (u^T(X_j - X_i))^2
\\&
\geq \sum\limits_{j=1}^n \1\left( \|X_i - X_j \| \leq (\gamma - 4\beta_1) h_k \right)  
(u^T(X_j - X_i))^2,
\end{align*}
so it holds
\begin{align*}
&
\p\ind{-i}\left( \sum\limits_{j=1}^n v_{ij} (u^T(X_j - X_i))^2 < 2t + 4 C' \eps^2 (\gamma 
+ 4\beta_1)^{d+2} nh_k^{d+2} \right)
\\&
\leq \p\ind{-i}\bigg( \sum\limits_{j=1}^n \1\left( \|X_i - X_j \| \leq (\gamma - 4\beta_1) 
h_k \right)  (u^T(X_j - X_i))^2
\\&
\qquad\qquad< 2t + 4 C' \eps^2 (\gamma + 4\beta_1)^{d+2} 
nh_k^{d+2} \bigg) \, .
\end{align*}
Given $X_i$, the random variables $\1\left( \|X_i - X_j \| \leq h_k \right)  (u^T(X_j - 
X_i))^2$, $1\leq j\leq n$, are conditionally independent and identically distributed, and 
expectation of each of them can be bounded below by
\begin{align*}
&
\E\ind{-i} \1\left( \|X_i - X_1 \| \leq (\gamma - 4\beta_1) h_k \right)  (u^T(X_1 - X_i))^2
\\&
\geq \frac{p_0}{4} \int\limits_{\M^*\cap\B(X_i, h_k)\cap\{|u^T(X_i-x)|\geq\frac12\}} 
\|x - X_i\|^2 dW(x) \geq c (\gamma - 4\beta_1)^{d+2} h_k^{d+2}.
\end{align*}
At the same time, the variance of these random variables does not exceed
\begin{align*}
&
\E\ind{-i} \1\left( \|X_i - X_1 \| \leq (\gamma - 4\beta_1) h_k \right)  (u^T(X_1 - X_i))^4
\\&
\leq (\gamma - 4\beta_1)^4 h_k^4 \p\ind{-i}\big( \|X_i - X_1 \| \leq (\gamma - 
4\beta_1) h_k \big)
\\&
\leq C(\gamma - 4\beta_1)^{d+4} h_k^{d+4}
\leq C' (\gamma - 4\beta_1)^{d+4} h_k^{d+4},
\end{align*}
where $C' = C \vee 16/3$.
Again, using the Bernstein's inequality, we obtain that for any $\tilde t$ it holds
\begin{align*}
&
\p\ind{-i}\Big( \sum\limits_{j=1}^n \1\left( \|X_i - X_j \| \leq (\gamma - 4\beta_1) 
h_k \right)  (u^T(X_j - X_i))^2
\\&
< n \E\ind{-i} \1\left( \|X_i - X_j \| \leq (\gamma - 
4\beta_1) h_k \right)  (u^T(X_j - X_i))^2 - \tilde t \Big) 
\\&
\leq \exp\left\{-\frac{\tilde t^2}{2n C' (\gamma - 4\beta_1)^{d+4} h_k^{d+4} + 
	2(\gamma - 4\beta_1)^2 h_k^2 \tilde t/3} \right\} \, .
\end{align*}
Take $\tilde t = \delta n\E\ind{-i} \1\left( \|X_i - X_j \| \leq (\gamma - 
4\beta_1) h_k \right)  (u^T(X_j - X_i))^2$.
Then
\begin{align*}
&
\exp\left\{-\frac{\tilde t^2}{2n C' (\gamma - 4\beta_1)^{d+4} h_k^{d+4} + 2(\gamma - 
	4\beta_1)^2 h_k^2 \tilde t/3} \right\}
\\&
\leq \exp\left\{-\frac{c^2 n^2 \delta^2 (\gamma - 4\beta_1)^{2d+4} h_k^{2d+4}}{2n C' 
	(\gamma - 4\beta_1)^{d+4} h_k^{d+4} + \cfrac{2cn\delta}3 (\gamma - 
	4\beta_1)^{d+4} h_k^{d+4}}\right\}
\\&
= \exp\left\{-\frac{n c^2 \delta^2 (\gamma - 4\beta_1)^d h_k^d}{2C' + 
	\cfrac{2c\delta}3}\right\} \, .
\end{align*}
Choose $\delta$ satisfying the inequality
\[
\frac{c^2 \delta^2 (\gamma - 4\beta_1)^d }{2C' + \cfrac{2c\delta}3} \geq 3 \, .
\]
In particular,
\[
\delta = \frac{2}{c (\gamma - 4\beta_1)^d} + \sqrt{\frac{6C'}{c^2(\gamma - 
		4\beta_1)^d}}
\]
is a suitable choice.
Then
\begin{align*}
&
\exp\left\{-\frac{n c^2 \delta^2 (\gamma - 4\beta_1)^d h_k^d}{2C' + 
	\frac{2c\delta}3}\right\}
\leq e^{-3 nh_k^d}
\\&
\leq e^{-3\log n}
\leq e^{-2\log n - \log|\N_\eps|}
= \frac1{|\N_\eps|n^2} \, .
\end{align*}
Thus, with probability at least $1 - (|\N_\eps|n^2)^{-1}$, it holds
\begin{align*}
&
\sum\limits_{j=1}^n \1\left( \|X_i - X_j \| \leq (\gamma - 4\beta_1) 
h_k \right)  (u^T(X_j - X_i))^2
\\&
\geq \left(1 - \frac{2}{c (\gamma - 4\beta_1)^d} - \sqrt{\frac{6C'}{c^2(\gamma - 
		4\beta_1)^d}} \right) c (\gamma - 4\beta_1)^{d+2} h_k^{d+2} \, .
\end{align*}
By the union bound, on an event with probability at least $1 - n^{-2}$ it holds
\begin{align}
\label{e4}
&
\min\limits_{u \in \N_\eps} \sum\limits_{j=1}^n \1\left( \|X_i - X_j \| \leq (\gamma - 
4\beta_1) h_k \right)  (u^T(X_j - X_i))^2
\\&\notag
\geq \left(1 - \frac{2}{c (\gamma - 4\beta_1)^d} - \sqrt{\frac{6C'}{c^2(\gamma - 
		4\beta_1)^d}} \right) cn (\gamma - 4\beta_1)^{d+2} h_k^{d+2}
\end{align}
Then, due to \eqref{e3} and \eqref{e4}, on this event
\begin{align*}
&
\min\limits_{u\in \B(0,1) \cap \T_{X_i}\M^*} \sum\limits_{j=1}^n v_{ij} (u^T(X_j - X_i))^2
\\&
\geq \frac c2 \left(1 - \frac{2}{c (\gamma - 4\beta_1)^d} - 
\sqrt{\frac{6C'}{c^2(\gamma - 4\beta_1)^d}} \right) (\gamma - 4\beta_1)^{d+2} 
n h_k^{d+2}
\\&
- 2 C' \eps^2 (\gamma + 4\beta_1)^{d+2} nh_k^{d+2},
\end{align*}
and, together with \eqref{e2}, this yields
\begin{align*}
&
\lambda_d(\widehat\bxi_i) - \lambda_{d+1}(\widehat\bxi_i)
\\&
\geq \frac c4 \left(1 - \frac{2}{c (\gamma - 4\beta_1)^d} - 
\sqrt{\frac{6C'}{c^2(\gamma - 4\beta_1)^d}} \right) (\gamma - 4\beta_1)^{d+2} 
n h_k^{d+2}
\\&
- C' \eps^2 (\gamma + 4\beta_1)^{d+2} nh_k^{d+2}
- 16 C'\beta_1^2 (\gamma + 4\beta_1)^d n h_k^{d+2}
\\&
- \frac{C' (\gamma + 4\beta_1)^{d+4} nh_k^{d+4}}{\varkappa^2}.
\end{align*}
The choice $\eps = 3n^{-1/d}$ yields the claim of Lemma \ref{lem_spectral_gap}.
\end{proof}

\begin{Lem}
	\label{p1a1}
	Let $\E\ind{-i}$ denote the conditional expectation $\E(\cdot | (X_i, Y_i))$ and let 
	$r_d = 4h \sqrt{(d+2) \log h^{-1}}$.
	Then, for any $i$ from $1$ to $n$, it holds
	\[
		\E\ind{-i} e^{-\frac{\|X_j - X_i\|^2}{h^2}} u^T(X_j - X_i) \1\left( X_j \in \B(X_i, r_d) 
		\right) \lesssim \frac{d h^{d+2}}\varkappa.
	\]
\end{Lem}

\begin{proof}
	
We have
\begin{align*}
	&
	\E\ind{-i} e^{-\frac{\|X_j - X_i\|^2}{h^2}} u^T(X_j - X_i) \1\left( X_j \in \B(X_i, r_d) \right) 
	\\&
	= \int\limits_{\M^* \cap \B(X_i, r_d)} e^{-\frac{\|x- X_i\|^2}{h^2}} u^T(x - X_i) p(x) 
	dW(x).
\end{align*}
Due to \eqref{a1'}, the last expression does not exceed 
\begin{align*}
	&
	\leq p(X_i) \int\limits_{\M^* \cap \B(X_i, r_d)} e^{-\frac{\|x- X_i\|^2}{h^2}} u^T(x - X_i) 
	dW(x)
	\\&
	+ \frac{L}\varkappa \int\limits_{\M^* \cap \B(X_i, r_d)} e^{-\frac{\|x- 
	X_i\|^2}{h^2}} \|x - X_i\|^2 dW(x).
\end{align*}
Due to Lemma \ref{p1a4}, 
\[
	\frac{L}\varkappa \int\limits_{\M^* \cap \B(X_i, r_d)} e^{-\frac{\|x- X_i\|^2}{h^2}} \|x - 
	X_i\|^2 dW(x)
	\lesssim \frac{h^{d+2}}\varkappa,
\]
so it remains to prove that
\[
	p(X_i) \int\limits_{\M^* \cap \B(X_i, r_d)} e^{-\frac{\|x- X_i\|^2}{h^2}} u^T(x - X_i) 
	dW(x) \lesssim \frac{h^{d+2}}\varkappa.
\]

Let $\Exp_{X_i}(\cdot)$ be the exponential map of $\M^*$ at $X_i$ and denote 
$\widetilde\B(X_i, r_d) = \Exp^{-1}(\M^*\cap \B(X_i, r_d))$.
Note that $\Exp_{X_i}(\cdot)$ is a bijection on $\widetilde\B(X_i, r_d)$ (see, for instance, \citep[Lemma 1]{al19}), because $r_d \leq \varkappa/4$.
Then
\begin{align*}
	&
	\int\limits_{\M^* \cap \B(X_i, r_d)} e^{-\frac{\|x- X_i\|^2}{h^2}} u^T(x - X_i) 
	dW(x)
	\\&
	= \int\limits_{\widetilde\B(X_i, r_d)} e^{-\frac{\|\Exp_{X_i}(v) - \Exp_{X_i}(0)\|^2}{h^2}} 
	u^T(\Exp_{X_i}(v) - \Exp_{X_i}(0)) \sqrt{\det g(v)} dv.
\end{align*}
Introduce functions
\[
	\psi_{X_i}(v) = \Exp_{X_i}(v) - \Exp_{X_i}(0) - v
\]
and
\[
	\varphi_{X_i}(v) = \|\psi_{X_i}(v)\|^2 + 2 v^T\psi_{X_i}(v).
\]
Due to \citep[Lemma 1]{al19}, it holds that
\begin{equation}
	\label{d2}
	\| \Exp_{X_i}(v) - \Exp_{X_i}(0) - v\|  = \| \psi_{X_i}(v) \| \leq \frac{5\|p\|^2}{4\varkappa},
\end{equation}
which yields $\psi_{X_i}(v) = O\left( \frac{\|v\|^2}\varkappa\right)$, 
$\varphi_{X_i}(v) = O\left( \frac{\|v\|^3}\varkappa\right)$.
Now, consider $\sqrt{\det g(v)}$.
It is known (see, for instance, \citep[Equation 2.1]{tsay19}) that there exists an absolute 
constant $\overline C$, such that
\begin{equation}
	\label{d2''}
	\left| \sqrt{\det g(v)} - 1 \right| \leq \frac{\overline C d \|v\|^2}{\varkappa^2}.
\end{equation}

Taking \eqref{d2} and \eqref{d2''} into account, we obtain
\begin{align*}
	&\notag
	\int\limits_{\widetilde\B(X_i, r_d)} e^{-\frac{\|\Exp_{X_i}(v) - \Exp_{X_i}(0)\|^2}{h^2}} 
	u^T(\Exp_{X_i}(v) - \Exp_{X_i}(0)) \sqrt{\det g(v)} dv
	\\&
	= \int\limits_{\widetilde\B(X_i, r_d)} e^{-\frac{\|v\|^2 + \varphi_{X_i}(v)}{h^2}} 
	u^T(v + \psi_{X_i}(v)) \sqrt{\det g(v)} dv
	\\&\notag
	= \int\limits_{\widetilde\B(X_i, r_d)} e^{-\frac{\|v\|^2}{h^2}}\left(1 + 
	O\left(\frac{\|v\|^3}{h^2 \varkappa}\right)\right) \left(u^Tv + 
	O\left(\frac{\|v\|^2}\varkappa\right)\right)\left(1 
	+ O\left(\frac{\|v\|^2}{\varkappa^2}\right)\right)dv
	\\&\notag
	= \int\limits_{\widetilde\B(X_i, r_d)} e^{-\frac{\|v\|^2}{h^2}} u^Tv dv
	+ \int\limits_{\widetilde\B(X_i, r_d)} e^{-\frac{\|v\|^2}{h^2}} O\left(\frac{\|v\|^4}{h^2 
	\varkappa}\right) dv
	\\&
	+ \int\limits_{\widetilde\B(X_i, r_d)} e^{-\frac{\|v\|^2}{h^2}} O\left(\frac{\|v\|^2}{ 
	\varkappa}\right)dv
	+ \int\limits_{\widetilde\B(X_i, r_d)} e^{-\frac{\|v\|^2}{h^2}} O\left(\frac{\|v\|^3}{ 
	\varkappa^2}\right)dv.
\end{align*}
For the last three terms, we get
\begin{align*}
	&
	\int\limits_{\widetilde\B(X_i, r_d)} e^{-\frac{\|v\|^2}{h^2}} O\left(\frac{\|v\|^4}{h^2 
	\varkappa}\right) dv \lesssim \frac{h^{d+2}}\varkappa,
	\\&
	\int\limits_{\widetilde\B(X_i, r_d)} e^{-\frac{\|v\|^2}{h^2}} O\left(\frac{\|v\|^2}{ 
	\varkappa}\right)dv \lesssim \frac{h^{d+2}}\varkappa,
	\\&
	\int\limits_{\widetilde\B(X_i, r_d)} e^{-\frac{\|v\|^2}{h^2}} O\left(\frac{\|v\|^3}{ 
	\varkappa^2}\right)dv \lesssim \frac{h^{d+3}}{\varkappa^2} \lesssim 
	\frac{h^{d+2}}\varkappa.
\end{align*}
Thus,
\begin{align}
	&
	\int\limits_{\widetilde\B(X_i, r_d)} e^{-\frac{\|\Exp_{X_i}(v) - \Exp_{X_i}(0)\|^2}{h^2}} 
	u^T(\Exp_{X_i}(v) - \Exp_{X_i}(0)) \sqrt{\det g(v)} dv
	\\&\notag
	= \int\limits_{\widetilde\B(X_i, r_d)} e^{-\frac{\|v\|^2}{h^2}} u^Tv dv + 
	O\left(\frac{h^{d+2}}\varkappa\right),
\end{align}
and, in order to complete the proof, we have to show that
\[
	\int\limits_{\widetilde\B(X_i, r_d)} e^{-\frac{\|v\|^2}{h^2}} u^Tv dv = 
	O\left(\frac{h^{d+2}}\varkappa\right).
\]

Note that, for any fixed $u \in \R^d$, it holds
\[
	\int\limits_{\R^d} e^{-\frac{\|v\|^2}{h^2}} u^Tv dv = 0.
\]
Then
\[
	\int\limits_{\widetilde\B(X_i, r_d)} e^{-\frac{\|v\|^2}{h^2}} u^Tv dv
	= -\int\limits_{\R^d\backslash\widetilde\B(X_i, r_d)} e^{-\frac{\|v\|^2}{h^2}} u^Tv dv
\]
Remind that $\widetilde\B(X_i, r_d) = \Exp_{X_i}^{-1}(\B(X_i, r_d)) = \{v : \|\Exp_{X_i}(v) - 
\Exp_{X_i}(0)\| \leq r_d\}$.
By the definition of the exponential map,
\[
	\|\Exp_{X_i}(v) - \Exp_{X_i}(0)\| \leq d_{\M^*}(\Exp_{X_i}(v), \Exp_{X_i}(0)) = \|v\|.
\]
Then we conclude that $\widetilde\B(X_i, r_d) \supseteq \B(0, r_d)$.
This yields
\begin{align*}
	&
	\left|\int\limits_{\R^d\backslash\widetilde\B(X_i, r_d)} e^{-\frac{\|v\|^2}{h^2}} u^Tv dv 
	\right|
	\leq \int\limits_{\R^d\backslash\widetilde\B(X_i, r_d)} e^{-\frac{\|v\|^2}{h^2}} \|v\| dv
	\\&
	\leq \int\limits_{\R^d\backslash\B(0, r_d)} e^{-\frac{\|v\|^2}{h^2}} \|v\| dv
	\leq e^{-\frac{r_d^2}{2h^2}}\int\limits_{\R^d\backslash\B(0, r_d)} 
	e^{-\frac{\|v\|^2}{2h^2}} \|v\| dv
	\\&
	\leq e^{-\frac{r_d^2}{2h^2}}\int\limits_{\R^d} e^{-\frac{\|v\|^2}{2h^2}} \|v\| dv.
\end{align*}
By definition, $r_d = 2h \sqrt{2(d+2) \log h^{-1}}$.
This implies $e^{-\frac{r_d^2}{2h^2}} = h^{4(d+2)}$.
Moreover, 
\[
	\int\limits_{\R^d} e^{-\frac{\|v\|^2}{2h^2}} \|v\| dv \lesssim h^{d+1}.
\]
Thus, we conclude
\begin{equation}
	\label{d3}
	\int\limits_{\widetilde\B(X_i, r_d)} e^{-\frac{\|v\|^2}{h^2}} u^Tv dv \lesssim h^{d+1 + 
	4(d+2)} \lesssim \frac{h^{d+2}}\varkappa,
\end{equation}
and \eqref{d3} finishes the proof of Lemma \ref{p1a1}.
\end{proof}

\begin{Lem}
	\label{p1a4}
	Let $\E\ind{-i}$ denote the conditional expectation $\E(\cdot | (X_i, Y_i))$ and let 
	$r_d = 4h \sqrt{(d+2) \log h^{-1}}$.
	Then, for any $i$ from $1$ to $n$, it holds
	\[
		\E\ind{-i} e^{-\frac{\|X_i - X_j\|^2}{16h^2}} \|X_j - X_i\|^q \1\left( X_j \in \B(X_i, r_d) 
		\right) \lesssim h^{q+d}.
	\]
\end{Lem}

\begin{proof}

Using \eqref{a1'}, we obtain
\begin{align*}
	&
	\E\ind{-i} e^{-\frac{\|X_i - X_j\|^2}{16h^2}} \|X_j - X_i\|^q \1\left( X_j \in \B(X_i, r_d) 
	\right)
	\\&
	= \int\limits_{\M^*\cap \B(X_i, r_d)}  e^{-\frac{\|X_i - x\|^2}{16h^2}} \|x - X_i\|^q  p(x) 
	dW(x)
	\\&
	\leq p_1 \int\limits_{\M^*\cap \B(X_i, r_d)}  e^{-\frac{\|X_i - x\|^2}{16h^2}} \|x - X_i\|^q 
	dW(x).
\end{align*}
Using the exponential map, we get
\begin{align*}
	&
	\int\limits_{\M^*\cap \B(X_i, r_d)}  e^{-\frac{\|X_i - x\|^2}{16h^2}} \|x - X_i\|^q 
	dW(x)
	\\&
	= \int\limits_{\widetilde \B(X_i, r_d)}  e^{-\frac{\|\Exp_{X_i}(v) - 
	\Exp_{X_i}(0)\|^2}{16h^2}} \|\Exp_{X_i}(v) - \Exp_{X_i}(0)\|^q \sqrt{\det g(v)} dv
\end{align*}
Taking into account that $\|v\| = d_{\M^*}(\Exp_{X_i}(v), \Exp_{X_i}(0))$ and applying \citep[Lemma 2.5]{blw18}, we conclude
\[
	\frac{\|v\|}2
	\leq \|\Exp_{X_i}(v) - \Exp_{X_i}(0)\|
	\leq \|v\|.
\]
On the other hand, \eqref{d2''} yields $\sqrt{\det g(v)} \lesssim 1$ for all $v \in 
\widetilde\B(X_i, r_d)$.
Thus, we obtain
\begin{align*}
	&
	\int\limits_{\widetilde \B(X_i, r_d)}  e^{-\frac{\|\Exp_{X_i}(v) - 
	\Exp_{X_i}(0)\|^2}{16h^2}} \|\Exp_{X_i}(v) - \Exp_{X_i}(0)\|^q \sqrt{\det g(v)} dv
	\\&
	\lesssim \int\limits_{\widetilde \B(X_i, r_d)} e^{-\frac{\|v\|^2}{32h^2}} \|v\|^q dv
	\leq \int\limits_{\R^d} e^{-\frac{\|v\|^2}{32h^2}} \|v\|^q dv
	\lesssim h^{d+q}.
\end{align*}

\end{proof}

\section{Pseudocode of the Manifold Blurring Mean Shift Algorithm}
\label{sec_mbms}

This section contains a pseudocode of the manifold blurring mean shift algorithm 
\citep[MBMS]{wcp10}.

\newpage

\begin{algorithm}[t]
	\caption{Manifold blurring mean shift algorithm (with full graph), \cite{wcp10}}
	\label{mbms}
	\begin{algorithmic}[1]
		\State The sample of noisy observations $\Y_n = (Y_1, \dots, Y_n)$, a 
		bandwidth $\sigma > 0$, and positive integers $\mathsf{k}$ and $d$ are given.
		\State Initialize $\widehat X_1 = Y_1, \dots, \widehat X_n = Y_n$.
		\Repeat
		\State Compute the increments
		\[
			\partial \widehat X_i
			= - \widehat X_i + \frac{\sum\limits_{j = 1}^n \K\left(\frac{\|\widehat X_j - 
			\widehat X_i\|^2}{2\sigma^2}\right) \widehat X_j}{\sum\limits_{j = 1}^n 
			\K\left(\frac{\|\widehat X_j - \widehat X_i\|^2}{2\sigma^2}\right)},
			\quad 1 \leq i \leq n,
		\]
		where $\K(t) = e^{-t}$.
		\State For each $i$ from $1$ to $n$, find $\mathsf{k}$ nearest neighbors $\mathcal N_i$ 
		of $\widehat X_i$.
		\State For all $i$ from $1$ to $n$, perform local PCA, that is compute
		\begin{align*}
			\mu_i
			&
			= \frac1{\mathsf{k}} \sum\limits_{j \in \mathcal N_i} \widehat X_j, \quad 1 \leq i \leq n,
			\\
			\bsigma_i
			&
			= \frac1{\mathsf{k}} \sum\limits_{j \in \mathcal N_i} (\widehat X_j - \mu_i)(\widehat X_j 
			- \mu_i)^T, \quad 1 \leq i \leq n,
		\end{align*}
		and put $\bpi_i$ a projector onto a linear span of eigenvectors of 
		$\widehat{\bsigma}_{i}$, corresponding to the largest $d$ eigenvalues.
		\State Update the increments
		\[
			\partial \widehat X_i \gets (\bid - \bpi_i) \partial \widehat X_i, \quad 1 \leq i 
			\leq n.
		\]
		\State Update the estimates
		\[
			\widehat X_i \gets \widehat X_i + \partial X_i, \quad 1 \leq i \leq n.
		\]
		\Until stop\\
		\Return the estimates $\widehat{X}_1, \dots, \widehat{X}_n$.
	\end{algorithmic}
\end{algorithm}

\vskip 0.2in
\bibliography{references.bib}

\end{document}